\pdfoutput=1                                     
\documentclass[twocolumn,secthm]{autart}

\usepackage{amssymb,amsmath,amsfonts}  % assumes amsmath package installed
\usepackage[pdftex]{graphicx}
\usepackage{epstopdf}

\usepackage{mathtools} 
\usepackage[ruled,norelsize,algo2e]{algorithm2e}
\usepackage{comment}
\usepackage{bm}

\usepackage{tikz}
\usetikzlibrary{positioning}
\usepackage[normalem]{ulem}

\usepackage{setspace}
\usepackage{lineno}
\usepackage{cite}
\usepackage{url}
\usepackage{booktabs}
\usepackage{multirow}
\usepackage{dsfont}
\usepackage{float} % for placing figures
\usepackage{soul}
\usepackage{schemata}
\usepackage{stfloats}
\usepackage{enumitem}

\usepackage{caption}% http://ctan.org/pkg/caption
\allowdisplaybreaks

\newcommand{\defeq}{:=}
\newcommand{\dd}{\text{d}}

\newcommand{\innerproduct}[2]{\langle#1,#2\rangle}
\newcommand{\norm}[1]{\left\lVert#1\right\rVert}
\newcommand{\abs}[1]{\left|#1\right|}
\newcommand{\costeval}[2]{J_{#1}(#2)}
\newcommand{\optcosteval}[2]{J^*_{#1}(#2)}
\DeclareMathOperator*{\esssup}{ess\,sup}
\newcommand{\trace}[1]{\text{Tr}(#1)}

\newcommand{\commentout}[1]{}

\newcommand{\revision}[1]{{\color{black} #1}}

\begin{document}
% \author{Sheng Cheng and Derek A. Paley
% % \thanks{Work reported in this article was sponsored under DARPA Contract No. HR001118C0142.}% <-this % stops a space
% \thanks{Sheng Cheng is a graduate student in the Department of Electrical and Computer Engineering and the Institute for Systems Research at the University of Maryland, College Park, MD 20742, USA. (email:{\tt\small \{cheng@terpmail.umd.edu\}})}
% \thanks{Derek A. Paley is the Willis H. Young Jr. Professor in the Department of Aerospace Engineering and the Institute for Systems Research at the University of Maryland, College Park, MD 20742, USA. (email:{\tt\small \{dpaley@umd.edu\}})}
% %\thanks{Artur Wolek and Derek A. Paley are with the Department of Aerospace Engineering and the Institute for Systems Research at the University of Maryland, College Park, MD 20742, USA. (email:{\tt\small \{wolek@umd.edu\} and \{dpaley@umd.edu\}})}
% }

\begin{frontmatter}
\date{15 June 2021}
\title{Optimal control of a 2D diffusion-advection process with a team of mobile actuators under jointly optimal guidance\thanksref{footnoteinfo}} % Title, preferably not more 

\thanks[footnoteinfo]{This paper was not presented at any IFAC 
meeting. Corresponding author Sheng Cheng (Tel. +1 301 335 2995).}

\author[Sheng]{Sheng Cheng}\ead{cheng@terpmail.umd.edu},    % Add the 
\author[Derek]{Derek A. Paley}\ead{dpaley@umd.edu}               % e-mail address 
% \author[Baiae]{Publius Maro Vergilius}\ead{vergilius@culture.ir}  % (ead) as shown

\address[Sheng]{University of Maryland, College Park}  % Please supply                                              
\address[Derek]{University of Maryland, College Park}             % full addresses
% \address[Baiae]{The White House, Baiae}        % here.

                             % keyword list or with the 
                                          % help of the Automatica 
                                          % keyword wizard

\begin{abstract}                          % Abstract of not more than 200 words.
This paper describes an optimization framework to control a distributed parameter system (DPS) using a team of mobile actuators. The framework simultaneously seeks optimal control of the DPS and optimal guidance of the mobile actuators such that a cost function associated with both the DPS and the mobile actuators is minimized subject to the dynamics of each. The cost incurred from controlling the DPS is linear-quadratic, which is transformed into an equivalent form as a quadratic term associated with an operator-valued Riccati equation. This equivalent form reduces the problem to seeking for guidance only because the optimal control can be recovered once the optimal guidance is obtained. We establish conditions for the existence of a solution to the proposed problem. Since computing an optimal solution requires approximation, we also establish the conditions for convergence to the exact optimal solution of the approximate optimal solution. \revision{That is, when evaluating these two solutions by the original cost function, the difference becomes arbitrarily small as the approximation gets finer.} Two numerical examples demonstrate the performance of the optimal control and guidance obtained from the proposed approach.

\end{abstract}

\end{frontmatter}

\normalsize  

\section{Introduction}
Recent development of mobile robots (unmanned aerial vehicles, terrestrial robots, and underwater vehicles) has greatly extended the type of distributed parameter system (DPS) over which mobile actuation and sensing can be deployed. Such a system is often modeled by a partial differential equation (PDE), which varies in both time and space. Exemplary applications of mobile control and estimation of a DPS can be found in thermal manufacturing \cite{demetriou2003scanning},
%displacement estimation of a flexible structure \cite{demetriou2014guidance}, 
monitoring and neutralizing groundwater contamination \cite{demetriou2009estimation},
and wildfire monitoring~\cite{kumar2011cooperative}.

We propose an optimization framework that simultaneously solves for the guidance of a team of mobile actuators and the control of a DPS. We consider a 2D diffusion-advection process as the DPS for its capability of modeling a variety of processes governed by continuum mechanics and the convenience of the state-space representation. The framework minimizes an integrated cost function, evaluating both the control of the DPS and the guidance of the actuators, subject to the dynamics of the DPS and the mobile actuators. The problem addresses the mobile actuator and the DPS as a unified system, instead of solely controlling the DPS. Furthermore, the additional degree of freedom endowed by mobility yields improved control performance in comparison to using stationary actuators.

\revision{The study of the control of a PDE-modeled DPS dates back to the 1960s (see the survey \cite{robinson1971survey}). For fundamental results, see the textbooks \cite{curtain2012introduction,omatu1989distributed,bensoussan1992representation}. Although it is possible to categorize prior work by whether the input operator is bounded, our literature review is categorized by the location of the actuation.}
When actuation acts on the boundary of the spatial domain, it is called boundary control. \revision{The main complexity in boundary control is that the input operator, which associates with the PDE's boundary condition, is unbounded. This is addressed in the tutorial \cite{emirsjlow2000pdes}. Recent developments on the design of boundary control uses the backstepping method} \cite{smyshlyaev2010adaptive}, where a Volterra transformation determines a stabilizing control by transforming the system into a stable target system. 
When actuation acts in the interior of the spatial domain, it is called distributed control. For distributed control, the DPS is actuated by in-domain actuators that are either stationary or mobile. 

The problem of determining the location of stationary actuators is called the actuator placement problem. Actuator placement has been studied for optimality in the sense of linear-quadratic (LQ) \cite{morris2010linear}, H$_2$~\cite{morris2015using}, and H$_{\infty}$~\cite{kasinathan2013Hinfinity}.
The author of \cite{morris2010linear} studies the actuator placement problem with the LQ performance criterion. The actuators' locations are chosen to minimize the operator norm or trace norm of the Riccati operator solved from an algebraic Riccati equation associated with an infinite-dimensional system. If the input operator is compact and continuous with respect to actuator location in the operator norm, then a solution exists for the problem minimizing the operator norm of the Riccati operator \cite[Theorem~2.6]{morris2010linear}, under stabilizability and detectability assumptions. When computing the optimal actuator locations, if the approximated Riccati operator converges to the original Riccati operator at each actuator location, then the approximate optimal actuator locations converge to the exact optimal locations \cite[Theorem~3.5]{morris2010linear}. For the above results to hold when minimizing the trace-norm of the Riccati operator, the input and output spaces have to be finite-dimensional \cite[Theorems~2.10 and~3.9]{morris2010linear} in addition to the assumptions stated above.

The authors of \cite{morris2015using} design optimal actuator placement by minimizing the H$_2$-control performance criterion, which minimizes the L$_2$-norm of the linear output of a linear system, subject state disturbances. Roughly speaking, H$_2$-control performance reduces the response to the disturbances while setting a zero initial condition, whereas the LQ performance reduces the response to the initial condition in a disturbance-free setting. For disturbances with known or unknown spatial distribution, the trace of the Riccati solution (scaled by the disturbance's spatial distribution) or operator norm of the Riccati solution are minimized, respectively, where the existence of a solution and convergence to the exact optimal solution of the approximate optimal solution are guaranteed.
In \cite{kasinathan2013Hinfinity}, the H$_{\infty}$-performance criterion is minimized for actuator placement. Specifically, the actuators are placed in the locations that yield infimum of the attenuation bound (upper bound of the infinity norm of the closed-loop disturbance-to-output transfer function). The conditions for the existence of a solution and convergence to the exact optimal placement of the approximate optimal placement are provided.

Geometric approaches have also been proposed for actuator placement. For example, a modified centroidal Voronoi tessellation (mCVT) yields locations of actuators and sensors that yields least-squares approximate control and estimation kernels for a parabolic PDE \cite{demetriou2017using}. The input operator is designed by maximizing the H$_2$-norm of the input-to-state transfer function, whereas the kernel of the state feedback is obtained using the Riccati solution. Next, mCVT determines the partition such that the actuator and sensor locations achieve optimal performance (in the sense of least-squares) to the input operator and state feedback kernel, respectively. %A dual problem for sensor placement is handled using mCVT in \cite{demetriou2018integrated}.
A comparison of various performance criteria for actuator placement has been conducted for controlling a simply supported beam \cite{morris2015comparison} and a diffusion process \cite{morris2016study}. It has been analyzed that maximizing the minimum eigenvalue of the controllability gramian is not a useful criterion. Because the lower bound of the eigenvalue is zero, the minimum eigenvalue approaches zero as the dimension of approximation increases \cite{morris2015comparison,morris2016study}.

The guidance of mobile actuators is designed to improve the control performance in comparison to stationary actuators. Various performance criteria have been proposed for guidance.
In \cite{demetriou2003scanning}, a mobile heat source is steered to maintain a spatially uniform temperature distribution in 1D using the optimal control method. The formulation uses a finite-dimensional approximation for modeling the process and evaluating performance. Additionally, the admissible locations of the heat source are chosen within a discrete set that yields approximate controllability requirements. Algorithms are provided to solve the proposed problem with considerations on real-time implementation and hardware constraints. Experimental results demonstrate the performance of the proposed scheme. %Note that the actuator considered in this paper, although mobile, is different in terms of mobility compared to the ones in our paper: the admissible locations are discretized into a finite set while we consider continuous trajectories. Also, the problem here is not applicable to the problem where mobility is constrained by onboard resources.
The authors of \cite{dubljevic2010discrete} propose an optimization framework that steers mobile actuators to control a reaction-diffusion process. A specific cost function, consisting of norms of control input and measurement of the DPS and the steering force, is minimized subject to dynamics of the actuator's motion and the PDE, and bounds on the control input and state of the DPS. The implementation of the framework is emphasized by discrete mechanics and model predictive control which yield computationally tractable solutions, in addition to an approximation of the PDE and a discrete set of admissible actuator locations.

The problem of ultraviolet curing using a mobile radiant actuator is investigated in \cite{zeng2010estimation}, where the curing process is modeled by a 1D nonlinear PDE. Both the radiant power and scanning velocity of the actuator are computed for reaching a target curing state. A dual extended Kalman filter is applied to estimate the state and parameters of the curing process for feedback control, based on the phases of curing.
In \cite{demetriou2020navigating}, a navigation problem is studied in which a mobile agent moves through a diffusion process represented by a hazardous field with given initial and terminal positions. Such a formulation may be applied to emergency evacuation guidance from an indoor environment with carbon monoxide. Both problems with minimum time and minimum accumulated effects of hazards are formulated, and closed-form solutions are derived using the Hamiltonian. 
A Lyapunov-based guidance strategy for collocated sensors and actuators to estimate and control a diffusion process is proposed in \cite{demetriou2010guidance}. The decentralized guidance of the actuators for controlling a diffusion process to a zero state is derived by constructing suitable Lyapunov functions. The same methodology is applied to construct a distributed consensus filter via the network among agents to improve state estimation. A follow-up work \cite{demetriou2012adaptive} incorporates nonlinear vehicle dynamics in addition.

% Geometric approaches have also been applied to actuator guidance. 
% An actuator guidance strategy is proposed in \cite{chen2005actuation,chen2006optimal}, where the centroidal Voronoi tessellation (CVT) partitions a parabolic PDE modeled concentration field. Obstacle avoidance for the mobile actuators in complex environments is incorporated using CVT with a repulsive potential field.

% {\color{red}
% \begin{itemize}
%     \item What's been done? What is relation to prior work? 
%     \item What novelty is in the proposed work? 
%     \item What's the significance of the proposed work? 
% \end{itemize}
% }

The problem formulation in this paper includes a cost function that simultaneously evaluates controlling the PDE-modeled DPS, referred to as the \textit{PDE cost}, and steering the mobile actuators, referred to as the \textit{mobility cost}.
The PDE cost is a quadratic cost of the PDE state and control, whose optimal value can be obtained by solving an operator-valued differential Riccati equation. Our results are based on the related work \cite{Burns2015Solutions}, which establishes Bochner integrable solutions of finite-horizon Riccati integral equations (with values in Schatten p-classes) associated with infinite-dimensional systems. The existence conditions for the solution of exact and approximate Riccati integral equations are established in \cite{Burns2015Solutions}. The significance of the Bochner integrable solution is that it allows the implementation of simple numerical quadratures for computing the approximated solution of Riccati integral equations. In \cite{Burns2015Solutions}, the Riccati solution is applied in a sensor placement problem, which computes optimal sensor locations that minimize the trace of the covariance operator of the Kalman filter of a diffusion-advection process. The same cost has been used in an optimization framework for mobile sensor's motion planning in \cite{burns2015infinite}. The existence of a solution of the optimization problem is established under the assumption that the integral kernel of the output operator is continuous with respect to the location of the sensor \cite[Definition~4.5]{burns2015infinite}. This assumption permits the construction of a compact set of operators \cite[Lemma 4.6]{burns2015infinite} over which the cost function is continuous, and hence establishes the existence of a solution. The assumption is also made on the input operator in this paper, which allows the derivation of a vital result on the Riccati operator's continuity with respect to the actuator trajectory (Lemma~\ref{lemma: continuity of equivalent cost wrt actuator location}). The continuity property plays a crucial role in establishing the existence of the proposed problem's solution and the convergence to the exact optimal solution of the approximate optimal solution.
The existence of a solution is established in using the fact that a weakly sequentially lower semicontinuous function attains its minimum on a weakly sequentially compact set over a normed linear space%\cite[Theorem~6.1.4]{werner2013optimization}
. In addition to the assumptions made for the existence of a solution, a stringent (yet with reasonable physical interpretation) requirement is placed on the admissible set to yield compactness, which leads to convergence of the approximate optimal solution. \revision{The convergence is in the sense that when evaluating the exact and approximate optimal solutions by the original cost function, the difference becomes arbitrarily small as the dimension of approximation increases.} %{\color{red} Think about the following sentence. I feel it shouldn't be here.} Although mobility cost has been formulated as a quadratic function of the steering force in \cite{dubljevic2010discrete}, our formulation evaluates steering using generic functions, which allows for the incorporation of various types of cost or penalty functions.

% \setlength{\tabcolsep}{5pt} % Default value: 6pt
% \renewcommand{\arraystretch}{1} % Default value: 1
%   \captionsetup{%size=footnotesize,
% 	%justification=centering, %% not needed
% 	skip=5pt, position = bottom}
% \begin{table*}
% 	\centering
% 	\small

% 	%\captionsetup{font=small}
% 	\caption{Cost breakdown of control and guidance in comparison. All costs are normalized with respect to the total cost of the naive control and guidance.}
% 	\begin{tabular}{lcccc}
% 		\toprule[1pt]
% 		 & assumption & \multirow{2}{*}{continuity of the PDE cost} & \multirow{2}{*}{Existence of a solution} & \multirow{2}{*}{Convergence of approximate solution} \\
% 		 & continuity of the input operator & & & \\
% 		 		 \midrule
% 		 \cite{morris2010linear} & wrt actuator location \cite[Theorem 2.6]{morris2010linear} & wrt acutator location \cite[Theorem 2.6]{morris2010linear} & \cite[Theorem 2.6]{morris2010linear} & \cite[Theorem 3.5]{morris2010linear} \\
% 		 Proposed work & wrt actuator trajectory: assumption (A4) & wrt actuator trajectory Lemmas~\ref{lemma: continuity of the total cost wrt guidance} and \ref{lemma: continuity of equivalent cost wrt actuator location} & Theorem~\ref{thm: existence of a solution of IOCA} & Theorem~\ref{thm: convergence of approximate solution}\\
% 		\bottomrule[1pt]
% 	\end{tabular}
% \end{table*}
% \normalsize

The contributions of this paper are threefold. 
First, we propose an optimization framework for controlling a PDE-modeled system using a team of mobile actuators. The framework incorporates both controlling the process and steering the mobile actuators. The former is handled by the linear-quadratic regulator of the PDE. The latter is taken care of by designing generic cost functions that address the constraints and limitations of the vehicles carrying the actuators.
Second, existence conditions of a solution of the proposed problem are established. It turns out that the conditions are generally satisfied in engineering problems, which allows the results to be applied to a wide range of applications. 
Third, conditions are also established under which the optimal solution computed using approximations converges to the exact optimal solution. \revision{The convergence is in the sense that the cost function of the exact problem evaluated at these two solutions becomes arbitrarily close as the dimensional of the approximation goes to infinity.} The convergence \revision{is verified in numerical studies and} confirms the appropriateness of the optimal solution of the approximation.

The proposed framework is well-suited for the limited onboard resources of mobile actuators in the following two aspects: (1) it adopts a finite-horizon optimization scheme that characterizes the resource limitation more precisely than the alternative approaches that do not specify a terminal time, such as an infinite-horizon optimization or Lyapunov-based method; and (2) it provides an intermediate step for the optimization problem that characterizes the limited resources as inequality constraints, because the constraints can be used to augment the cost function and turned into the proposed form using the method of Lagrange multipliers.
Potential applications of this work include forest firefighting using unmanned aerial vehicles and oil spill removal or harmful algae containment using autonomous skimmer boats. %Also, the simple formulation allows it to be combined with the model predictive control for implementation. 
A preliminary version of this paper \cite{cheng2020optimalACC} considered controlling a 1D diffusion process by a team of mobile actuators. The results in this paper extends the controlled process in \cite{cheng2020optimalACC} to a 2D diffusion-advection process and generalizes the mobility cost therein. Furthermore, the proofs of the existence of a solution and convergence of the approximated optimal solution are presented for the first time in this paper. The results for a dual estimation framework can be found in \cite{Cheng2021estimationInReivew}.

% The contributions of this paper are summarized as follows: (1) the problem of regulating a diffusion process with a \textit{team} of mobile actuators is formulated as an integrated optimization problem; (2) conditions are showed under which optimal guidance is guaranteed to restrict the mobile actuators to the spatial domain even without explicit constraints; % OR:   (2) conditions are proved under which any guidance that steers the mobile actuators out of the spatial domain will be non-optimal;
% %the mobile actuators will not wander out of the spatial domain of the diffusion process even without explicit constraints;
% and (3) a gradient-descent method is applied to numerically solve the approximation of the formulated problem for the jointly optimal guidance of the mobile actuators and the (feedback) control of the DPS. The proposed framework provides a new approach to simultaneously design the guidance and actuation of a team of mobile actuators to control a DPS. Potential applications include wildland firefighting using unmanned aerial vehicles and oil spill removal using autonomous skimmer boats. %Also, the simple formulation allows it to be combined with the model predictive control for implementation. 

% \subsection{Notation and terminology}

The paper adopts the following notation. The symbols $\mathbb{R}$, $\mathbb{R}^+$, and $\mathbb{N}$ denote the set of real numbers, nonnegative real numbers, and nonnegative integers, respectively. The boundary of a set $M$ is denoted by $\partial M$. The $n$-nary Cartesian power of a set $M$ is denoted by $M^n$. A continuous embedding is denoted by $\hookrightarrow$. % For better understanding of continuous embedding, see Definition 7 here: http://www-m6.ma.tum.de/~kuttler/script_reaktdiff.pdf
% Improvement for future revision: The notation $X_1 \hookrightarrow X_2$ means that the space $X_1$ is densely and continuously embedded in $X_2$.
We use $\abs{\cdot}$ and $\norm{\cdot}$ for the norm defined on a finite- and infinite-dimensional space, respectively, with subscript indicating type. The superscript ${}^*$ denotes an optimal variable or an optimal value, whereas ${}^{\star}$ denotes the adjoint of a linear operator. The transpose of a matrix $A$ is denoted by $A^\top$. An $n \times n$-dimensional identity matrix is denoted by $I_n$. We denote by $0_{n \times m}$ and $1_{n \times m}$ an $n \times m$-dimensional matrix with all entries being $0$ and $1$, respectively. 
%An $n \times n$-dimensional diagonal matrix with elements of vector $[a_1,a_2,\dots,a_n]$ on the main diagonal is denoted by $\text{diag}(a_1,a_2,\dots,a_n)$. 
The term \textit{guidance} refers to the steering of the mobile actuators, whereas the term \textit{control} refers to the control input to the DPS. For an optimization problem (P0) that minimizes cost function $J(\cdot)$ over variable $x$ subject to constraints, we use $\costeval{\text{(P0)}}{x}$ to denote the cost function of (P0) evaluated at $x$. Specifically, $\optcosteval{\text{(P0)}}{x}$ indicates that the optimal value of (P0) is attained when the cost function is evaluated at $x$.
% For an optimization problem
% \begin{equation}
% \begin{aligned}
%     & \underset{x \in X}{\text{minimize}} && J(x) \\ %\nonumber \\
%     & \text{subject to} & & x \in C,
% \end{aligned}
% \tag{P0}
% \label{prob: notation example problem}
% \end{equation}
% we use $\costeval{\eqref{prob: notation example problem}}{x}$ to denote the cost function of \eqref{prob: notation example problem} evaluated at $x$. Specifically, $\optcosteval{\eqref{prob: notation example problem}}{x}$ indicates that the optimal value of \eqref{prob: notation example problem} is attained when the cost function is evaluated at $x$.

% \subsection{Paper organization}
Section~II introduces relevant mathematical background, including representation of a diffusion-advection equation by an infinite-dimensional system, the associated LQ optimal control, and its finite-dimensional approximation. Section~III introduces the proposed optimization problem and its equivalent problem. Conditions for the existence of a solution are stated. Section~IV details the computation of an optimal solution using finite-dimensional approximations. Conditions for the convergence to the exact optimal solution of the approximate optimal solution are stated. A gradient-based method is applied to find an optimal solution. Section~V provides two numerical examples to illustrate optimal guidance and control solved by the proposed method. Section~VI summarizes the paper and discusses ongoing work.

%\singlespacing
% \iffalse This is not. \fi This is typeset again.

\section{Background}
This paper is motivated by the problem of controlling the following diffusion-advection process on a two-dimensional spatial domain $\Omega = [0,1] \times [0,1]$ with a team of $m_a$ mobile actuators:
\begin{align}
    % \frac{\partial z(x,y,t)}{\partial t} = & a (\frac{\partial^2 z(x,y,t)}{\partial x^2} + \frac{\partial^2 z(x,y,t)}{\partial y^2}) - (v_x \frac{\partial z(x,y,t)}{\partial x} + \frac{ v_y \partial z(x,y,t)}{\partial y}) \nonumber \\
    \frac{\partial z(x,y,t)}{\partial t} = & \ a \nabla^2 z(x,y,t) - \mathbf{v} \cdot \nabla z(x,y,t) \nonumber \\
    & \ + \sum_{i=1}^{m_a} (\mathcal{B}_i u_i) (x,y,t), \label{eq: motivating dynamical system 1}\\
    z(\cdot,\cdot,t)|_{\partial \Omega} = & \ 0, \label{eq: motivating dynamical system 2}\\
    z(x,y,0) = & \ z_0(x,y), \label{eq: motivating dynamical system 3}
\end{align}
where $z(\cdot,\cdot,t)$ is the state at time $t$, $\mathbf{v} \in \mathbb{R}^2$ is the velocity field for advection, and $u_i$ is the control implemented by actuator $i$, with the actuation characterized spatially by~$\mathcal{B}_i$. The state $z$ lives in the state space $L^2(\Omega)$. A representative model of the actuation dispensed by each actuator is Gaussian-shaped and centered at the actuator $i$'s location $(x_i,y_i)$ with a bounded support such that \begin{equation}\label{eq: example input operator}
    \mathcal{B}_i(x,y
    ) = 
    \left\{ \begin{aligned}
        \frac{1}{2 \pi \sigma_i^2} & \text{exp}  \left(-\frac{(x-x_i)^2}{\sigma_i^2} - \frac{(y-y_i)^2}{\sigma_i^2} \right)\\
        & \quad \text{if } |x-x_i|\leq \sigma_i \text{ and } |y-y_i| \leq \sigma_i,\\
        0 & \quad \text{otherwise},
    \end{aligned}\right.
\end{equation}
where the parameter $\sigma_i$ determines the spatial influence of the actuation, which is concentrated mostly at the location of the actuator and disperses to the surrounding with an exponential decay.

\subsection{Dynamics of the mobile actuators} 
Assume the mobile actuators have linear dynamics, so that the dynamics of actuator $i$ are
\begin{equation}\label{eq: agentwise linear dynamics}
    \dot{\xi}_i(t) = \alpha_i \xi_i(t) + \beta_i p_i(t), \quad \xi_i(0) = \xi_{i,0},
\end{equation}
where $\xi_i(t) \in \mathbb{R}^{n}$ $(n \geq 2)$ and $p_i(t) \in P_i \subset \mathbb{R}^{m}$ are the state and guidance at $t$, respectively. Assume that system \eqref{eq: agentwise linear dynamics} is controllable. The first two elements of $\xi_i(t)$ are the horizontal and vertical position, $x_i(t)$ and $y_i(t)$, of the actuator in the 2D domain. One special case would be a single integrator, where $\xi_i(t) \in \mathbb{R}^2$ is the position, $p_i(t) \in \mathbb{R}^2$ is the velocity commands, $\alpha_i =  0_{2 \times 2}$, and $\beta_i =  I_2$.

% The last two elements are the horizontal and vertical velocities $\dot{x}_i(t)$ and $\dot{y}_i(t)$. The guidance $p_i(t)$ contains $p_{i,x}(t) \in \mathbb{R}$ along the horizontal axis and $p_{i,y}(t) \in \mathbb{R}$ along the vertical axis.

For conciseness, we concatenate the states and guidance of all actuators, respectively, and use one dynamical equation to characterize the dynamics of all agents:
\begin{equation}\label{eq: general dynamics of the mobile actuator}
    \dot{\xi}(t) = \alpha \xi(t) + \beta p(t), \quad \xi(0) = \xi_0,
\end{equation}
where matrices $\alpha$ and $\beta$ are assembled from $\alpha_i$ and $\beta_i$ for $i \in \{1,2,\dots,m_a\}$, respectively and are consistent with the concatenation for $\xi$ and $p$.
% \begin{align}
%     \xi^\top(t) = & \begin{bmatrix}
%     \xi_1^\top(t) & \xi_2^\top(t) & \dots & \xi_{m_a}^\top(t)
%     \end{bmatrix}, \nonumber \\
%     \xi_0^\top = & \begin{bmatrix}
%     \xi_{1,0}^\top & \xi_{2,0}^\top & \dots & \xi_{m_a,0}^\top
%     \end{bmatrix}, \nonumber \\
%     p^\top(t) = & \begin{bmatrix}
%         p_1^\top(t) & p_2^\top(t) & \dots & p_{m_a}^\top(t)
%     \end{bmatrix}, \nonumber \\
%     a = & \text{diag}(a_1,a_2,\dots,a_{m_a}), \nonumber \\
%     b = & \text{diag}(b_1,b_2\dots, b_{m_a}). %{\color{red} keep it?}
% \end{align}
With a slight abuse of notation, we use $n$ for the dimension of $\xi(t)$ and $m$ for the dimension of $p(t)$. Define the admissible set of guidance $P \defeq P_1 \times P_2 \times \dots \times P_{m_a}$ such that $p(t) \in P$ for $t \in [0,t_f]$. Let $M \in \mathbb{R}^{2 m_a \times n}$ be a matrix such that $M \xi(t)$ is a vector of locations of the actuators.

\subsection{Abstract linear system and linear-quadratic regulation} To describe the dynamics of PDE \eqref{eq: motivating dynamical system 1}--\eqref{eq: motivating dynamical system 3}, consider the following abstract linear system:
\begin{equation}
    \dot{\mathcal{Z}}(t) = \mathcal{A} \mathcal{Z}(t) + \mathcal{B}(M\xi(t),t)u(t), \qquad \mathcal{Z}(0) = \mathcal{Z}_0, \label{eq: new abstract linear system}
\end{equation}
where $\mathcal{Z}(\cdot)$ is the state within state space $\mathcal{H} \defeq L^2(\Omega)$ and $u(\cdot)$ is the control within the control space ${u(t) \in U \subseteq \mathbb{R}^{m_a}}$ for $t \in [0,t_f]$. 
In the case of diffusion-advection process \eqref{eq: motivating dynamical system 1}, for $\phi \in \mathcal{H}$, 
\begin{equation}\label{eq: definition of diffusion-advection operator}
    (\mathcal{A} \phi ) (x,y) = a \nabla^2 \phi(x,y) - \mathbf{v} \cdot \nabla \phi(x,y),
\end{equation}
where the operator $\mathcal{A}$ has domain $\text{Dom}(\mathcal{A}) = H^2(\Omega) \cap H_0^1(\Omega)  $. %$\text{Dom}(\mathcal{A}) = \{\psi \in H_0^1(\Omega) \cap H^2(\Omega) | \psi(\cdot,\cdot,t)|_{\partial \Omega} = 0 , t \in [0,t_f]\}$. 
The input operator $\mathcal{B}(M\xi(t),t) \in \mathcal{L}(U,\mathcal{H})$ is a function of the actuator locations such that $\mathcal{B}(M\xi(t),t) = [\mathcal{B}_1(M\xi_1(t),t), \dots,\mathcal{B}_{m_a}(M\xi_{m_a}(t),t)]^\top$, where $\mathcal{B}_{i}(\cdot,t) \in L^2(\Omega)$ for all $t \in [0,t_f]$ and $i \in \{1,2,\dots,m_a \}$. A special case is the time-invariant input operator in \eqref{eq: example input operator}. Since the actuator state $\xi(t)$ is a function of time $t$, we sometimes use $\mathcal{B}(t)$ for brevity. %Assume $\mathcal{B}(\cdot) \in L^2([0,t_f];\mathcal{L}(U,\mathcal{H}))$.

The operator $\mathcal{A}: \text{Dom}(\mathcal{A}) \rightarrow \mathcal{H}$ is an infinitesimal generator of a strongly continuous semigroup $\mathcal{S}(t)$ on $\mathcal{H}$.
% Assume that $\mathcal{A}: \text{Dom}(\mathcal{A}) \rightarrow \mathcal{H}$ is an infinitesimal generator of a strongly continuous semigroup $\mathcal{S}(t)$ on $\mathcal{H}$.
Subsequently, the dynamical system \eqref{eq: new abstract linear system} has a unique mild solution $\mathcal{Z} \in C([0,t_f];\mathcal{H})$ for any $\mathcal{Z}_0 \in \mathcal{H}$ and any $u \in L^2([0,t_f];U)$ such that $\mathcal{Z}(t) = \mathcal{S}(t)\mathcal{Z}_0 + \int_0^t \mathcal{S}(t-\tau) \mathcal{B}(\xi(\tau),\tau)u(\tau) \dd \tau$.

Similar to a finite-dimensional linear system, a linear-quadratic regulator (LQR) problem can be formulated with respect to \eqref{eq: new abstract linear system}, which looks for a control $u(\cdot) \in L^2([0,t_f];U)$ that minimizes the following quadratic cost:
\begin{align}
    J(\mathcal{Z},u) \defeq & \int_0^{t_f} \innerproduct{\mathcal{Z}(t)}{\mathcal{Q}(t)\mathcal{Z}(t)} + u(t)^\top R u(t) \dd t  \nonumber \\
    &+ \innerproduct{\mathcal{Z}(t_f)}{\mathcal{Q}_f \mathcal{Z}(t_f)}, \label{eq: new quadratic PDE cost}
\end{align}
where $\mathcal{Q}(t) \in \mathcal{L}(\mathcal{H})$ and $\mathcal{Q}_f \in \mathcal{L}(\mathcal{H})$ are self-adjoint and nonnegative, which evaluates the running cost and terminal cost of the PDE state. The coefficient $R$ is an $m_a \times m_a$-dimensional symmetric and positive definite matrix that evaluates the control effort. We refers to $J(\mathcal{Z},u)$ as the \textit{PDE cost}.

Analogous to the finite-dimensional LQR, an optimal control $u^*$ that minimizes the quadratic cost \eqref{eq: new quadratic PDE cost} is 
\begin{equation}
    u^*(t) = -R^{-1} \mathcal{B}^{\star}(t) \Pi(t) \mathcal{Z}(t), \label{eq: new optimal LQR feedback control}
\end{equation}
where $\Pi$ is an operator that associates with the following backward differential operator-valued Riccati equation:
\begin{multline}
    \dot{\Pi}(t) =  -\mathcal{A}^{\star} \Pi(t) - \Pi(t) \mathcal{A} - \mathcal{Q}(t)  \\
     + \Pi(t) \bar{\mathcal{B}} \bar{\mathcal{B}}^{\star}(t) \Pi(t) \label{eq: new operator Riccati equation}
\end{multline}
with terminal condition $\Pi(t_f) = \mathcal{Q}_f$, where $\bar{\mathcal{B}} \bar{\mathcal{B}}^{\star}(t)$ is short for $\mathcal{B}(t)R^{-1} \mathcal{B}^{\star}(t)$. Before we proceed to state the conditions for the existence of a unique solution of \eqref{eq: new operator Riccati equation}, we introduce the $\mathcal{J}_q$-class as follows.

Denote the trace of a nonnegative operator $A \in \mathcal{L}(\mathcal{H})$ by $\trace{A}$, where $\trace{A} \defeq \sum_{k=1}^{\infty} \innerproduct{\phi_k}{A \phi_k}$ for any orthonormal basis $\{\phi_k \}_{k=1}^{\infty}$ of $\mathcal{H}$ (the trace is independent of the choice of basis functions). % \cite{burns2015SIAM}
For $1 \leq q < \infty$, let $\mathcal{J}_q(\mathcal{H})$ denote the set of all bounded operators $\mathcal{L}(\mathcal{H})$ such that $\trace{(\sqrt{A^{\star}A})^q} < \infty$ \cite{Burns2015Solutions}. %[Definition~3.2]  
If $A \in \mathcal{J}_{q}(\mathcal{H})$, then the $\mathcal{J}_q$-norm of $A$ is defined as $\norm{A}_{\mathcal{J}_q(\mathcal{H})} \defeq (\trace{(\sqrt{A^{\star}A})^q})^{1/q}< \infty$.
The class $\mathcal{J}_1(\mathcal{H})$ and $\mathcal{J}_2(\mathcal{H})$ are known as the space of trace operators and the space of Hilbert-Schmidt operators, respectively. Note that a continuous embedding ${\mathcal{J}_{q_1}(\mathcal{H}) \hookrightarrow \mathcal{J}_{q_2}(\mathcal{H})}$ holds if $1 \leq q_1 < q_2 \leq \infty$. In other words, if $A \in \mathcal{J}_{q_1}(\mathcal{H})$, then $A \in \mathcal{J}_{q_2}(\mathcal{H})$ and ${ \norm{A}_{\mathcal{J}_{q_2}(\mathcal{H})} \leq \norm{A}_{\mathcal{J}_{q_1}(\mathcal{H})} }$.%\cite{Burns2015Infinitedimensional}.

The existence of a mild solution of \eqref{eq: new operator Riccati equation} is established via Lemma~\ref{lemma: existence of Riccati mild solution}. We omit the proof of this lemma because it is a direct consequence of \cite[Theorem~3.6]{Burns2015Solutions}.

Consider the following assumptions with $1 \leq q < \infty$:
\begin{enumerate}
        \item[(A1)] $\mathcal{Q}_f \in \mathcal{J}_q(\mathcal{H})$ and $\mathcal{Q}_f$ is nonnegative.
        \item[(A2)] $\mathcal{Q}(\cdot) \in L^1([0,t_f];\mathcal{J}_q(\mathcal{H}))$ and $\mathcal{Q}(t)$ is nonnegative for all $t \in [0,t_f]$.
        \item[(A3)] $\bar{\mathcal{B}} \bar{\mathcal{B}}^{\star}(\cdot) \in L^{\infty}([0,t_f];\mathcal{L}(\mathcal{H}))$ and $\bar{\mathcal{B}} \bar{\mathcal{B}}^{\star}(t)$ is nonnegative for $t \in [0,t_f]$.
\end{enumerate}
% Our result should hold for q=1, the trace class, which is also the Hilbert-Schmidt class that yield kernel representation.

\begin{lem}\label{lemma: existence of Riccati mild solution}
    Let $\mathcal{H}$ be a separable Hilbert space and let $\mathcal{S}(t)$ be a strongly continuous semigroup on $\mathcal{H}$. Suppose assumptions (A1)--(A3) hold. Then, the equation 
    \begin{multline}
        \Pi(t) = \mathcal{S}^{\star}(t_f-t) \mathcal{Q}_f \mathcal{S}(t_f-t) + \int_t^{t_f} \mathcal{S}^{\star}(\tau-t) \\ \left( \mathcal{Q}(\tau) - \Pi(\tau) \bar{\mathcal{B}} \bar{\mathcal{B}}^{\star}(\tau) \Pi(\tau) \right) \mathcal{S}(\tau-t) \dd \tau \label{eq: new mild solution of operator Riccati equation}
    \end{multline}
provides a unique mild solution to \eqref{eq: new operator Riccati equation} in the space $L^2([0,t_f];\mathcal{J}_{2q}(\mathcal{H}))$. The solution also belongs to $C([0,t_f];$ $\mathcal{J}_q(\mathcal{H}))$ and is pointwise self-adjoint and nonnegative.
Furthermore, if $\mathcal{Q}(\cdot) \in C([0,t_f];\mathcal{J}_q(\mathcal{H}))$ and $\bar{\mathcal{B}} 
\bar{\mathcal{B}}^{\star}(\cdot) \in C([0,t_f];\mathcal{L}(\mathcal{H}))$, then $\Pi$ is a weak solution to \eqref{eq: new operator Riccati equation}.
\end{lem}

The equality introduced next in Lemma~\ref{lemma: equivalent PDE cost} allows for turning the optimal quadratic PDE cost into a quadratic term associated with the initial condition of the PDE and the Riccati operator. We state it without proof because it can be established by integrating $\dd  \innerproduct{\mathcal{Z}(t)}{\Pi(t) \mathcal{Z}(t)} / \dd t$ from $0$ to $t_f$; the differentiability of $\innerproduct{\mathcal{Z}(t)}{\Pi(t) \mathcal{Z}(t)}$ is proven in \cite[Theorem~6.1.9]{curtain2012introduction}.

\begin{lem}\label{lemma: equivalent PDE cost}
    Suppose $\Pi(t)$ is a mild solution to \eqref{eq: new operator Riccati equation}, given by \eqref{eq: new mild solution of operator Riccati equation}. For every $\mathcal{Z}_0 \in \mathcal{H}$, the optimal PDE cost \eqref{eq: new quadratic PDE cost} satisfies the equality $ J(\mathcal{Z}^*,u^*) = \innerproduct{\mathcal{Z}_0}{\Pi(0) \mathcal{Z}_0}$, where $\mathcal{Z}^*$ is the state that follows the dynamics \eqref{eq: new abstract linear system} under optimal control $u^*$ of \eqref{eq: new optimal LQR feedback control}, and $\Pi(0)$ is the solution \eqref{eq: new mild solution of operator Riccati equation} evaluated at $t=0$.
\end{lem}

% \begin{pf*}{Proof}
%     See Appendix~\ref{prf: establish equivalent PDE cost}.
% \qed 
% \end{pf*}

The following assumption is vital to the main results in this paper.

\begin{enumerate}
\item[(A4)] The input operator $\mathcal{B}_i(x,t)$ is continuous with respect to location $x \in \mathbb{R}^2$ \cite[Definition~4.5]{burns2015infinite}, that is, there exists a continuous function $l: \mathbb{R}^+ \rightarrow \mathbb{R}^+$ such that $l(0) = 0$ and $ \norm{\mathcal{B}_i(x,t)-\mathcal{B}_i(y,t) }_{L^2(\Omega)} \leq l(|x-y|_2)$ for all $t \in [0,t_f]$, all $x,y \in \mathbb{R}^2$, and all $i \in \{1,2,\dots,m_a \}$.
\end{enumerate}

% Since $\mathcal{B}(M\xi(t),t)$ maps locations of the actuators $M\xi(t)$ to a bounded linear operator, the solution $\Pi(0)$ of \eqref{eq: new mild solution of operator Riccati equation} at $t = 0$ is a mapping that maps the actuators' trajectory $M\xi$ to a $\mathcal{J}_q$-valued operator. Furthermore, the optimal PDE cost $\innerproduct{\mathcal{Z}_0}{\Pi(0) \mathcal{Z}_0}$ is a mapping that maps the actuators' trajectory to a nonnegative real number.

\revision{The actuators' locations determine where the input is actuated and, furthermore, how $\Pi(\cdot)$ evolves through~\eqref{eq: new mild solution of operator Riccati equation}. %We characterize this relation by a composite mapping. 
Since the input operator $\mathcal{B}(\cdot,t)$ is a mapping of the actuators' locations at time $t$, the composite input operator $\bar{\mathcal{B}} \bar{\mathcal{B}}^\star (\cdot)$ is a mapping of the actuator state in $[0,t_f]$ and so is $\Pi(0)$ by \eqref{eq: new mild solution of operator Riccati equation}, although the actuator state is not explicitly reflected in the notation of $\bar{\mathcal{B}} \bar{\mathcal{B}}^\star (\cdot)$ or $\Pi(0)$. % The value of $\Pi(0)$ depends on the actuator state $\xi$ despite the fact that $\Pi(0)$ is propagated backwards in time.
Hence, we can define the optimal PDE cost $\innerproduct{\mathcal{Z}_0}{\Pi(0) \mathcal{Z}_0}$ as a mapping of the actuator state. Let $K:C([0,t_f];\mathbb{R}^n) \rightarrow \mathbb{R}^+$ such that $K(\zeta) \defeq \innerproduct{\mathcal{Z}_0}{\Pi(0) \mathcal{Z}_0}$. 
Assumption (A4) plays an important role in yielding the continuity of the mapping $K(\cdot)$ stated below in Lemma~\ref{lemma: continuity of equivalent cost wrt actuator location}, whose proof is in the supplementary material% \cite{Cheng2020optimalControlSupplement}.
.}
% Lemma~\ref{lemma: continuity of trace cost wrt sensor trajectory} below shows when the uncertainty cost varies continuously with respect to the sensor state.}

\begin{lem}\label{lemma: continuity of equivalent cost wrt actuator location}
Suppose $\mathcal{Z}_0 \in \mathcal{H}$. Let assumptions (A1)--(A3) hold with $q = 1$ and $\Pi \in C([0,t_f];\mathcal{J}_1(\mathcal{H}))$ be defined as in \eqref{eq: new mild solution of operator Riccati equation}. If assumption (A4) holds, then the mapping $K: C([0,t_f];\mathbb{R}^n) \rightarrow \mathbb{R}^+$ such that $K(\xi) \defeq \innerproduct{\mathcal{Z}_0}{\Pi(0) \mathcal{Z}_0}$ is continuous.
\end{lem}

% \begin{pf*}{Proof}
%     See Appendix~\ref{prf: continuity of Riccati operator wrt location}. 
% \qed 
% \end{pf*}

Approximations to \eqref{eq: new abstract linear system} and \eqref{eq: new mild solution of operator Riccati equation} permit numerical computation. Consider a finite-dimensional subspace $\mathcal{H}_N \subset \mathcal{H}$ with dimension $N$. The inner product and norm of $\mathcal{H}_N$ are inherited from that of $\mathcal{H}$. Let $P_N: \mathcal{H} \to \mathcal{H}_N$ denote the orthogonal projection of $\mathcal{H}$ onto $\mathcal{H}_N$. Let $Z_N(t) \defeq P_N \mathcal{Z}(t)$ and $S_N(t) \defeq P_N \mathcal{S}(t) P_N$ denote the finite-dimensional approximation of $\mathcal{Z}(t)$ and $\mathcal{S}(t)$, respectively. A finite-dimensional approximation of \eqref{eq: new abstract linear system} is
\begin{align}
    \dot{Z}_N(t) = & \ A_N Z_N(t) + B_N(M\xi(t),t) u(t), \label{eq: finite dim approx of abstract linear system}\\ 
    Z_N(0) = & \ Z_{0,N} \defeq P_N \mathcal{Z}_0, \label{eq: initial condition finite dim approx of abstract linear system}
\end{align}
where $A_N \in \mathcal{L}(\mathcal{H}_N)$ and $B_N(M\xi(t),t) \in \mathcal{L}(U,\mathcal{H}_N)$ are approximations of $\mathcal{A}$ and $\mathcal{B}(M\xi(t),t)$%\cite{morris2010linear}
, respectively. Since the actuator state $\xi(t)$ is a function of time $t$, we sometimes use $B_N(t)$ for brevity. %Similar to the usage of $\mathcal{B}(t)$ for brevity, we sometimes use $B_N(t)$ for the same cause. 
Correspondingly, the finite-dimensional approximation of \eqref{eq: new mild solution of operator Riccati equation} is 
\begin{multline}\label{eq: mild solution of approximate Riccati}
        \Pi_{N}(t) =  S^{\star}_N(t_f-t) Q_{f,N} S_N(t_f-t) 
         + \int_t^{t_f} S^{\star}_N(\tau-t) \\
         \left( Q_N(\tau) - \Pi_{N}(\tau) \bar{B}_N \bar{B}_N^{\star}(\tau) \Pi_{N}(\tau) \right) S_N(\tau-t) \dd \tau ,
    \end{multline}
% \begin{align}
%     \dot{S}_N(t) = & -(A_N)^T S_N(t) - S_N(t) A_N - Q_N \nonumber \\
%     & + S_N(t)B_{\xi,N}(t)R^{-1}B_{\xi,N}^T(t)S_N(t), \label{eq: finite-dimensional approximation of the operator riccati equation} \\
%     S_N(t_f) =& \ Q_{fN},
% \end{align}
where $Q_N = P_N \mathcal{Q} P_N$, $Q_{fN} = P_N \mathcal{Q}_f P_N$, and $\bar{B}_N \bar{B}_N^{\star}(\tau)$ is short for $B_N(\tau) R^{-1} B_N^{\star}(\tau)$.

The optimal control $u_N^*$ that minimizes the approximated PDE cost
\begin{multline}
    J_N(Z_N,u_N) \defeq \innerproduct{Z_N(t_f)}{Q_{f,N} Z_N(t_f)} \\+ \int_0^{t_f} \innerproduct{Z_N(t)}{Q_N(t)Z_N(t)} + u_N^\top(t) R u_N(t) \dd t   
     \label{eq: approximated quadratic PDE cost}
\end{multline}
is analogous to \eqref{eq: new optimal LQR feedback control}: 
\begin{equation}
    u_N^* = -R^{-1} B_N^{\star}(t) \Pi_N(t) Z_N(t),   \label{eq: approximated optimal lqr feedback control} 
\end{equation}
where $\Pi_N(t)$ is a solution of \eqref{eq: mild solution of approximate Riccati}.

The following assumptions are associated with the approximations:
\begin{enumerate}
    \item[(A5)] Both $\mathcal{Q}_f$ and sequence $\{Q_{f,N} \}_{N=1}^{\infty}$ are elements of {$\mathcal{J}_q(\mathcal{H})$}. Both $\mathcal{Q}_f$ and $Q_{f,N}$ are nonnegative for all $N \in \mathbb{N}$ and $\norm{\mathcal{Q}_f-Q_{f,N}}_{\mathcal{J}_q(\mathcal{H})} \rightarrow 0$ as $ N \rightarrow \infty$.
    \item[(A6)] Both $\mathcal{Q}(\cdot)$ and sequence $\{Q_N(\cdot) \}_{N=1}^{\infty}$ are elements of $L^1([0,t_f];\mathcal{J}_q(\mathcal{H}))$. Both $\mathcal{Q}(\tau)$ and $Q_N(\tau)$ are nonnegative for all $\tau \in [0,t_f]$ and all $N \in \mathbb{N}$ and satisfy $\int_0^t \norm{\mathcal{Q}(\tau)- Q_N(\tau)}_{\mathcal{J}_q(\mathcal{H})} \dd \tau \rightarrow 0$ for all $t \in [0,t_f]$ as $N \rightarrow \infty$.
    \item[(A7)] Both $\bar{\mathcal{B}} \bar{\mathcal{B}}^{\star}(\cdot)$ and sequence $\{\bar{B}_N \bar{B}_N^{\star}(\cdot)\}_{N=1}^{\infty}$ are elements of $L^{\infty}([0,t_f];\mathcal{L}(\mathcal{H}))$. Both $\bar{\mathcal{B}} \bar{\mathcal{B}}^{\star}(t) $ and $\bar{B}_N \bar{B}_N^{\star}(t) $ are nonnegative for all $t \in [0,t_f]$ and all $N \in \mathbb{N}$ and satisfy 
        \begin{equation}
            \underset{t \in [0,t_f]}{\esssup} \norm{\bar{\mathcal{B}}\bar{\mathcal{B}}^{\star}(t) - \bar{B}_N \bar{B}_N^{\star}(t)}_{\text{op}} \rightarrow 0
        \end{equation}
        as $N \rightarrow \infty$ ($\norm{\cdot}_{\text{op}}$ denotes the operator norm).
\end{enumerate}

Note that the assumptions (A1), (A2), and (A3) are contained in (A5), (A6), and (A7), respectively.

The next theorem states the convergence of an approximate solution of the Riccati equation, which is reproduced from \cite[Theorem~3.5]{Burns2015Solutions} and hence stated without a proof.

\begin{thm}\label{thm: convergence of riccati operator}
    Suppose $\mathcal{S}(t)$ is a strongly continuous semigroup of linear operators over a Hilbert space $\mathcal{H}$ and that $\{S_N(t)\}$ is a sequence of uniformly continuous semigroup over the same Hilbert space that satisfy, for each $\phi \in \mathcal{H}$
    \begin{equation}
        \norm{\mathcal{S}(t) \phi - S_N(t) \phi} \rightarrow 0, \quad \norm{\mathcal{S}^{\star}(t) \phi - S_N^{\star}(t) \phi} \rightarrow 0 % p. 213 of \cite{Burns2015infinte} justifies these two conditions
    \end{equation}
    as $N \rightarrow \infty$, uniformly in $[0,t_f]$. Suppose assumptions (A5)--(A7) hold. If $\Pi(\cdot) \in C([0,t_f];\mathcal{J}_q(\mathcal{H}))$ is a solution of \eqref{eq: new mild solution of operator Riccati equation} and $\Pi_{N}(\cdot) \in C([0,t_f];\mathcal{J}_q(\mathcal{H}))$ is the sequence of solution of \eqref{eq: mild solution of approximate Riccati}, then
    \begin{equation}\label{eq: convergence of the Riccati operator}
        \underset{t \in [0,t_f]}{ \sup} \norm{\Pi(t)-\Pi_{N}(t)}_{\mathcal{J}_q(\mathcal{H})} \rightarrow 0
    \end{equation}
    as $N \rightarrow \infty$.
\end{thm}

The following assumption and lemma are analogous to (A4) and Lemma~\ref{lemma: continuity of equivalent cost wrt actuator location}, respectively:
\begin{enumerate}
    \item[(A8)] The approximated input operator $B_{i,N}(x,t)$ is continuous with respect to location $x \in \mathbb{R}^2$, that is, there exists a continuous function $l_N: \mathbb{R}^+ \rightarrow \mathbb{R}^+$ such that $l_N(0) = 0$ and $\norm{B_{i,N}(x,t) - B_{i,N}(y,t)}_{L^2(\Omega)} \leq l_N(|x-y|_2)$ for all $t \in [0,t_f]$, all $x,y \in \mathbb{R}^2$, and all $i \in \{1,2,\dots,m_a \}$.
\end{enumerate}

Similar to the mapping $K(\cdot)$ in Lemma~\ref{lemma: continuity of equivalent cost wrt actuator location}, the optimal approximated PDE cost can be characterized as a mapping of the actuator state through \eqref{eq: mild solution of approximate Riccati}, where the continuity is established in Lemma~\ref{lem: continuity of approximated equivalent cost wrt actuator location}, whose proof is in the supplementary material. % \cite{Cheng2020optimalControlSupplement}.

\begin{lem}\label{lem: continuity of approximated equivalent cost wrt actuator location}
Suppose $Z_{0,N} \in \mathcal{H}_N$. Let assumptions (A5)--(A7) hold and $\Pi_{N}(t)$ be defined as in \eqref{eq: mild solution of approximate Riccati}. If assumption (A8) holds, then the mapping $K_N: C([0,t_f];\mathbb{R}^n) \rightarrow \mathbb{R}^+$ such that $K_N(\xi) \defeq \innerproduct{Z_{0,N}}{\Pi_{N}(0) Z_{0,N}}$ is continuous.
\end{lem}

\section{Problem formulation}

% \subsection{Integrated optimal control problem}
% This subsection introduces the formulation of the integrated optimal control problem, which simultaneously solves for optimal (open-loop) guidance of the mobile actuators and optimal (open-loop) control of the DPS. Specifically, the cost functions evaluating actuator motions and DPS regulation are integrated into one cost function. Consequently, the dynamics of the DPS and that of the mobile actuator are both constraints. The integrated problem is 
  
This paper seeks to derive the guidance and control input of each actuator such that the state $\mathcal{Z}$ of the abstract linear system \eqref{eq: new abstract linear system} can be driven to zero. Specifically, consider the following problem:
\begin{equation}
	\begin{aligned}
	& \underset{ \substack{u \in L^2([0,t_f];U)\\p \in L^2([0,t_f];P)}}{\text{minimize}}
	& & J(\mathcal{Z},u) + J_{\text{m}}(\xi,p) \\
	& \text{subject to} 
	& &  \dot{\mathcal{Z}}(t) = \mathcal{A} \mathcal{Z}(t) + \mathcal{B}(t)u(t), \quad \mathcal{Z}(0) = \mathcal{Z}_0,\\
	&
	& & \dot{\xi}(t) = \alpha \xi(t) + \beta p(t), \quad \xi(0) = \xi_0,\\
	\end{aligned}
	\tag{P}
	\label{prob: new IOCA}
\end{equation}
where $J_{\text{m}}(\xi,p) \defeq \int_0^{t_f} h(\xi(t),t) + g(p(t),t) \dd t + h_f(\xi(t_f))$ is the cost associated with the motion of the actuators, named the \textit{mobility cost}, such that the mappings $h: \mathbb{R}^{n} \times [0,t_f] \rightarrow \mathbb{R}^+$ and $g: \mathbb{R}^{m} \times [0,t_f] \rightarrow \mathbb{R}^+$ evaluate the running state cost and running guidance cost, respectively, and the mapping $h_f: \mathbb{R}^{n} \rightarrow \mathbb{R}^+$ evaluates the terminal state cost. 

The running state cost $h(\cdot,\cdot)$ may characterize restrictions to actuator state. For example, a Gaussian-type function with its peak in the center of the spatial domain, i.e.,
\begin{multline}
    h(\left[\begin{smallmatrix}
    x \\ y
    \end{smallmatrix} \right],t) = \\ \frac{1}{2 \pi \sigma_x(t) \sigma_y(t)}  \text{exp}  \left(-\frac{(x-0.5)^2}{\sigma_x^2(t)} - \frac{(y-0.5)^2}{\sigma_y^2(t)} \right),
\end{multline}
where $\sigma_x(t), \sigma_y(t) > 0$ and $x,y \in [0,1]$, can model a hazardous field that may shorten the life span of an actuator. The integral of this function in the interval $[0,t_f]$ evaluates the accumulated exposure of the mobile actuator along its trajectory, which may need to be contained as small as possible (see \cite{demetriou2020navigating}%\cite{demetriou2018incorporating,demetriou2020navigating}
). Another example is the artificial potential field \cite{hoy2015algorithms}, cast as a soft constraint, that penalizes the trajectory when it passes an inaccessible region such as an obstacle. 
The running guidance cost $g(\cdot,\cdot)$ may be the absolute value or a quadratic function of the guidance, which characterizes the total amount (of fuel) or energy for steering, respectively. 
And the terminal state cost $h_f(\cdot)$ may characterize restrictions of the terminal state of the mobile actuators. For example, if an application specifies terminal positions, then $h_f(\cdot)$ may be a quadratic function that penalizes the deviation of the actual terminal positions.

The formulation in \eqref{prob: new IOCA} provides an intermediate step for minimizing the PDE cost subject to mobility constraints, in addition to the dynamics constraints. The mobility constraints are characterized by inequalities of $h_f(\cdot)$ and the integrals of $h(\cdot,\cdot)$ and $g(\cdot,\cdot)$, because these constraints can be used to augment the cost function and turned into the form of \eqref{prob: new IOCA} using the method of Lagrange multipliers.

An equivalent problem of \eqref{prob: new IOCA} can be derived using Lemma~\ref{lemma: equivalent PDE cost}. For an arbitrary admissible guidance $p$, the actuator trajectory $\xi$ is determined following the dynamics \eqref{eq: general dynamics of the mobile actuator}, which also determines the input operator $\mathcal{B}(\xi(\cdot),\cdot)$. By Lemma~\ref{lemma: equivalent PDE cost}, the control $u$ that minimizes the cost function of \eqref{prob: new IOCA}---specifically, the PDE cost $J(\mathcal{Z},u)$---is given by \eqref{eq: new optimal LQR feedback control}, and the minimum PDE cost is $\innerproduct{\mathcal{Z}_0}{\Pi(0) \mathcal{Z}_0}$, where $\Pi(0)$ is the mild solution of \eqref{eq: new mild solution of operator Riccati equation} with actuator trajectory steered by guidance $p$. Hence, we derive the following problem equivalent to \eqref{prob: new IOCA}:
\begin{equation}
	\begin{aligned}
	& \underset{p \in L^2([0,t_f];P)}{\text{minimize}}
	& & \innerproduct{\mathcal{Z}_0}{\Pi(0) \mathcal{Z}_0} + J_{\text{m}}(\xi,p) \\
	& \text{subject to} 
	& &  \dot{\xi}(t) = \alpha \xi(t) + \beta p(t), \quad \xi(0) = \xi_0,
% 	&
% 	& & \Pi(0) = \mathcal{S}^{\star}(t_f) \mathcal{Q}_f \mathcal{S}(t_f) + \int_0^{t_f} \mathcal{S}^{\star}(\tau) \left( \mathcal{Q}(\tau) - \Pi(\tau) \mathcal{B}^{\star} (\tau) R^{-1} \mathcal{B}(\tau) \Pi(\tau) \right) \mathcal{S}(\tau) \dd \tau.
	\end{aligned}
	\tag{P1}
	\label{prob: new equivalent IOCA}
\end{equation}
where $\Pi(0)$ is defined in \eqref{eq: new mild solution of operator Riccati equation} with $t = 0$.

To prove the existence of a solution to \eqref{prob: new equivalent IOCA}, we make the following assumptions on the admissible set of guidance and the functions composing the mobility cost:
\begin{enumerate}
    %\item The set $X_0 \subset \mathbb{R}^{n_a m_a}$, where $\xi_0 \in X_0$, is convex and closed.
    \item[(A9)] The set of admissible guidance $P \subset \mathbb{R}^{m}$ is closed and convex.
    \item[(A10)] The mappings $h: \mathbb{R}^{n} \times [0,t_f] \rightarrow \mathbb{R}^+$, $g: \mathbb{R}^{m} \times [0,t_f] \rightarrow \mathbb{R}^+$, and $h_f: \mathbb{R}^{n} \rightarrow \mathbb{R}^+$ are continuous. For every $t \in [0,t_f]$, the function $g(\cdot,t)$ is convex.
    % The old A9 is removed because it is not very useful in the proof. Instead, we require all the cost functions in the mobility cost are nonnegative.
    %\item[(A9)] There exists $b_1,b_2,c_1,c_2 \in \mathbb{R}$ such that $h_f(\xi) \geq b_1 + c_1 |\xi|_2$ (\textbf{Euclidean norm $|\cdot|_2$}) for all $\xi \in \mathbb{R}^{n_a m_a}$ and $h(\xi,t) \geq b_2 + c_2 |\xi|_2$ for all $(\xi,t) \in \mathbb{R}^{n_a m_a} \times [0,t_f]$. Without loss of generality, we may assume $c_1,c_2 \leq 0$.
    \item[(A11)] There exists a constant $d_1 >0$ with $g(p,t) \geq d_1 |p|_2^2$  for all $(p,t) \in P \times [0,t_f]$.
    % \item[(A10)] The operator $\mathcal{B}_i(x,t)$ is continuous with respect to location for $i \in \{1,2,\dots, m_a \}$.
\end{enumerate}

Assumptions (A9)--(A11) are generally satisfied in applications with vehicles carrying the actuators.
Assumption (A9) is a physically reasonable characterization of the steering of a vehicle, where the admissible steering is generally a continuum with attainable limits within its range. Assumption (A10) places a general continuity requirement on the cost functions and a convexity requirement on the steering cost function. Assumption (A11) requires the function $g(p,t)$ to be bounded below by a quadratic function of the guidance $p$ for all $t$, which is generally satisfied, e.g., with $g$ itself being a quadratic function of $p$. These assumptions are applied in Theorem~\ref{thm: existence of a solution of IOCA} below regarding the existence of a solution of \eqref{prob: new equivalent IOCA}, whose proof is in Appendix~\ref{prf: existence of a solution of IOCA}. Subsequently, the solution to \eqref{prob: new equivalent IOCA} can be used to reconstruct the solutions to \eqref{prob: new IOCA}, which is stated in Theorem~\ref{thm: equivalence between P and P1} with its proof in Appendix~\ref{prf: equivalence between P and P1}.

\begin{thm}\label{thm: existence of a solution of IOCA}
    Consider problem \eqref{prob: new equivalent IOCA} and let assumptions (A1)--(A4) and (A9)--(A11) hold. %If there exists guidance $p_0$ in the admissible set $\mathcal{P}$ such that ${\costeval{\eqref{prob: new equivalent IOCA}}{p_0} < \infty}$, 
    Then \eqref{prob: new equivalent IOCA} has a solution.
    % Furthermore, the guidance $p^*(\cdot)$ constitutes an optimal solution of \eqref{prob: new IOCA}, whose optimal control $u^*(\cdot)$ is given by \eqref{eq: new optimal LQR feedback control} with actuators steered by $p^*(\cdot)$.
\end{thm}

% \begin{pf*}{Proof}
%     See Appendix~\ref{prf: existence of a solution of IOCA}.
% \qed 
% \end{pf*}

\begin{thm}\label{thm: equivalence between P and P1}
    Consider problems \eqref{prob: new IOCA} and \eqref{prob: new equivalent IOCA}. Let assumptions (A4) and (A9)--(A11) hold. Let $p^*$ be the optimal solution of \eqref{prob: new equivalent IOCA} and $u^*$ be the optimal control obtained from \eqref{eq: new optimal LQR feedback control} with actuator trajectory steered by $p^*$. Then $u^*$ and $p^*$ minimize problem \eqref{prob: new IOCA}.
\end{thm}

% \begin{pf*}{Proof}
%     See Appendix~\ref{prf: equivalence between P and P1}.
% \end{pf*}

The equivalent problem \eqref{prob: new equivalent IOCA} allows us to search for an optimal guidance $p$ such that the mobility cost plus the optimal PDE cost is minimized. The control is no longer an optimization variable, because it is determined by the LQR of the abstract linear system for arbitrary trajectories of the mobile actuators.

\section{Computation of optimal control and guidance}\label{sec: computation of optimal solution}
Approximation of the infinite-dimensional terms in problem \eqref{prob: new IOCA} is necessary when computing the optimal control and guidance. Hence, we replace the PDE cost and dynamics of \eqref{prob: new IOCA} by \eqref{eq: approximated quadratic PDE cost} and \eqref{eq: finite dim approx of abstract linear system}, respectively, and obtain the following approximate problem \eqref{prob: finite dim approx integrated optimization problem}:

\begin{equation}
	\begin{aligned}
	& \underset{ \substack{u \in L^2([0,t_f];U)\\p \in L^2([0,t_f];P)}}{\text{minimize}}
	& & J_N(Z_N,u) + J_\text{m}(\xi,p) \\
	& \text{subject to} 
	& &  \dot{Z}_N(t) = A_N Z_N(t) + B_N(M \xi(t),t)u(t) \\
	&
	& &Z_{N}(0) = Z_{0,N},\\
	&
	& & \dot{\xi}(t) = \alpha \xi(t) + \beta p(t), \\
	& 
	& & \xi(0) = \xi_0.
	\end{aligned}
	\tag{AP}\label{prob: finite dim approx integrated optimization problem}
\end{equation}

%% the following remark is not applicable because the mobility cost is no longer the LQR version of the actuators' motion.
% \begin{rem}
%     Problem \eqref{prob: new IOCA} is not an LQR because $B_N(x,t)$ is not linear in $x$, i.e., $B_N(x_1,t)+B_N(x_2,t) \neq B_N(x_1+x_2,t)$ for all $x_1,x_2 \in \mathbb{R}^{2m_a}$ and for all $t$.
% \end{rem}

Similar to \eqref{prob: new IOCA}, problem \eqref{prob: finite dim approx integrated optimization problem} can be turned into an equivalent form using LQR results for a finite-dimensional system:
\begin{equation}
	\begin{aligned}
	& \underset{p \in L^2([0,t_f];P)}{\text{minimize}}
	& & \innerproduct{Z_{0,N}}{\Pi_{N}(0) Z_{0,N}} + J_\text{m}(\xi,p) \\
	& \text{subject to} 
	& & \dot{\xi}(t) = \alpha \xi(t) + \beta p(t), \quad \xi(0) = \xi_0,
% 	& 
% 	& & \Pi_{N}(0) = S_N^T(t_f) Q_{f,N} S_N(t_f) + \int_0^{t_f} S^T(\tau) \left( Q_N(\tau) - \Pi_{N}(\tau) B_N^T (\tau) R^{-1} B_N(\tau) \Pi_{N}(\tau) \right) S_N(\tau) \dd \tau.
	\end{aligned}
	\tag{AP1}\label{prob: equivalent finite dim approx integrated optimization problem}
\end{equation}
where $\Pi_{N}(0)$ is defined in \eqref{eq: mild solution of approximate Riccati} with $t = 0$. Analogous to Theorems~\ref{thm: existence of a solution of IOCA} and \ref{thm: equivalence between P and P1}, the existence of a solution of \eqref{prob: equivalent finite dim approx integrated optimization problem} and how to use its solution to reconstruct a solution for \eqref{prob: finite dim approx integrated optimization problem} are stated in Theorem~\ref{thm: existence of a solution of approximated IOCA} below, whose proof is presented in Appendix~\ref{proof: existence of a solution of approximated IOCA}.

\begin{thm}\label{thm: existence of a solution of approximated IOCA}
    Consider problem \eqref{prob: equivalent finite dim approx integrated optimization problem} and let assumptions (A5)--(A8) and (A9)--(A11) hold. Then \eqref{prob: equivalent finite dim approx integrated optimization problem} has a solution, denoted by $p_N^*$. 
    Let $u_N^*$ be the optimal control obtained from \eqref{eq: approximated optimal lqr feedback control} with actuator trajectory steered by $p_N^*$. Then $u_N^*$ and $p_N^*$ minimize problem \eqref{prob: finite dim approx integrated optimization problem}.
    
    %This solution is also a solution of \eqref{prob: finite dim approx integrated optimization problem}, whose optimal control $u_N^*(\cdot)$ is given by \eqref{eq: new optimal LQR feedback control} with actuators steered by $p_N^*(\cdot)$.
\end{thm}

% \begin{pf*}{Proof}
%     See appendix~\ref{proof: existence of a solution of approximated IOCA}.
% \end{pf*}

An extension to Theorem~\ref{thm: existence of a solution of approximated IOCA} is that an optimal feedback control can be obtained from \eqref{eq: approximated optimal lqr feedback control} whenever the optimal guidance is solved from \eqref{prob: finite dim approx integrated optimization problem} or \eqref{prob: equivalent finite dim approx integrated optimization problem}. Basically, when the trajectory is determined via the optimal guidance, a feedback control can be implemented.

% Define the mapping $T:C([0,t_f];\mathbb{R}^m) \rightarrow C([0,t_f];\mathbb{R}^n)$ such that
%     \begin{equation}
%         (Tp)(t) \defeq e^{\alpha t}\xi_0 + \int_0^t e^{\alpha (t - \tau)} \beta p(\tau) \dd \tau = \xi(t). \label{eq: define trajectory state as a mapping}
%     \end{equation}
% The mapping $T$ is continuous because
%     \begin{align}
%         & |(Tp_1)(t) - (Tp_2)(t)|_1 \nonumber \\
%         %= & \int_0^t e^{\alpha (t-\tau)} \beta (p_1(\tau)-p_2(\tau)) \dd \tau \nonumber \\
%         \leq & \int_0^t |e^{\alpha (t - \tau)} \beta |_1 \dd \tau \sup_{\tau \in [0,t]} |p_1(\tau) - p_2(\tau)|_1,
%     \end{align}
%     which implies the existence of $d_2 >0$ such that
%     \begin{align}
%         \norm{Tp_1 - Tp_2}_{C([0,t_f];\mathbb{R}^n)} \leq  d_2 \norm{p_1 - p_2}_{C([0,t_f];\mathbb{R}^m)}. \label{eq: boundedness of trajectory norm}
%     \end{align}
% Furthermore, define the mapping $\bar{J}_{\text{m}}: C([0,t_f];\mathbb{R}^m) \rightarrow \mathbb{R}^+$ by 
% \begin{equation}
%     \bar{J}_{\text{m}}(p) \defeq J_{\text{m}}(Tp,p). \label{eq: define mobility cost as a mapping}
% \end{equation}

% Before we show that the solution of \eqref{prob: equivalent finite dim approx integrated optimization problem} converges to that of \eqref{prob: new equivalent IOCA}, we state assumption (A12) and Lemma~\ref{lemma: continuity of the total cost wrt guidance} below regarding the continuity of the mapping $\bar{J}_{\text{m}}$.
\revision{To establish convergence to the solution of \eqref{prob: new equivalent IOCA} of \eqref{prob: equivalent finite dim approx integrated optimization problem}'s solution, we need to restrict the set of admissible guidance to a smaller set as introduced below in assumption (A12).
}

\begin{enumerate}
    \item[(A12)] There exist $p_{\max}>0$ and $a_{\max}>0$ such that the set of admissible guidance is $\mathcal{P}(p_{\max},a_{\max}) \defeq \{p \in C([0,t_f];P): |p(t)|$ is uniformly bounded by $p_{\max}$ and $|p(t_1) - p(t_2)| \leq a_{\max} |t_1-t_2|, \ \forall t_1,t_2 \in [0,t_f] \}$.
\end{enumerate}

% Since the guidance functions defined in the set $\mathcal{P}(p_{\max},a_{\max})$ are uniformly equicontinuous and uniformly bounded, by the Arzel\`a-Ascoli Theorem \cite{royden2010real}, the set is sequentially compact.

\revision{There are two perspectives to interpreting the assumption (A12). Mathematically, (A12) requires the admissible guidance to be a continuous function that is uniformly bounded and uniformly equicontinuous. These two properties yield the sequential compactness of the set $\mathcal{P}(p_{\max},a_{\max})$ by the Arzel\`a-Ascoli Theorem \cite{royden2010real}. Practically, (A12) requires the input signal to be continuous and have bounds $p_{\max}$ and $a_{\max}$ on the magnitude and the rate of change, respectively. This requirement is reasonable and checkable because a continuous signal is commonly used for smooth operation, and the bounds on magnitude and changing rate are due to the physical limits of the motion of a platform.}
For example, in the case of single integrator dynamics where $p$ is the velocity command, $p_{\max}$ and $a_{\max}$ refer to the maximum speed and maximum acceleration, respectively.
\revision{Moreover, since time discretization of the signal is applied when computing the optimal guidance, as long as the bound $p_{\max}$ on the magnitude of the signal is determined, then the changing rate is bounded by $a_{\max} :=2p_{\max}/\Delta t_{\min}$ for the smallest discrete interval length $\Delta t_{\min}$.} Theorem~\ref{thm: convergence of approximate solution} below states the convergence of the approximate optimal solution with its proof in Appendix~\ref{prf: convergence of approximate solution}. 

% The parameters $p_{\max}$ and $a_{\max}$ can be determined by the vehicles carrying actuators. For example, in the case of single integrator dynamics where $p$ is the velocity command, $p_{\max}$ and $a_{\max}$ refer to the maximum speed and maximum acceleration, respectively.

% \begin{lem}\label{lemma: continuity of the total cost wrt guidance}
% Consider the mapping $\bar{J}_{\text{m}}$ defined as in \eqref{eq: define mobility cost as a mapping} and suppose assumptions (A9)--(A12) hold. Then the mapping $\bar{J}_{\text{m}}$ is continuous.
% \end{lem}
% \begin{pf*}{Proof}
%     See Appendix~\ref{proof: continuity of mobility cost wrt guidance}.
% \end{pf*}

\begin{thm}\label{thm: convergence of approximate solution}
    Consider problem \eqref{prob: new equivalent IOCA} and its finite-dimensional approximation \eqref{prob: equivalent finite dim approx integrated optimization problem}. Let assumptions (A4)--(A12) hold and let $p^*$ and $p_N^*$ denote the optimal guidance of \eqref{prob: new equivalent IOCA} and \eqref{prob: equivalent finite dim approx integrated optimization problem}, respectively. Then 
    \begin{equation}\label{eq: convergence of the approximated optimal cost}
        \lim_{ N \rightarrow \infty} |\optcosteval{\eqref{prob: equivalent finite dim approx integrated optimization problem}}{p_N^*} - \optcosteval{\eqref{prob: new equivalent IOCA}}{p^*}| = 0.
    \end{equation} 
    Furthermore, the cost function of \eqref{prob: new equivalent IOCA} evaluated at the guidance $p_N^*$ converges to the optimal cost of \eqref{prob: new equivalent IOCA}
    \begin{equation}\label{eq: convergence of the approximated optimal guidance}
        \lim_{ N \rightarrow \infty} |\costeval{\eqref{prob: new equivalent IOCA}}{p_N^*} - \optcosteval{\eqref{prob: new equivalent IOCA}}{p^*}| = 0.
    \end{equation} 
\end{thm}
% Assumptions (A4), (A5)-(A7), and (A9)-(A11) are needed for the existence of a solution of (P1).
% Assumptions (A8) and (A9)-(A11) are needed for the existence of a solution of (AP1).
% Assumptions (A12) is needed for the convergence of the optimal solutions.

% \begin{pf*}{Proof}
%     See Appendix~\ref{prf: convergence of approximate solution}.
% \qed 
% \end{pf*}

\begin{rem}
Two implications of Theorem~\ref{thm: convergence of approximate solution} follow. First, \eqref{eq: convergence of the approximated optimal cost} implies that the optimal cost of the approximated problem \eqref{prob: equivalent finite dim approx integrated optimization problem} converges to that of the exact problem \eqref{prob: new equivalent IOCA}, which justifies the approximation in \eqref{prob: finite dim approx integrated optimization problem}. Second, \eqref{eq: convergence of the approximated optimal guidance} implies that the approximate optimal guidance $p_N^*$, when evaluated by the cost function of \eqref{prob: new equivalent IOCA}, yields a cost that is arbitrarily close to the exact optimal cost of \eqref{prob: new equivalent IOCA}. Since $p_N^*$ is computable and $p^*$ is not, the convergence in \eqref{eq: convergence of the approximated optimal guidance} qualifies $p_N^*$ as an appropriate optimal guidance.
\end{rem}

% {\newContents
% \begin{remark}
%     Another way of solving \eqref{prob: integrated optimal control problem} is to construct a coupled PDE-ODE system by compositing the state variable $\mathcal{Z}$ with the state of the mobile actuators $\xi$ to form a composite state variable that belongs to the composite state space $\mathcal{H} \oplus \mathbb{R}^m$. However, the dynamics of the composite state become a control-affine system, because the control input operator $\mathcal{B}_{\xi}(\cdot)$ now contains $\xi(\cdot)$, which is part of the composite state. 
%     %To the authors' knowledge, there is no available result for solving an optimal control problem with infinite-dimensional control-affine dynamics.
% \end{remark}}

The convergence stated in Theorem~\ref{thm: convergence of approximate solution} is established based on several earlier stated results, including
\begin{enumerate}
    \item the input operator's continuity with respect to location (assumption (A4)), which leads to the continuity of the PDE cost with respect to actuator trajectory (Lemma~\ref{lemma: continuity of equivalent cost wrt actuator location});
    \item existence of the Riccati operator (Lemma~\ref{lemma: existence of Riccati mild solution}) and convergence of its approximation (Theorem~\ref{thm: convergence of riccati operator}); and
    \item sequential compactness of the set of admissible guidance (assumption (A12)), which leads to the continuity of the cost function with respect to guidance (Lemma~\ref{lemma: continuity of the total cost wrt guidance}).
\end{enumerate}
Notice that these key results, in an analogous manner, are also required in \cite{morris2010linear} when establishing the convergence to the exact optimal actuator locations of the approximate optimal locations \cite[Theorem~3.5]{morris2010linear}, i.e., 
\begin{enumerate}
    \item continuity with respect to location and compactness of the input operator \cite[Theorem~2.6]{morris2010linear}, which lead to continuity of the Riccati operator with respect to actuator locations \cite[Theorem~2.6]{morris2010linear};
    \item existence of the Riccati operator \cite[Theorem~2.3]{morris2010linear} and the convergence of its approximation \cite[Theorem~3.1]{morris2010linear}; and
    \item sequential compactness of the set of admissible locations, which is inherited from the setting that the spatial domain is closed and bounded in a finite-dimensional space.
\end{enumerate}
Although the establishment of convergence is similar to the one in \cite{morris2010linear}, the cost function and type of Riccati equation are different: we have the quadratic PDE cost plus generic mobility cost and differential Riccati equation in this paper for control and actuator guidance versus the Riccati operator's norm as cost function and algebraic Riccati equation in \cite{morris2010linear} for actuator placement. \revision{The similarity comes from the infinite-dimensional nature of PDEs such that approximation is necessary for computation, and convergence in approximation qualifies the approximate optimal solutions.}

\revision{
\subsection{Checking assumptions (A4)--(A12)}
For the approximated optimal guidance to be a good proxy of the exact optimal guidance, by Theorem~\ref{thm: convergence of approximate solution}, assumption (A4)--(A12) have to be checked to ensure the convergence. We summarize methods for checking these assumptions here. (A4): examine the explicit form of $\mathcal{B}$; (A5)--(A7): examine the explicit form of the operators $\mathcal{Q}$, $\mathcal{Q}_f$, and $\bar{\mathcal{B}}\bar{\mathcal{B}}^\star$ and their approximations; (A8): examine the explicit form of $B_N$; (A9)--(A11): examine the explicit form of $J_{\text{m}}$; (A12): examining the bounds on the magnitude and changing rate of the admissible guidance.
}

% \subsection{Galerkin approximation}\label{sec: Galerkin}
% We use the Galerkin scheme to approximate the infinite-dimensional variables. The orthonormal set of eigenfunctions of the Laplacian operator $\nabla^2$ (with zero Dirichlet boundary condition) over the spatial domain $\Omega = [0,1] \times [0,1]$ is $\phi_{i,j}(x,y) = 2 \sin(\pi i x) \sin(\pi j y)$. We introduce a single index $k \defeq (i-1)N + j$ such that $\phi_k \defeq \phi_{i,j}$. For brevity, we use $\mathcal{H}_N$ to denote the $N^2$-dimensional space spanned by the basis functions $\{\phi_k \}_{k=1}^{N^2}$. Recall the orthogonal projection $P_N: \mathcal{H} \rightarrow \mathcal{H}_N$. It follows that $P_N^{\star} = P_N$ and $P_N^{\star} P_N \rightarrow I$ strongly \cite{Burns2015Solutions}. Let $\Phi_{N} \defeq [\phi_1 \ \phi_2 \ \dots \ \phi_{N^2}]^\top$. The approximated terms can be represented by coordinates in the space spanned by $\{\phi_1 \ \phi_2 \ \dots \ \phi_{N^2} \}$:
% \begin{align*}
%     Z_N(t) & \defeq P_N \mathcal{Z}(t) =\Phi_{N}^\top \mathring{Z}_N(t)  \\
%     \Pi_N(t) & \defeq P_N \Pi(t) P_N = \Phi_{N}^\top \mathring{\Pi}_N(t) \Phi_{N} \\
%     A_N(t) & \defeq P_N \mathcal{A} P_N = \Phi_{N}^\top \mathring{A}_N \Phi_{N} \\
%     B_N(t) & \defeq P_N \mathcal{B}(t) = \Phi_{N}^\top \mathring{B}_N(t)  \\
%     Q_N(t) & \defeq P_N \mathcal{Q}(t) P_N = \Phi_{N}^\top \mathring{Q}_N(t) \Phi_{N} \\
%     Q_{f,N} & \defeq P_N \mathcal{Q}_f P_N = \Phi_{N}^\top \mathring{Q}_{f,N} \Phi_{N},
% \end{align*}
% where the notation $\mathring{ \ }$ indicates the coordinates.

\subsection{Gradient-descent for solving problem \eqref{prob: finite dim approx integrated optimization problem}} \label{sec: GD method}
% To find a local minimum of \eqref{prob: finite dim approx integrated optimization problem}, we use Pontryagin's minimum principle \cite{liberzon2011calculus}. 
Define the costates $\lambda(t) \in \mathcal{H}_N$ and ${\mu(t) \in \mathbb{R}^n}$ associated with $Z_N(t)$ and $\xi(t)$, respectively, for $t \in [0,t_f]$ and define the Hamiltonian:
\begin{align}
    & H(Z_N(t), \xi(t), u(t), p(t), \lambda(t), \mu(t)) \nonumber \\
    = & \ \innerproduct{Z_N(t)}{Q_N(t) Z_N(t)} + u^\top(t) R u(t) + h(\xi(t),t) \nonumber \\ 
    & \ + g(p(t),t)  + \lambda^\top(t) \left( A_N Z_N(t) + B_N(M\xi(t),t)u(t)  \right) \nonumber \\
    & \ + \mu^\top(t) \left( \alpha \xi(t) + \beta p(t) \right). \label{eq: hamiltonian for PMP}
\end{align}
By Pontryagin's minimum principle \cite{liberzon2011calculus}, we can solve a two-point boundary value problem originated from \eqref{eq: hamiltonian for PMP} to find a local minimum of \eqref{prob: finite dim approx integrated optimization problem}. The iterative procedure for solving the two-point boundary value problem can be implemented in a gradient-descent manner \cite{kirk2012optimal,lewis2012optimal}.

\section{Numerical examples}
We demonstrate the performance of the optimal guidance and control in two numerical examples. The first example uses the diffusion-advection process with zero Dirichlet boundary condition \eqref{eq: motivating dynamical system 1}--\eqref{eq: motivating dynamical system 3}. The second example uses the same process but with zero Neumann boundary condition.

\revision{
The examples are motivated by and simplified from practical applications, e.g., removal of harmful algal blooms (HAB). In this case, the distribution of the HAB's concentration on the water surface can be modeled by a 2D diffusion-advection process. The cases of zero Dirichlet and Neumann boundary conditions correspond  to the scenarios where the surface is circumvented by absorbent and nonabsorbent materials, respectively. The control to the process is implemented by the surface vehicles that use physical methods (e.g., emitting ultrasonic waves or hauling algae filters) or chemical methods (by releasing algal treatment) \cite{schroeder2018mitigating}, whose impact on the process can be characterized by the input operator \eqref{eq: example input operator}. The magnitude of the control determines how fast the concentration is reduced at the location of the actuator. The optimal control and guidance minimize the cost such that the HAB concentration is reduced while the vehicles do not exercise too much control nor conduct aggressive maneuvers. And vehicles' low-level control can track the optimal trajectories despite the model mismatch between the dynamics of the vehicles and those applied in the optimization problem \eqref{prob: new IOCA}.
}

We apply the following values in the numerical examples: $\Omega = [0,1] \times [0,1]$, $z_0(x,y) = 320(x-x^2)(y-y^2)$, $N = 13$, $m_a = 4$, $t_f = 1$, $\mathbf{v} = [0.1, -0.1]^\top$, $a = 0.05$, $U = \mathbb{R}^4$, $P_i = [-100,100]$, $p_{\max} = a_{\max} = 100$, $R = 0.1I_4$, $\mathcal{Q} = \mathcal{Q}_f = \chi(x,y)$, $h(\xi(t),t) = h_f(\xi(t_f)) =  0$, $g(p(t),t) = 0.1p^\top (t)p(t)$, %q(x,y) = 5 \mathds{1}(x=y),q_f(x,y) = \mathds{1}(x=y),\
$\xi_1(0) = [0.1,0.1]^\top$, $\xi_2(0) = [0.125,0.1]^\top$, $\xi_3(0) = [0.125,0.125]^\top$, $\xi_4(0) = [0.1,0.125]^\top$, $\sigma_i=0.05$, $\alpha_i = 0_{2 \times 2}$, and $\beta_i = I_2$
% \interdisplaylinepenalty=10000
% \begin{gather*}
%     \Omega = [0,1] \times [0,1], z_0(x,y) = 320(x-x^2)(y-y^2), \\
%     N = 13, m_a = 4, t_f = 1,\mathbf{v} = [0.1, -0.1]^\top\\
%     a = 0.05, b_i = 1,\sigma_i=0.05 ,\alpha_i = 0_{2 \times 2},\beta_i = I_2,\\
%     U = \mathbb{R}^m, P_i = [-100,100],p_{\max} = a_{\max} = 100,\\
%     R = 0.1I_4, \mathcal{Q} = \mathcal{Q}_f = \chi(x,y),\\
%     h(\xi(t),t) = h_f(\xi(t_f)) =  0, g(p(t),t) = 0.1p^\top (t)p(t),\\
%      %q(x,y) = 5 \mathds{1}(x=y),q_f(x,y) = \mathds{1}(x=y),\
%      \xi_1(0) = [0.1,0.1]^\top,\xi_2(0) = [0.125,0.1]^\top, \\
%       \xi_3(0) = [0.125,0.125]^\top,\xi_4(0) = [0.1,0.125]^\top,
% \end{gather*}
% \allowdisplaybreaks
for $i \in \{ 1,2,3,4\}$, where the indicator function $\chi(x,y) = 1$ if $x=y$, and $\chi(x,y) = 0$ if $x \neq y$. We use \eqref{eq: example input operator} for the input operator of each actuator. \revision{The P\'eclect number of the process is $|\mathbf{v}|_2/a \approx 2.83$, which implies neither the diffusion or the advection dominates the process. %(see the state evolution in Figs.~\ref{pic: Dirichlet PDE evolution} and \ref{pic: Neumann PDE evolution}).
}

\subsection{Diffusion-advection process with Dirichlet boundary condition}\label{sec: Dirichlet example}
We use the dynamics in \eqref{eq: motivating dynamical system 1}--\eqref{eq: motivating dynamical system 3} with the Dirichlet boundary condition. We use the Galerkin scheme to approximate the infinite-dimensional variables. The orthonormal set of eigenfunctions of the Laplacian operator $\nabla^2$ (with zero Dirichlet boundary condition) over the spatial domain $\Omega = [0,1] \times [0,1]$ is $\phi_{i,j}(x,y) = 2 \sin(\pi i x) \sin(\pi j y)$. We introduce a single index $k \defeq (i-1)N + j$ such that $\phi_k \defeq \phi_{i,j}$. For brevity, we use $\mathcal{H}_N$ to denote the $N^2$-dimensional space spanned by the basis functions $\{\phi_k \}_{k=1}^{N^2}$. Recall the orthogonal projection $P_N: \mathcal{H} \rightarrow \mathcal{H}_N$. It follows that $P_N^{\star} = P_N$ and $P_N^{\star} P_N \rightarrow I$ strongly \cite{Burns2015Solutions}. Let $\Phi_{N} \defeq [\phi_1 \ \phi_2 \ \dots \ \phi_{N^2}]^\top$. \revision{We choose $N=13$ because it is the smallest dimension such that the resulting optimal cost is within the 1\% of the optimal cost evaluated with the maximum dimension $N=20$ in the numerical studies (see Fig.~\ref{pic: convergence figure}).}

Assumption (A4) holds for the choice of input operator. With the Galerkin approximation using the orthonormal eigenfunctions of the Laplacian operator $\nabla^2$ with zero Dirichlet boundary condition, it can be shown that assumption (A8) holds for $l_N (\cdot) = N^2 l(\cdot)$. 
% The following derivation uses the notation for the 1D case for brevity (with a total of $N$ basis functions). But it also applies to 2D. Recall that 
% \begin{equation}
%     \int_{\Omega} \Phi_N (z) \mathcal{B}(\xi,t)(z) \dd z = \int_{\Omega} \Phi_N(z) \Phi_N^\top(z) \mathring{B}_N(t) \dd z,
% \end{equation}
% i.e., $\mathring{B}_N(t) = (\int_{\Omega} \Phi_N(z) \Phi_N^\top(z) \dd z)^{-1} \int_{\Omega} \Phi_N (z) \mathcal{B}(\xi,t)(z) \dd z  = \int_{\Omega} \Phi_N (z) \mathcal{B}(\xi,t)(z) \dd z$ (orthonormal basis yields the inverse term to be identity).
% \begin{align*}
%     & \norm{B_N(x,t) - B_N(y,t)}_{L^2(\Omega)} \\
%     = & \left( \int_{\Omega} (\mathring{B}_N^\top(x,t)(z) \Phi_N(z) - \mathring{B}_N^\top (y,t) \Phi_N(z)  )^2 \dd z \right)^{1/2} \\
%     \leq & \left( \int_{\Omega} |\mathring{B}_N(x,t) - \mathring{B}_N(y,t)|^2 |\Phi_N(z)|^2 \dd z \right)^{1/2} \\
%     = & N^{1/2} \left| \int_{\Omega} (\mathcal{B}(x,t)(z) - \mathcal{B}(y,t)(z))\Phi_N(z) \dd z   \right| \\
%     \leq & N^{1/2} \left(\int_{\Omega} (\mathcal{B}(x,t)(z) - \mathcal{B}(y,t)(z))^2 \dd z \right)^{1/2} \\ \nonumber
%     & \left( \int_{\Omega} |\Phi_N(z)|^2 \dd z \right)^{1/2} \\
%     = & N \norm{\mathcal{B}(x,t) - \mathcal{B}(y,t)}_{L^2(\Omega)} \leq N l(|x-y|_2).
% \end{align*}
Assumptions (A5)--(A7) hold with $q=1$ under the Galerkin approximation with aforementioned basis functions $\Phi_N$ \cite{Burns2015Solutions}. Assumptions (A9)--(A11) and (A12) hold for the choice of functions in the mobility cost and parameters of the set of admissible guidance, respectively. 

We use the forward-backward sweeping method \cite{mcasey2012convergence} to solve the two-point boundary value problem originated from the Hamiltonian \eqref{eq: hamiltonian for PMP}.
The forward propagation of $Z_N$ and $\xi$ and backward propagation of $\lambda$ and $\mu$ are computed using the Runge-Kutta method. The same method is also applied to propagate the approximate Riccati solution $\Pi(t)$. Spatial integrals are computed using Legendre-Gauss quadrature.
\revision{To verify the convergence of the approximate optimal cost $\optcosteval{\eqref{prob: equivalent finite dim approx integrated optimization problem}}{p_N^*}$ stated in \eqref{eq: convergence of the approximated optimal cost}, we compute $\optcosteval{\eqref{prob: equivalent finite dim approx integrated optimization problem}}{p_N^*}$ for $N \in \{6,7,\dots,20 \}$. Note that the total number of basis functions is $N^2$. The result is shown in Fig.~\ref{pic: convergence figure}, where exponential convergence can be observed.}

% \begin{figure*}[t]
% 	\centering
% 	\includegraphics[width=\textwidth]{Dirichlet_PDEstate.eps}
% 	\caption{Evolution of the diffusion-advection process with Dirichlet boundary condition under the optimal feedback control $\bar{u}^*$. The actuators are steered by the optimal guidance $p^*$. Snapshots at $t=0,\ 0.05,\ 0.1,$ and $0.2$ s show the transient stage, whereas those at $t=0.6$ and $1$ s show the relatively steady stage. The mobile disturbance is shown in gray.}
% 	\label{pic: Dirichlet PDE evolution}
% \end{figure*}

\begin{figure*}[t]
	\centering
	\includegraphics[width=\textwidth]{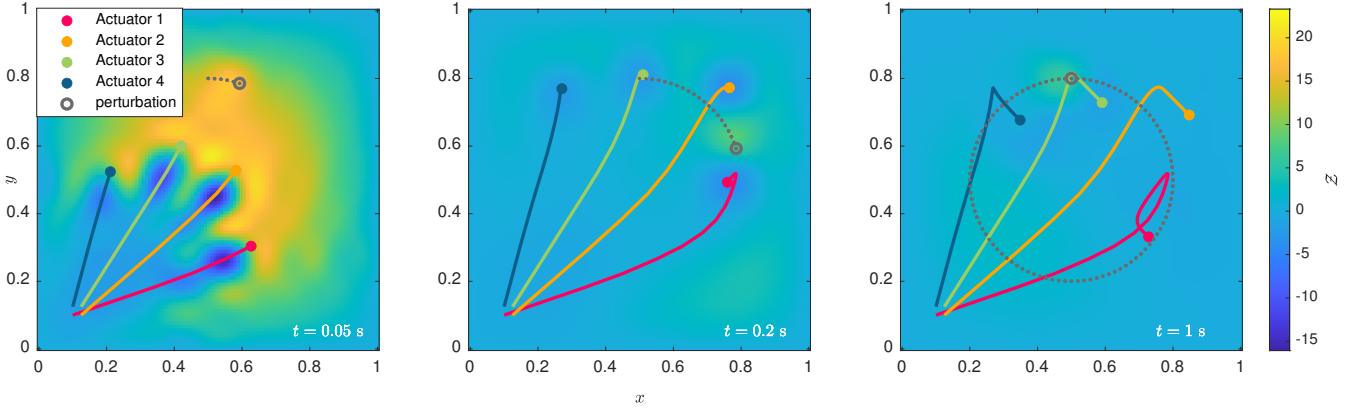}
	\caption{Evolution of the diffusion-advection process with Dirichlet boundary condition under the optimal feedback control $\bar{u}^*$. The actuators are steered by the optimal guidance $p^*$. Snapshots at $t=0.05$ and $0.2$ s show the transient stage, whereas the one at $t=1$ s shows the relatively steady stage. The mobile disturbance is shown by the gray circle.}
	\label{pic: Dirichlet PDE evolution}
\end{figure*}

\begin{figure}[t]
	\centering
	
    \includegraphics[width=\columnwidth]{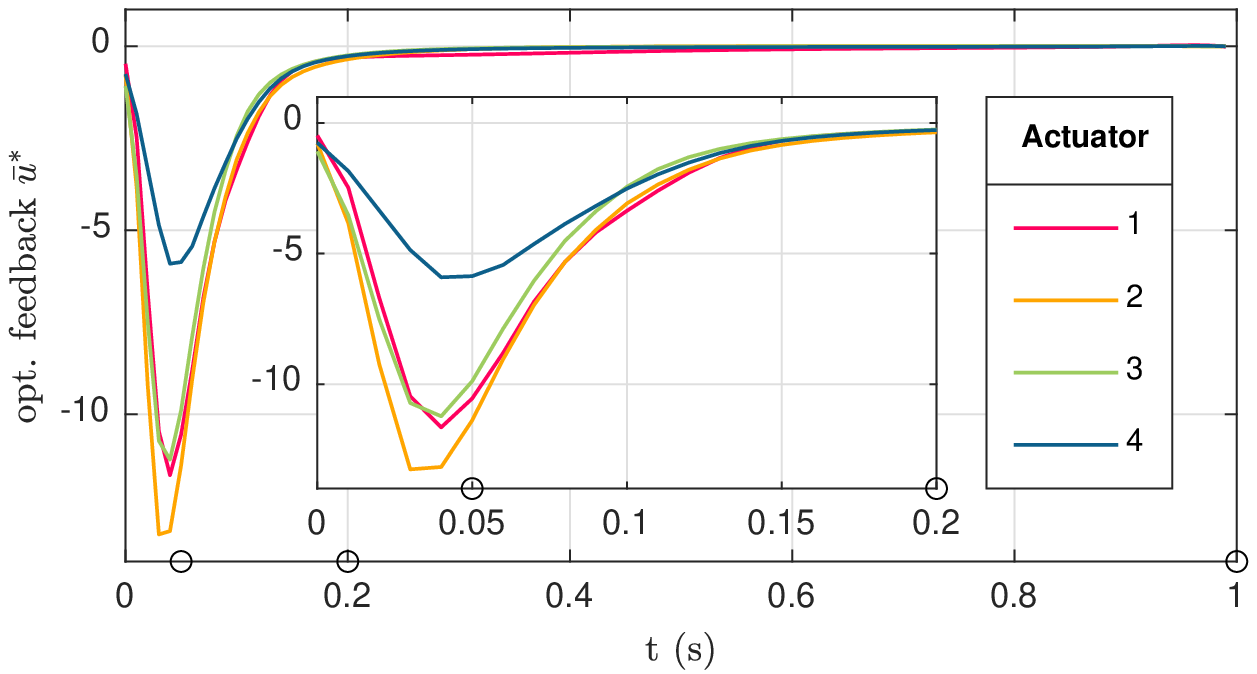}
	
	\caption{Optimal feedback control $\bar{u}^*$ of each actuator in the case of Dirichlet boundary condition. The circles along the horizontal axis correspond to the snapshots in Fig.~\ref{pic: Dirichlet PDE evolution}.}
	\label{pic: Dirichlet optimal feedback control}
\end{figure}

\begin{figure}[h]
	\centering
	
	\includegraphics[width=\columnwidth]{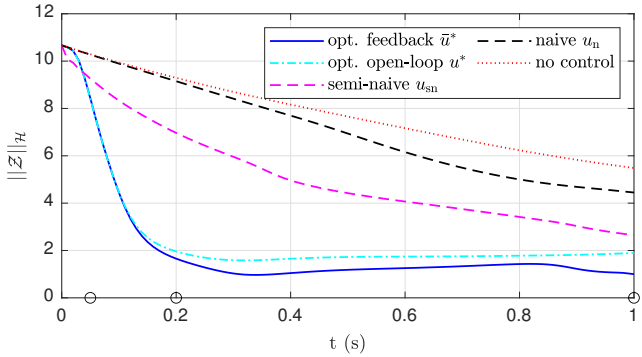}
	
	\caption{Norm of the PDE state in the case of Dirichlet boundary condition with pairs of control and guidance in Table~\ref{tb: the table for Dirichlet}. The circles along the horizontal axis correspond to the snapshots in Fig.~\ref{pic: Dirichlet PDE evolution}.}
	\label{pic: Dirichlet PDE norm figure}
\end{figure}

% It can be verified that assumptions (A4)--(A12) hold ({\color{red} expand it?}) with $q = 1$ under the Galerkin approximation scheme with basis functions using the orthonormal basis of operator $\mathcal{A}$. % Proposition 6.1 of Burns2015infinite covers (A5) and (A6).

In the simulation, a mobile disturbance $0.5\mathcal{B}(x_d(t),t)$, whose trajectory is $x_d(t) = [ 0.5 + 0.3 \sin(2 \pi t), 0.5 + 0.3 \cos(2 \pi t)]^\top$ is added to the right-hand side of the dynamics \eqref{eq: motivating dynamical system 1}. %The same type of disturbance has been applied in \cite{demetriou2012adaptive}. 

%We compare the costs of the optimal open-loop control~$u^*$, optimal state feedback control $\bar{u}^*$, semi-naive control $u_{sn}$ and naive control $u_n$. 
Denote the optimal open-loop control and optimal guidance solved using the gradient-descent method in Section~\ref{sec: GD method} by $u^*$ and $p^*$, respectively. The trajectory steered by $p^*$ is denoted by $\xi^*$. Recall that an optimal feedback control, denoted by $\bar{u}^*$, can be synthesized using \eqref{eq: approximated optimal lqr feedback control} based on the optimal trajectory $\xi^*$ of the actuators.

Fig.~\ref{pic: Dirichlet PDE evolution} shows the evolution of the process controlled by the optimal feedback control and the optimal trajectories of the actuators. The actuation concentrates in the first $0.2$~s, which is shown in Fig.~\ref{pic: Dirichlet optimal feedback control}. Meanwhile, the actuators quickly pass the peak of the initial PDE at the center of the spatial domain and spread evenly in space. Subsequently, the actuators 2--4 cease active steering and dispensing actuation. The flow field causes the actuators to drift until the terminal time.
% Revert the comment out of the following statement when we get chance to place 6 snapshots in Fig. 1.
% However, actuator 1 turns back and dispenses a small amount of actuation to get rid of the remaining peak from the upper flow from $0.2$ to $0.6$ s, before it ceases steering and dispensing actuation like the other actuators.

To demonstrate the performance of the optimal feedback control $\bar{u}^*$, we compare it with semi-naive control $u_{\text{sn}}$ and naive control $u_{\text{n}}$ defined as local feedback controls: $u_{\text{sn}}(t) = -0.1z_{\text{sn}}(\xi^*(t),t)$ and $u_{\text{n}}(t) =  - 0.1 z_{\text{n}}(\xi_{\text{n}}(t),t)$.
The semi-naive actuators follow the optimal trajectory~$\xi^*$, whereas the naive actuators follow the trajectory $\xi_{\text{n}}$, which moves at a constant speed from $\xi_0$ to $1_{n \times 1}-\xi_0$.
Table~\ref{tb: the table for Dirichlet} compares the cost breakdown of all the control and guidance strategies. 
The optimal feedback control yields a smaller cost than the optimal open-loop control due to the capability of feedback control in rejecting disturbances. Simulations with a disturbance-free model (not shown) yield identical total cost for optimal open-loop control and optimal feedback control, which justifies the correctness of the synthesis. Fig.~\ref{pic: Dirichlet PDE norm figure} compares the norm of the PDE state controlled by pairs of control and guidance listed in Table~\ref{tb: the table for Dirichlet}. %, closed-loop control $\bar{u}^*$ and 
As can be seen, the PDE is effectively regulated using optimal feedback control. %Fig.~\ref{pic: trajectory and actuation figure} shows the optimal actuation input $\bar{u}^*$ of each actuator. 
As a comparison, the norm associated with optimal open-loop control grows slowly after $0.3$~s due the influence of the disturbance, although its reduction in the beginning is indistinguishable from that of the optimal feedback control. 

\setlength{\tabcolsep}{3pt} % Default value: 6pt
\renewcommand{\arraystretch}{1} % Default value: 1
  \captionsetup{%size=footnotesize,
	%justification=centering, %% not needed
	skip=5pt, position = bottom}
\begin{table}[t]
	\centering
	\small
	%\captionsetup{font=small}
	\caption{Cost comparison of control and guidance strategies in the case of Dirichlet boundary condition. All costs are normalized with respect to the total cost of the case with no control.}
	\begin{tabular}{lcccccc}
		\toprule[1pt]
		\multicolumn{3}{c}{Control (C) and Guidance (G)}  & & \multicolumn{3}{c}{Cost}  \\
		\cmidrule{1-3}  \cmidrule{5-7} 
		& C & G & & $J_N$  &  $J_\text{m}$ & Total \\
		\midrule
        opt. feedback & $\bar{u}^*$ & $\xi^*$ & & 13.7\%	& 3.0\%	 & 16.7\% \\ 
opt. open-loop & $u^*$ & $\xi^*$ & & 17.5\% & 3.0\% & 20.5\% \\ 
semi-naive & $u_{\text{sn}}$ & $\xi^*$ & & 42.5\% & 3.0\% & 45.5\% \\ 
naive & $u_{\text{n}}$ &  $\xi_{\text{n}}$ & & 78.8\% & 0.5\% & 79.3\%  \\ 
no control & - &  - & & 100.0\% & 0.0\% & 100.0\%  \\ 
		\bottomrule[1pt]
	\end{tabular}\label{tb: the table for Dirichlet}
\end{table}
\normalsize

% Fig.~\ref{pic: PDE disturbance figure} shows the spatiotemporal distribution of the disturbed DPS with optimal feedback control $\bar{u}^*$ dispensed by actuators following the optimal trajectory $\xi^*$.

% {\newContents
% Two additional types of disturbances are considered. The first type is spatially distributed noise \cite{morris2015using}, $\mathcal{D} \nu(t)$, added to the right-hand side of \eqref{eq: diffusion equation}, where $\mathcal{D}$ denotes the spatial distribution (such as uniform distribution: $\mathcal{D}(x) = 1$ if $x \in \Omega$ and is $0$ otherwise) of the noise $\nu(t)$, which is a random variable satisfying a probability distribution (such as zero-mean Gaussian). The second type is parameter perturbation where the diffusion coefficient $a$ is perturbed to $a(1+\Delta)$ for $\Delta > -1$. If $\Delta < -1$, then a negative diffusion coefficient yields reverse diffusion, which is not considered in this paper. We observe that the optimal feedback control effectively regulates the DPS and performs better than the optimal open-loop control in simulations (not shown) that involves these two types of disturbances separately. }

\subsection{Diffusion-advection process with Neumann boundary condition}
The results derived in this paper also apply to the operator $\mathcal{A}$ defined in \eqref{eq: definition of diffusion-advection operator} with a Neumann boundary condition (BC), because a general second-order and uniformly elliptic operator with Neumann BC yields a strongly continuous analytic semigroup on $L^2(\Omega)$ \cite{lasiecka2000control}. % see the first paragraph on p. 197.
% See Definition 4.6 of https://www.math.ucdavis.edu/~hunter/pdes/ch4.pdf for uniformly elliptic operator
In this example, we consider the diffusion-advection process \eqref{eq: motivating dynamical system 1} with initial condition \eqref{eq: motivating dynamical system 3} and zero Neumann BC: ${\partial z (x,y,t)}/{ \partial \mathbf{n}} = 0$,
% \begin{equation}\label{eq: Neumann BC}
%     \frac{\partial z (\cdot,\cdot,t)}{ \partial \mathbf{n}} |_{\partial \Omega} = 0,
% \end{equation}
where $\mathbf{n}$ is the normal to the boundary $\partial \Omega$ and $(x,y) \in \partial \Omega$. Notice that the basis functions applied for Galerkin approximation in this case are the eigenfunctions of the Laplacian with zero Neumann BC, $\phi_{i,j}(x,y) = 2 \cos(\pi i x) \cos(\pi j y)$ for $i,j \in \{0,1,\dots\}$. All the parameters, disturbance, and pairs of control and guidance for comparison applied in this example are identical to those in Section~\ref{sec: Dirichlet example}. \revision{Exponential convergence in the approximate optimal cost can be observed in Fig.~\ref{pic: convergence figure}.}

\begin{figure*}[t]
	\centering
	\includegraphics[width=\textwidth]{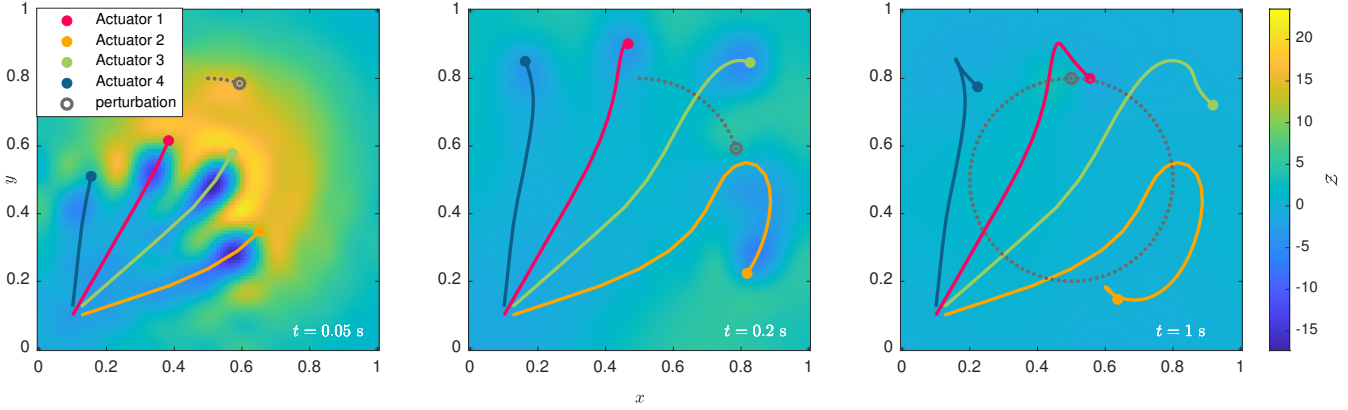}
	\caption{Evolution of the process with Neumann boundary condition under the optimal feedback control $\bar{u}^*$. The actuators are steered by the optimal guidance $p^*$. Snapshots at $t=0.05$ and $0.2$ s show the transient stage, whereas the one at $t=1$~s shows the relatively steady stage. The mobile disturbance is shown by the gray circle.}
	\label{pic: Neumann PDE evolution}
\end{figure*}

\begin{figure}[h]
	\centering
    \includegraphics[width=\columnwidth]{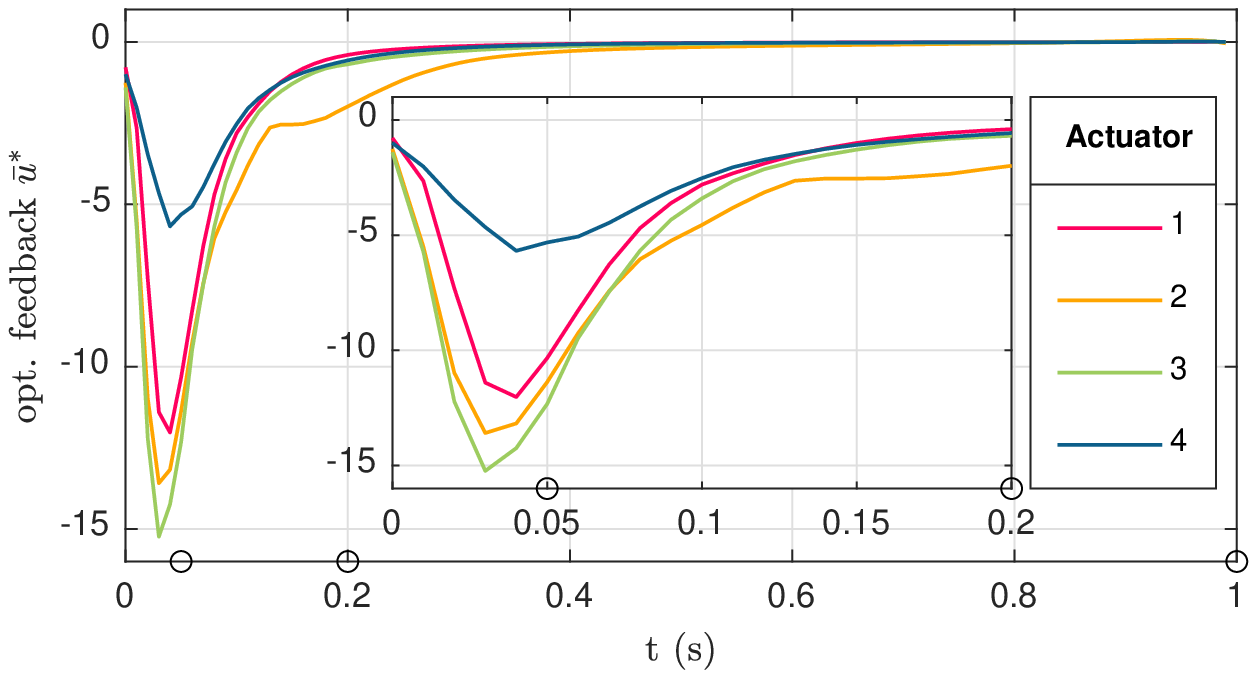}
	\caption{Optimal feedback control $\bar{u}^*$ of each actuator in the case of Neumann boundary condition. The circles along the horizontal axis correspond to the snapshots in Fig.~\ref{pic: Neumann PDE evolution}.}
	\label{pic: Neumann optimal feedback control}
\end{figure}

\begin{figure}[h]
	\centering
	\includegraphics[width=\columnwidth]{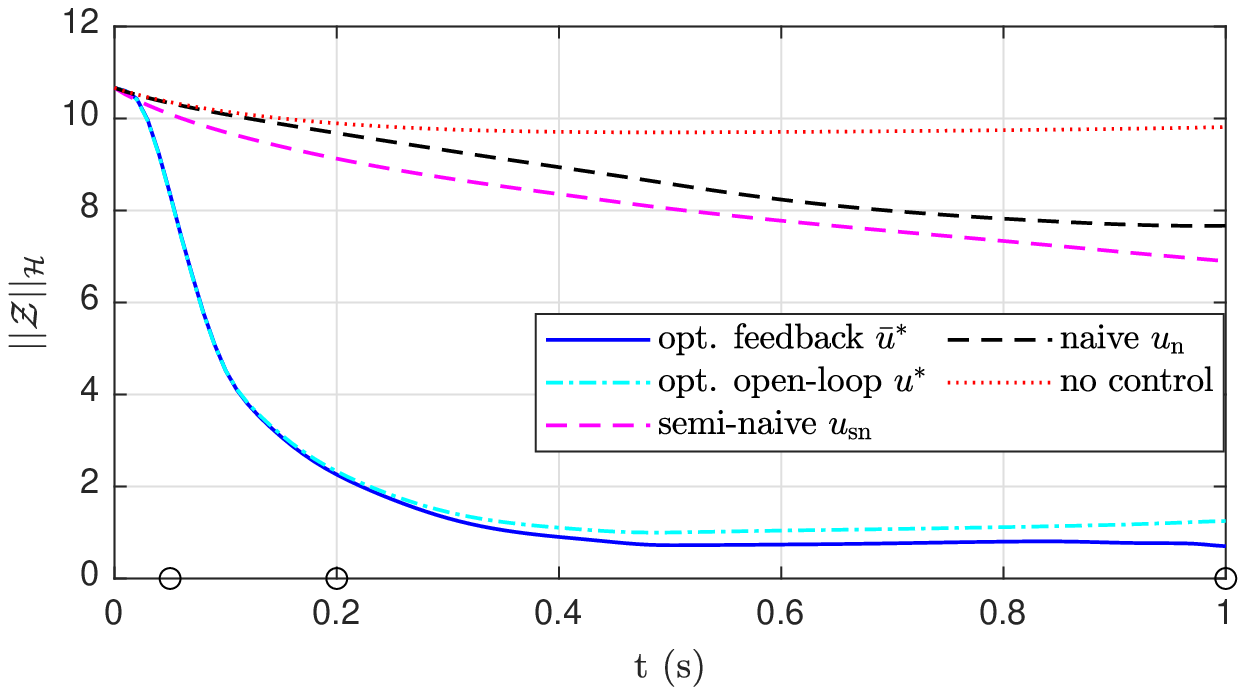}
	\caption{Norm of the PDE state in the case of Neumann boundary condition with pairs of control and guidance in Table~\ref{tb: the table for Neumann}. The circles along the horizontal axis correspond to the snapshots in Fig.~\ref{pic: Neumann PDE evolution}.}
	\label{pic: Neumann PDE norm figure}
\end{figure}

Fig.~\ref{pic: Neumann PDE evolution} shows the evolution of the process and the optimal trajectory of the actuators. Similar to the case of Dirichlet BC, the actuators spread out to cover most of the domain in the initial $0.2$ s, with most of the actuation implemented during the same interval, seen in Fig.~\ref{pic: Neumann optimal feedback control}. However, the actuators span a slightly larger area (Fig.~\ref{pic: Neumann PDE evolution}) and the maximum amplitude of actuation is bigger (Fig.~\ref{pic: Neumann optimal feedback control}), compared to the case of Dirichlet BC in Fig.~\ref{pic: Dirichlet PDE evolution} and Fig.~\ref{pic: Dirichlet optimal feedback control}, respectively. The difference is a consequence of the fact that the zero Neumann BC does not contribute to the regulation of the process because it insulates the process from the outside. Contrarily, the zero Dirichlet BC acts as a passive control that can essentially regulate the process to a zero state when there is no inhomoegeneous term in the dynamics \eqref{eq: motivating dynamical system 1}. This difference can be observed when comparing the norm of the PDE state in Fig.~\ref{pic: Neumann PDE norm figure} with Fig.~\ref{pic: Dirichlet PDE norm figure}. The norm of the uncontrolled state reduces slightly in the case of Neumann BC (Fig.~\ref{pic: Neumann PDE norm figure}) compared to the almost linear reduction in the case of Dirichlet BC (Fig.~\ref{pic: Dirichlet PDE norm figure}). Fig.~\ref{pic: Neumann PDE norm figure} also shows the difference of norm reduction between the optimal feedback control and optimal open-loop control. Once again, the former yields a smaller terminal norm than the latter due to the feedback's capability of disturbance rejection. The cost breakdown of the pairs of control and guidance in comparison is shown in Table.~\ref{tb: the table for Neumann}.

% \begin{figure*}[t]
% 	\centering
% 	\includegraphics[width=\textwidth]{Neumann_stateEvolution.eps}
% 	\caption{Evolution of the process with Neumann boundary condition under the optimal feedback control $\bar{u}^*$. The actuators are steered by the optimal guidance $p^*$. Snapshots at $t=0,\ 0.05,\ 0.1,$ and $0.2$ s show the transient stage, whereas those at $t=0.6$ and $1$ s show the relatively steady stage. The mobile disturbance is shown in gray.}
% 	\label{pic: Neumann PDE evolution}
% \end{figure*}

\setlength{\tabcolsep}{3pt} % Default value: 6pt
\renewcommand{\arraystretch}{1} % Default value: 1
  \captionsetup{%size=footnotesize,
	%justification=centering, %% not needed
	skip=5pt, position = bottom}
\begin{table}[t]
	\centering
	\small
	%\captionsetup{font=small}
	\caption{Cost comparison of control and guidance strategies in the case of Neumann boundary condition. All costs are normalized with respect to the total cost of the case with no control.}
	\begin{tabular}{lcccccc}
		\toprule[1pt]
		\multicolumn{3}{c}{Control (C) and Guidance (G)}  & & \multicolumn{3}{c}{Cost}  \\
		\cmidrule{1-3}  \cmidrule{5-7} 
		& C & G & & $J_N$  &  $J_\text{m}$ & Total \\
		\midrule
opt. feedback & $\bar{u}^*$ & $\xi^*$ & & 6.4\%	& 1.6\%	 & 8.0\% \\ 
opt. open-loop & $u^*$ & $\xi^*$ & & 7.1\% & 1.6\% & 8.7\% \\ 
semi-naive & $u_{\text{sn}}$ & $\xi^*$ & & 63.7\% & 1.6\% & 65.3\% \\ 
naive & $u_{\text{n}}$ &  $\xi_{\text{n}}$ & & 65.9\% & 0.2\% & 66.1\%  \\ 
no control & - &  - & & 100.0\% & 0.0\% & 100.0\%  \\ 
		\bottomrule[1pt]
	\end{tabular}\label{tb: the table for Neumann}
\end{table}
\normalsize

\revision{
% \subsection{Convergence of the approximate optimal cost}
% To verify the convergence of the approximate optimal cost $\optcosteval{\eqref{prob: equivalent finite dim approx integrated optimization problem}}{p_N^*}$ stated in \eqref{eq: convergence of the approximated optimal cost}, we compute $\optcosteval{\eqref{prob: equivalent finite dim approx integrated optimization problem}}{p_N^*}$ for $N \in \{6,7,\dots,20 \}$. Note that the total number of basis functions is $N^2$. The result is shown in Fig.~\ref{pic: convergence figure}, where exponential convergence for both the Dirichlet BC \eqref{eq: motivating dynamical system 2} and Neumann BC \eqref{eq: Neumann BC} can be observed. 

\begin{figure}[h]
	\centering
	\includegraphics[width=\columnwidth]{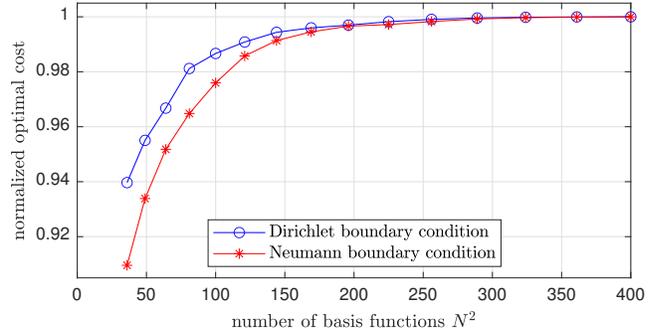}
	\caption{Approximate optimal costs $\optcosteval{\eqref{prob: equivalent finite dim approx integrated optimization problem}}{p_N^*}$ normalized with respect to the optimal cost for $N^2=400$.}
	\label{pic: convergence figure}
\end{figure}
}

\section{Conclusion}
This paper proposes an optimization framework that steers a team of mobile actuators to control a DPS modeled by a 2D diffusion-advection process. Specifically, jointly optimal control of the DPS and guidance of the mobile actuators are solved such that the sum of a quadratic PDE cost and a generic mobility cost is minimized subject to the dynamics of the DPS and of the mobile actuators. We obtain an equivalent problem using LQR of an abstract linear system, which reduces the problem to search for optimal guidance only. The optimal control can be synthesized once the optimal guidance is obtained. Conditions on the existence of a solution are established based on the equivalent problem. We use the Galerkin approximation scheme to reduce the problem to a finite-dimensional one and apply a gradient-descent method to compute optimal guidance and control numerically. We prove conditions under which the approximate optimal guidance converges to that of the exact optimal guidance \revision{in the sense that when evaluating these two solutions by the original cost function, the difference becomes arbitrarily small as the dimension of approximation increases}. The convergence justifies the appropriateness of both the approximate problem and its solution. The performance of the proposed optimal control and guidance is illustrated with two numerical examples\revision{, where exponential convergence of the approximate optimal cost is observed.}

Ongoing and future work includes \revision{establishing the convergence rate of the approximate optimal cost and} studying problems with other types of PDE cost, such as the operator norm of the Riccati operator \cite{morris2010linear,Cheng2021optimalWunknown} to characterize unknown initial conditions and H$_2$- or H$_{\infty}$-performance criteria for different types of perturbation \cite{morris2015using,kasinathan2013Hinfinity}.
% \commentout{The continuity of the PDE costs mentioned above can be established under suitable construction of the space of Riccati operator and the critical assumption that the input operator is continuous with respect to location.} 
% The study of a dual estimation framework is undergoing. 
On the actuator side, decentralized guidance design may be incorporated in future work to enable more autonomy of the team than a centralized implementation. \revision{Actuators that travel along the boundary may be considered as well which would result in boundary controller design.}
%Generally, methods for numerical integration include the Simpson's rule and the Gaussian quadrature. Because the above integrands all have explicit expressions, we will save the numerical integration for later.

% \section*{Acknowledgments} The authors would like to thank Artur Wolek and Debdipta Goswami for helpful suggestions. %Work reported in this article was sponsored under DARPA Contract No. HR001118C0142.

\appendix

\section{Proof of Theorem~\ref{thm: existence of a solution of IOCA}}\label{prf: existence of a solution of IOCA}

The proof uses Theorem~\ref{thm: optimization existence of solution} (stated below) to establish the existence of an optimal solution of \eqref{prob: new equivalent IOCA}.

% \begin{defn}\label{def: applied for existence of optimization solution}
%     Suppose $(X,\norm{\cdot})$ is a normed linear space. 
%     \begin{enumerate}
%         \item A sequence $\{x_k \} \subset X$ is weakly convergent to an $x \in X$, denoted by $x_k \rightharpoonup x$, if $\lim \limits_{k \rightarrow \infty} (x^{\star},x_k) = (x^{\star},x)$ for all $x^{\star}$ belonging to the dual space $X^{\star}$. 
%         \item A subset $A \subset X$ is weakly sequentially closed if $\{x_k \} \subset A$ and $x_k \rightharpoonup x$ implies $x \in A$.
%         \item A subset $A \subset X$ is weakly sequentially compact if for every sequence $\{x_k \} \subset A$ there exists a subsequence $\{x_{k_i} \} \subset \{x_k \}$ and an $x \in A$ with $x_{k_i} \rightharpoonup x$.
%         \item Suppose $A \subset X$ and $f:A \rightarrow \mathbb{R}$. The mapping $f$ is weakly sequentially lower semicontinuous on $A$ if $\{ x_k\} \subset A$ and $x_k \rightharpoonup x \in A$ implies $f(x) \leq \liminf_{k \rightarrow \infty} f(x_k)$.
%     \end{enumerate}
% \end{defn}

\begin{thm}{\cite[Theorem 6.1.4]{werner2013optimization}}\label{thm: optimization existence of solution}
    Suppose $(X,\norm{\cdot})$ is a normed linear space, $M_0 \subset X$ is weakly sequentially compact and $f : M_0 \rightarrow \mathbb{R}$ is weakly sequentially lower semicontinuous on $M_0$. Then there exists an $\bar{x} \in M_0$ such that $f(\bar{X}) = \inf \{f(x): x \in M_0 \}$.
\end{thm}

\begin{pf*}{Proof of Theorem~\ref{thm: existence of a solution of IOCA}}
%The proof is based on \cite[p.~219]{werner2013optimization}. 
Without loss of generality, we consider the case of one mobile actuator, i.e., $m_a = 1$. The case of $m_a \geq 2$ follows naturally. 

% First, we define the norm on $L^2([0,t];\mathbb{R}^m)$. For $p \in L^2([0,t];\mathbb{R}^m)$,
% \begin{equation}
%     \norm{p}_{L^2([0,t];\mathbb{R}^m)} =  \left( \int_0^{t} |p(\tau)|^2_2 \dd \tau \right)^{1/2}.
% \end{equation}
% The space $L^2([0,t];\mathbb{R}^m)$ equipped with norm $\norm{\cdot}_{L^2([0,t];\mathbb{R}^m)}$ is a reflexive Banach space for $t \in [0,t_f]$ \cite{yosida1988functional}. % see Example 3 on p. 115.

We want to apply Theorem~\ref{thm: optimization existence of solution} to prove that the minimum of the cost function of \eqref{prob: new equivalent IOCA} is achieved on a subset $\mathcal{P}_0$ (defined below) of the admissible set in which the cost of guidance is upper bounded. Consider problem \eqref{prob: new equivalent IOCA}'s admissible set of guidance functions $\mathcal{P} \defeq \{p \in L^2([0,t_f];\mathbb{R}^m): p(t) \in P, t \in [0,t_f] \}$. Assume there exists $p_0 \in \mathcal{P}$ such that $\costeval{\eqref{prob: new equivalent IOCA}}{p_0} < \infty$ and let $\mathcal{P}_0 \defeq \{p \in \mathcal{P}: \ \costeval{\eqref{prob: new equivalent IOCA}}{p} \leq \costeval{\eqref{prob: new equivalent IOCA}}{p_0} \}$. We wish to prove \underline{Condition-1}, \underline{Condition-2}, and \underline{Condition-3} stated below:

    \vspace{0.1cm}
    \noindent \revision{\underline{Condition-1}: The set $\mathcal{P}_0$ is bounded.}
    
    \noindent \revision{\underline{Condition-2}: The set $\mathcal{P}_0$ is weakly sequentially closed.}
    
    \noindent \underline{Condition-3}: The mapping $\costeval{\eqref{prob: new equivalent IOCA}}{\cdot}: \mathcal{P} \rightarrow \mathbb{R}$ is weakly sequentially lower semicontinuous on $\mathcal{P}_0$.
    
\revision{\underline{Condition-1} and \underline{Condition-2} imply that $\mathcal{P}_0$ is weakly sequentially compact. By Theorem~\ref{thm: optimization existence of solution}, problem \eqref{prob: new equivalent IOCA} has a solution when \underline{Condition-1}--\underline{Condition-3} hold.}

Before proving these three conditions, we define a mapping $T: L^2([0,t_f];\mathbb{R}^m) \rightarrow C ([0,t_f];\mathbb{R}^n)$ by $(Tp)(t) \defeq \xi(t) = e^{\alpha t} \xi_0 + \int_0^t e^{\alpha (t-\tau)} \beta p(\tau) \dd \tau $ for $t \in [0,t_f]$. For $p_1,p_2 \in L^2([0,t_f];\mathbb{R}^m)$ and $t \in [0,t_f]$, we have
\begin{align}
    & |Tp_1(t) - Tp_2(t)|_1 \nonumber \\
     \leq & \textstyle \int_0^{t} |e^{\alpha (t-\tau)} \beta|_1 |p_1(\tau) - p_2(\tau)|_1 \dd \tau \nonumber\\
    \leq & c_5  \textstyle \int_0^{t} |e^{\alpha (t-\tau)} \beta|_1 |p_1(\tau) - p_2(\tau)|_2 \dd \tau \nonumber\\
    % \leq & c_5 \left( \int_0^{t} |e^{\alpha (t-\tau)} \beta|_1^2 \dd \tau \right)^{1/2} \left(\int_0^t |p_1(\tau) - p_2(\tau)|_2^2 \dd \tau \right)^{1/2} \nonumber \\
    = & c_5 (  \textstyle \int_0^{t} |e^{\alpha (t-\tau)} \beta|_1^2 \dd \tau )^{1/2} \norm{p_1-p_2}_{L^2([0,t];\mathbb{R}^m)} \nonumber \\
    \leq & c_5 c_6 \norm{p_1-p_2}_{L^2([0,t];\mathbb{R}^m)}
\end{align}
for  $c_5$ and $c_6 > 0$. Hence, $\norm{Tp_1 - Tp_2}_{C([0,t_f];\mathbb{R}^n)} = \sup_{t \in [0,t_f]} |Tp_1(t) - Tp_2(t)|_1 \leq  c_5 c_6 \norm{p_1-p_2}_{L^2([0,t_f];\mathbb{R}^m)}$, which also shows that $T$ is a continuous mapping, i.e., $\norm{Tp}_{C([0,t_f];\mathbb{R}^n)} \leq c_5 c_6 \norm{p}_{L^2([0,t_f];\mathbb{R}^m)}$ for all $p \in L^2([0,t_f];\mathbb{R}^m)$.

\vspace{0.1cm}
\noindent \textit{Proof of \underline{Condition-1}}: Suppose $p \in \mathcal{P}_0$, then
\begin{align}
    \costeval{\eqref{prob: new equivalent IOCA}}{p_0} \geq & \ \costeval{\eqref{prob: new equivalent IOCA}}{p} \nonumber \\
    = & \ h_f(Tp(t_f)) +  \textstyle \int_0^{t_f} h(Tp(t),t) + g(p(t),t) \dd t\nonumber \\
    &  + \innerproduct{\mathcal{Z}_0}{\Pi(0) \mathcal{Z}_0}    \nonumber \\
    \geq & \  \textstyle \int_0^{t_f} d_1 |p(t)|^2_2 \dd t \nonumber \\
    = & \ d_1 \norm{p}_{L^2([0,t_f];\mathbb{R}^m)}^2,
\end{align}
where the second inequality follows from the nonnegativity of $h_f(\cdot), h(\cdot,\cdot)$, and $\innerproduct{\mathcal{Z}_0}{\Pi(0) \mathcal{Z}_0}$.
Since $d_1>0$, the boundedness of $\mathcal{P}_0$ follows.

\noindent \textit{Proof of \underline{Condition-2}}:
Suppose $\{p_k \} \subset \mathcal{P}_0$ and $\{p_k\}$ converges weakly to $p$ (denoted by $p_k \rightharpoonup p$). We want to show $p \in \mathcal{P}_0$. We start with proving that $\mathcal{P}$ is weakly sequentially closed and, hence, $p \in \mathcal{P}$. Subsequently, we show $\costeval{\eqref{prob: new equivalent IOCA}}{p} \leq \costeval{\eqref{prob: new equivalent IOCA}}{p_0}$ to conclude \underline{Condition-2}.

\revision{To show that the set $\mathcal{P}$ is weakly sequentially closed, by \cite[Theorem~2.11]{troltzsch2010optimal}, %see google books:
% Theorem~2.11: Every convex and closed subset of a Banach space is weakly sequentially closed. If the space is reflexive and the set is in addition bounded, then it is weakly sequentially compact.
it suffices to show that $\mathcal{P}$ is closed and convex. Let $\{q_k\} \subset \mathcal{P}$ and $q_k \rightarrow q$. We want to show $q \in \mathcal{P}$, i.e., $q \in L^2([0,t_f],\mathbb{R}^m)$ and $q(t) \in P$ for $t \in [0,t_f]$. Since $L^2([0,t_f],\mathbb{R}^m)$ is complete, we can choose a subsequence $\{q_{k_j} \} \subset \mathcal{P}$ that converges to $q$ pointwise almost everywhere on $[0,t_f]$ \cite[p.~53]{yosida1988functional}. Since $P$ is closed (assumption (A9)), $q(t) \in P$ for almost all $t \in [0,t_f]$. Hence, $\mathcal{P}$ is closed. The convexity of $\mathcal{P}$ follows from that of $P$ (assumption (A9)), i.e., if $p_1,p_2 \in \mathcal{P}$, then $\lambda p_1 + (1- \lambda)p_2 \in L^2([0,t_f];\mathbb{R}^m)$ and $\lambda p_1(t) + (1- \lambda)p_2(t) \in P$ for $t \in [0,t_f]$ and $\lambda \in [0,1]$.}

What remain to be shown is $\costeval{\eqref{prob: new equivalent IOCA}}{p} \leq \costeval{\eqref{prob: new equivalent IOCA}}{p_0}$. Since $p_k \rightharpoonup p$, by definition, we have $Tp_k \rightarrow Tp$. We now show that the sequence $\{Tp_k\}$ contains a uniformly convergent subsequence in $C([0,t_f];\mathbb{R}^n)$. 
The sequence $\{Tp_k \} \subset C([0,t_f];\mathbb{R}^n)$ is uniformly bounded and uniformly equicontinuous for the following reasons: Since $\norm{Tp_k}_{C([0,t_f];\mathbb{R}^n)} \leq c_5 c_6 \norm{p_k}_{L^2([0,t_f];\mathbb{R}^m)}$, it follows that $\norm{Tp_k}_{C([0,t_f];\mathbb{R}^n)}$ is uniformly bounded, because $\{p_k\} \subset \mathcal{P}_0$ which is a bounded set. For $s,t \in [0,t_f]$, we have%there exists $c_7>0$ and $c_8>0$ such that 
\begin{align*}
    & \ |Tp_k(t) - Tp_k(s)|_1  \nonumber \\
    = &\ \left| \textstyle \int_s^t \alpha Tp_k(\tau) +  \beta  p_k(\tau) \dd \tau \right|_1 \\
    \leq &\ |t-s||\alpha|_1\norm{Tp_k}_{C([0,t_f];\mathbb{R}^n)} \nonumber \\
    &\ + |t-s|^{1/2} |\beta|_2 \norm{p_k}_{L^2([0,t_f];\mathbb{R}^m)} .
    % \leq &\ c_7 |t-s| + c_8 |t-s|^{1/2}.
\end{align*}
Since $\{\norm{p_k}_{L^2([0,t_f];\mathbb{R}^m)} \}$ and $\{ \norm{Tp_k}_{C([0,t_f];\mathbb{R}^n)} \}$ both are uniformly bounded for all $p_k \in \mathcal{P}_0$, $\{Tp_k \}$ is uniformly equicontinuous. By the Arzel\`a-Ascoli Theorem \cite{royden2010real}, there is a uniformly convergent subsequence $\{ Tp_{k_j}\} \subset \{Tp_k \}$.

Without loss of generality, we assume $p_k \rightharpoonup p$ and $Tp_k \rightarrow Tp$ uniformly on $[0,t_f]$, and $\costeval{\eqref{prob: new equivalent IOCA}}{p_k} \leq \costeval{\eqref{prob: new equivalent IOCA}}{p_0}$. We have $\costeval{\eqref{prob: new equivalent IOCA}}{p_0} - \costeval{\eqref{prob: new equivalent IOCA}}{p} 
    = \costeval{\eqref{prob: new equivalent IOCA}}{p_0} - \costeval{\eqref{prob: new equivalent IOCA}}{p_k} +     \costeval{\eqref{prob: new equivalent IOCA}}{p_k} - \costeval{\eqref{prob: new equivalent IOCA}}{p} 
    \geq  \costeval{\eqref{prob: new equivalent IOCA}}{p_k} - \costeval{\eqref{prob: new equivalent IOCA}}{p}$,
% \begin{align}
%     & \costeval{\eqref{prob: new equivalent IOCA}}{p_0} - \costeval{\eqref{prob: new equivalent IOCA}}{p} \nonumber \\
%     = & \costeval{\eqref{prob: new equivalent IOCA}}{p_0} - \costeval{\eqref{prob: new equivalent IOCA}}{p_k} +     \costeval{\eqref{prob: new equivalent IOCA}}{p_k} - \costeval{\eqref{prob: new equivalent IOCA}}{p} \nonumber \\
%     \geq & \costeval{\eqref{prob: new equivalent IOCA}}{p_k} - \costeval{\eqref{prob: new equivalent IOCA}}{p},
% \end{align}
by which, to show $\costeval{\eqref{prob: new equivalent IOCA}}{p} \leq \costeval{\eqref{prob: new equivalent IOCA}}{p_0}$, it suffices to show $\costeval{\eqref{prob: new equivalent IOCA}}{p} \leq \liminf_{k \rightarrow \infty} \costeval{\eqref{prob: new equivalent IOCA}}{p_k} $, which is to show
\begin{align}
    & h_f(Tp(t_f)) +  \textstyle \int_0^{t_f} h(Tp(t),t) + g(p(t),t) \dd t \nonumber \\
    & + \innerproduct{\mathcal{Z}_0}{\Pi(0) \mathcal{Z}_0} 
 \nonumber \\
    \leq & \liminf_{k \rightarrow \infty} h_f(Tp_k(t_f)) +  \textstyle \int_0^{t_f} h(Tp_k(t),t) + g(p_k(t),t) \dd t \nonumber \\
    & + \innerproduct{\mathcal{Z}_0}{\Pi^k(0) \mathcal{Z}_0}, \label{eq: inequality to be proved}
\end{align}
where $\Pi^k(0)$ is the solution of \eqref{eq: new mild solution of operator Riccati equation} associated with actuator state $Tp_k$. Since $\{Tp_k\}$ converges to $Tp$ uniformly on $[0,t_f]$, the continuity of $h_f(\cdot)$ implies
\begin{equation}
    h_f(Tp(t_f)) = \liminf_{k \rightarrow \infty} h_f(Tp_k(t_f)); \label{eq: intermediate step Fatou's lemma}
\end{equation}
Fatou's lemma \cite{royden2010real} implies
\begin{equation}
      \textstyle \int_0^{t_f} h(Tp(t),t) \dd t \leq \liminf_{k \rightarrow \infty}  \textstyle \int_0^{t_f} h(Tp_k(t),t) \dd t;
\end{equation}
and Lemma~\ref{lemma: continuity of equivalent cost wrt actuator location} implies
\begin{equation}
    \innerproduct{\mathcal{Z}_0}{\Pi(0) \mathcal{Z}_0}  = \liminf_{k \rightarrow \infty} \innerproduct{\mathcal{Z}_0}{\Pi^k(0) \mathcal{Z}_0}. \label{eq: intermediate equality on PDE cost}
\end{equation}
To prove \eqref{eq: inequality to be proved}, based on \eqref{eq: intermediate step Fatou's lemma}--\eqref{eq: intermediate equality on PDE cost}, it suffices to show $\textstyle \int_0^{t_f} g(p(t),t) \dd t \leq \liminf_{k \rightarrow \infty}  \textstyle \int_0^{t_f} g(p_k(t),t) \dd t$.
By contradiction, assume there is $\lambda > 0$ such that 
\begin{equation}
    \liminf_{k \rightarrow \infty}  \textstyle \int_o^{t_f} g(p_k(t),t) \dd t  < \lambda <  \textstyle \int_0^{t_f} g(p(t),t) \dd t. \label{eq: contradiction assumed for existence of solution}
\end{equation}
There exists a subsequence $\{p_{k_j} \} \subset \{p_k \}$ such that $O_{\lambda} \defeq \{q \in L^2([0,t_f];\mathbb{R}^{m}):  \textstyle \int_0^{t_f} g(q(t),t) \dd t \leq \lambda \} $
% \begin{equation*}
%     O_{\lambda} \defeq \{q \in L^2([0,t_f];\mathbb{R}^{m}):  \textstyle \int_0^{t_f} g(q(t),t) \dd t \leq \lambda \}   
% \end{equation*}
and $\{p_{k_j} \} \subset O_{\lambda} $.
We wish to show that $O_{\lambda}$ is weakly sequentially closed. By \cite[Theorem~2.11]{troltzsch2010optimal}, it suffices to show that $O_{\lambda}$ is convex and closed. Since $g(\cdot,t): \mathbb{R}^{m} \rightarrow \mathbb{R}$ is convex for all $t \in [0,t_f]$, it follows that $O_{\lambda}$ is convex. Let $\{q_k \} \subset O_{\lambda}$ and $\norm{q_k - q}_{L^2([0,t_f];\mathbb{R}^m)}$ converges to $0$ as $k \rightarrow \infty$. We can choose a subsequence $\{q_{k_j} \} \subset \{q_k  \}$ such that $q_{k_j}$ converges to $q$ pointwise almost everywhere on~$[0,t_f]$ \cite[p. 53]{yosida1988functional}. Now we have
\newline (1) $g(q_{k_j}(t),t) \geq 0$ for all $t \in [0,t_f]$ (assumption (A11));
\newline (2) $ \lim_{j \rightarrow \infty} g(q_{k_j}(t),t) = g(q(t),t)$ almost everywhere on $[0,t_f]$.
% \begin{enumerate}
%     \item $g(q_{k_j}(t),t) \geq 0$ for all $t \in [0,t_f]$ (assumption (A11));
%     \item $ \lim_{j \rightarrow \infty} g(q_{k_j}(t),t) = g(q(t),t)$ almost everywhere on $[0,t_f]$.
% \end{enumerate}

By Fatou's lemma \cite{royden2010real}, 
\begin{equation*}
     \textstyle \int_0^{t_f} g(q(t),t) \dd t  \leq \liminf_{k \rightarrow \infty}  \textstyle \int_0^{t_f} g(q_{k_j}(t),t) \dd t \leq \lambda,
\end{equation*}
where the last inequality holds due to $\{q_{k_j} \} \subset O_{\lambda}$. Hence, $q \in O_{\lambda}$ and $O_{\lambda}$ is closed.

Since $O_{\lambda}$ is weakly sequentially closed, $p_{k_j} \rightharpoonup p$ implies that $p \in O_{\lambda}$, which contradicts \eqref{eq: contradiction assumed for existence of solution}. Hence, $\costeval{\eqref{prob: new equivalent IOCA}}{p} \leq \costeval{\eqref{prob: new equivalent IOCA}}{p_0}$ is proved, and we conclude \underline{Condition-2}.

\noindent \textit{Proof of \underline{Condition-3}}: We now show that the mapping $\costeval{\eqref{prob: new equivalent IOCA}}{\cdot}: \mathcal{P} \rightarrow \mathbb{R}$ is weakly sequentially lower semicontinuous %\footnote{Suppose $(X,\norm{\cdot})$ is a normed linear space, $A \subset X$, and $f:A \rightarrow \mathbb{R}$. The function $f$ is weakly sequentially lower semicontinuous on $A$ if 
% \begin{equation}
%     \{x_k \} \subset A, \ x_k \rightharpoonup x \in A \Rightarrow f(x) \leq \liminf_{k \rightarrow \infty} f(x_k).
% \end{equation}
% } 
on $\mathcal{P}_0$. Suppose $\{p_k \} \subset \mathcal{P}_0$ and $p_k \rightharpoonup p \in \mathcal{P}_0$. We wish to show $\costeval{\eqref{prob: new equivalent IOCA}}{p} \leq \liminf_{k \rightarrow \infty} \costeval{\eqref{prob: new equivalent IOCA}}{p_k} $, which has been established when we proved $\costeval{\eqref{prob: new equivalent IOCA}}{p} \leq \costeval{\eqref{prob: new equivalent IOCA}}{p_0}$ in \underline{Condition-2} (starting from \eqref{eq: inequality to be proved}).

So we conclude that the existence of a solution of problem \eqref{prob: new equivalent IOCA}.
 \qed
\end{pf*}

\section{Proof of Theorem~\ref{thm: equivalence between P and P1}}\label{prf: equivalence between P and P1}

\begin{pf*}{Proof}
By contradiction, assume there are $p_0^*$ and $u_0^*$ minimizing \eqref{prob: new IOCA} and $p_0^* \neq p^*$ and $u_0^* \neq u^*$ such that $\optcosteval{\eqref{prob: new IOCA}}{u_0^*,p_0^*} < \costeval{\eqref{prob: new IOCA}}{u^*,p^*} = \costeval{\eqref{prob: new equivalent IOCA}}{p^*}$. Denote $\bar{u}_0^*$ the optimal control \eqref{eq: new optimal LQR feedback control} associated with actuator trajectory steered by $p_0^*$. It follows that $\optcosteval{\eqref{prob: new IOCA}}{u_0^*,p_0^*} = \costeval{\eqref{prob: new IOCA}}{\bar{u}_0^*,p_0^*}$, because {$\optcosteval{\eqref{prob: new IOCA}}{u_0^*,p_0^*} > \costeval{\eqref{prob: new IOCA}}{\bar{u}_0^*,p_0^*}$ violates the optimality of $u_0^*$} and {$\optcosteval{\eqref{prob: new IOCA}}{u_0^*,p_0^*} < \costeval{\eqref{prob: new IOCA}}{\bar{u}_0^*,p_0^*}$} contradicts the fact that $\bar{u}_0^*$ minimizes the quadratic cost $J(\mathcal{Z},u)$ (see Lemma~\ref{lemma: equivalent PDE cost}). Since $\costeval{\eqref{prob: new IOCA}}{\bar{u}_0^*,p_0^*} = \innerproduct{\mathcal{Z}_0}{\Pi_0^*(0) \mathcal{Z}_0} + J_{\text{m}}(\xi_0^*,p_0^*) =  \costeval{\eqref{prob: new equivalent IOCA}}{p_0^*} 
    <  \optcosteval{\eqref{prob: new equivalent IOCA}}{p^*}$,
% \begin{align}
%     \costeval{\eqref{prob: new IOCA}}{\bar{u}_0^*,p_0^*} = & \innerproduct{\mathcal{Z}_0}{\Pi_0^*(0) \mathcal{Z}_0} + J_{\text{m}}(\xi_0^*,p_0^*) \nonumber \\
%     = & \costeval{\eqref{prob: new equivalent IOCA}}{p_0^*} \nonumber \\
%     <&  \optcosteval{\eqref{prob: new equivalent IOCA}}{p^*},
% \end{align}
where $\Pi_0^*(0)$ associates with trajectory $\xi_0^*$ steered by $p_0^*$, it follows that $p^*$ is not an optimal solution of \eqref{prob: new equivalent IOCA}, which contradicts the optimality of $p^*$ for \eqref{prob: new equivalent IOCA}.\qed 
\end{pf*}

\section{Proof of Theorem~\ref{thm: existence of a solution of approximated IOCA}} \label{proof: existence of a solution of approximated IOCA}

\begin{pf*}{Proof}
    Since $\innerproduct{Z_{0,N}}{\Pi_{N}(0) Z_{0,N}} \geq 0$ and the mapping $K_N: C([0,t_f];\mathbb{R}^n) \rightarrow \mathbb{R}^+$ is continuous (see Lemma~\ref{lem: continuity of approximated equivalent cost wrt actuator location}), the proof is analogous to that of Theorem~\ref{thm: existence of a solution of IOCA}, where we use $\innerproduct{Z_{0,N}}{\Pi_{N}(0) Z_{0,N}}$ to substitute $\innerproduct{\mathcal{Z}_0}{\Pi(0) \mathcal{Z}_0} $. The proof that $u_N^*$ and $p_N^*$ minimize problem \eqref{prob: finite dim approx integrated optimization problem} follows from the same logic as the proof of Theorem~\ref{thm: equivalence between P and P1}. \qed
\end{pf*}

\section{Proof of Theorem~\ref{thm: convergence of approximate solution}} \label{prf: convergence of approximate solution}

\revision{
Before we prove Theorem~\ref{thm: convergence of approximate solution}, we first establish two intermediate results in Lemma~\ref{lemma: continuity of the total cost wrt guidance}, whose proof is in the supplementary material% \cite{Cheng2020optimalControlSupplement}
.}

% Consider the continuous mapping $T:L^2([0,t_f];\mathbb{R}^m) \rightarrow C([0,t_f];\mathbb{R}^n)$ as defined in the proof of Theorem~\ref{thm: existence of a solution of IOCA} such that $(Tp)(t) \defeq \xi(t) = e^{\alpha t} \xi_0 + \int_0^t e^{\alpha (t-\tau)} \beta p(\tau) \dd \tau $ for $t \in [0,t_f]$. Since the admissible guidance is in $C([0,t_f];\mathbb{R}^m) \subset L^2([0,t_f];\mathbb{R}^m)$, the continuity of $T$ is preserved, i.e., there exists $d_2 > 0$ such that for $p_1,p_2 \in C([0,t_f];\mathbb{R}^m)$
%     \begin{align}
%         \norm{Tp_1 - Tp_2}_{C([0,t_f];\mathbb{R}^n)} \leq  d_2 \norm{p_1 - p_2}_{C([0,t_f];\mathbb{R}^m)}. \label{eq: boundedness of trajectory norm}
%     \end{align}
% }

% Furthermore, define the mapping $\bar{J}_{\text{m}}: C([0,t_f];\mathbb{R}^m) \rightarrow \mathbb{R}^+$ by
% \begin{equation}
%     \bar{J}_{\text{m}}(p) \defeq J_{\text{m}}(Tp,p). \label{eq: define mobility cost as a mapping}
% \end{equation}
% {\color{red}The continuity of $\bar{J}_{\text{m}}$ is stated in Lemma~\ref{lemma: continuity of the total cost wrt guidance}, whose proof is in xxx.}

\revision{
\begin{lem}\label{lemma: continuity of the total cost wrt guidance}
Consider problem \eqref{prob: new equivalent IOCA} and its approximation \eqref{prob: equivalent finite dim approx integrated optimization problem}. If assumptions (A4)--(A7) and (A9)--(A12) hold, then the following two implications hold:
\newline 1. For $p \in C([0,t_f];P)$, $\lim_{N \rightarrow \infty} |\costeval{\eqref{prob: equivalent finite dim approx integrated optimization problem}}{p} - \costeval{\eqref{prob: new equivalent IOCA}}{p}| = 0$, where $N$ is the dimension of approximation applied in \eqref{prob: finite dim approx integrated optimization problem}.
\newline 2. The mapping $J_{\eqref{prob: new equivalent IOCA}}: C([0,t_f];P) \rightarrow \mathbb{R}^+$ is continuous, where $\costeval{\eqref{prob: new equivalent IOCA}}{p} = \innerproduct{\mathcal{Z}_0}{\Pi(0) \mathcal{Z}_0} + J_{\text{m}}(\xi,p)$. Here, the actuator state $\xi$ follows the dynamics \eqref{eq: general dynamics of the mobile actuator} steered by the guidance $p$, and $\Pi(0)$ follows \eqref{eq: new operator Riccati equation} with the actuator state~$\xi$.
\end{lem}
}

\begin{pf*}{Proof of Theorem~\ref{thm: convergence of approximate solution}}
    \revision{In the notation $\costeval{\eqref{prob: equivalent finite dim approx integrated optimization problem}}{p_N^*}$, the dimension of approximation in \eqref{prob: equivalent finite dim approx integrated optimization problem}, which is $N$ in this case, is indicated by its solution $p_N^*$. We append a subscript to indicate the dimension when it is not explicitly reflected by the argument, e.g., $\costeval{\eqref{prob: equivalent finite dim approx integrated optimization problem}_N}{p}$.}
    
    We first show \eqref{eq: convergence of the approximated optimal cost}, i.e., $|\optcosteval{\eqref{prob: equivalent finite dim approx integrated optimization problem}}{p_N^*} - \optcosteval{\eqref{prob: new equivalent IOCA}}{p^*}| \rightarrow 0$ as $N \rightarrow \infty$. First,
    \begin{align}
        \optcosteval{\eqref{prob: equivalent finite dim approx integrated optimization problem}}{p_N^*} = &\ \underset{p \in \mathcal{P}(p_{\max},a_{\max})}{\min} \costeval{\eqref{prob: equivalent finite dim approx integrated optimization problem}}{p} \nonumber \\
        \leq &\ \costeval{\eqref{prob: equivalent finite dim approx integrated optimization problem}}{p^*} \nonumber \\
        \leq &\ |\costeval{\eqref{prob: equivalent finite dim approx integrated optimization problem}}{p^*} - \optcosteval{\eqref{prob: new equivalent IOCA}}{p^*}| + \optcosteval{\eqref{prob: new equivalent IOCA}}{p^*}. \nonumber
    \end{align}
    Since $|\costeval{\eqref{prob: equivalent finite dim approx integrated optimization problem}}{p^*} - \optcosteval{\eqref{prob: new equivalent IOCA}}{p^*}| \rightarrow 0$ as $N \rightarrow 0$ (see Lemma~\ref{lemma: continuity of the total cost wrt guidance}-1), it follows that
    \begin{align}\label{eq: intermediate step 3 for solution convergence}
        \limsup_{N \rightarrow \infty} \optcosteval{\eqref{prob: equivalent finite dim approx integrated optimization problem}}{p_N^*} \leq &\ \optcosteval{\eqref{prob: new equivalent IOCA}}{p^*}.
    \end{align}

    To proceed with proving \eqref{eq: convergence of the approximated optimal cost}, in addition to \eqref{eq: intermediate step 3 for solution convergence}, we shall show $\liminf_{N \rightarrow \infty} \optcosteval{\eqref{prob: equivalent finite dim approx integrated optimization problem}}{p_N^*} \geq \optcosteval{\eqref{prob: new equivalent IOCA}}{p^*}$.
    Choose a convergent subsequence $\{\optcosteval{\eqref{prob: equivalent finite dim approx integrated optimization problem}}{p_{N_k}^*} \}_{k=1}^{\infty}$ such that $\lim_{k \rightarrow \infty} \optcosteval{\eqref{prob: equivalent finite dim approx integrated optimization problem}}{p_{N_k}^*} = \liminf_{N \rightarrow \infty} \optcosteval{\eqref{prob: equivalent finite dim approx integrated optimization problem}}{p_N^*}$.
    Since the guidance functions defined in the set $\mathcal{P}(p_{\max},a_{\max})$ are uniformly equicontinuous and uniformly bounded, by the Arzel\`a-Ascoli Theorem \cite{royden2010real}, there is a uniformly convergent subsequence of $\{ p_{N_k}^* \}_{k=1}^{\infty}$ which we use the same index $\{N_k\}_{k=1}^{\infty}$ to simplify notation and let the limit of $\{p_{N_k}^*\}_{k=1}^{\infty}$ be $p_{\inf}^*$, i.e.,
    \begin{equation}\label{eq: convergence of p_Nk to p_inf}
        \lim_{k \rightarrow \infty} \norm{p_{N_k}^* - p_{\inf}^*}_{C([0,t_f];\mathbb{R}^n)}  = 0.
    \end{equation}
    Now, $|\optcosteval{\eqref{prob: equivalent finite dim approx integrated optimization problem}}{p_{N_k}^*} - \costeval{\eqref{prob: new equivalent IOCA}}{p_{\inf}^*}| \leq | \optcosteval{\eqref{prob: equivalent finite dim approx integrated optimization problem}}{p_{N_k}^*} - \costeval{\eqref{prob: new equivalent IOCA}}{p_{N_k}^*}|+ | \costeval{\eqref{prob: new equivalent IOCA}}{p_{N_k}^*}- \costeval{\eqref{prob: new equivalent IOCA}}{p_{\inf}^*}|$,
    % \begin{align}
    %     & \ |\optcosteval{\eqref{prob: equivalent finite dim approx integrated optimization problem}}{p_{N_k}^*} - \costeval{\eqref{prob: new equivalent IOCA}}{p_{\inf}^*}| \nonumber \\
    %     \leq & \ | \optcosteval{\eqref{prob: equivalent finite dim approx integrated optimization problem}}{p_{N_k}^*} - \costeval{\eqref{prob: new equivalent IOCA}}{p_{N_k}^*}| \nonumber \\
    %     & \ + | \costeval{\eqref{prob: new equivalent IOCA}}{p_{N_k}^*}- \costeval{\eqref{prob: new equivalent IOCA}}{p_{\inf}^*}|,
    % \end{align}
    which implies
    \begin{align}
        & \ \limsup_{k \rightarrow \infty}|\optcosteval{\eqref{prob: equivalent finite dim approx integrated optimization problem}}{p_{N_k}^*} - \costeval{\eqref{prob: new equivalent IOCA}}{p_{\inf}^*}| \nonumber \\
        \leq & \ \lim_{k \rightarrow \infty}| \optcosteval{\eqref{prob: equivalent finite dim approx integrated optimization problem}}{p_{N_k}^*} - \costeval{\eqref{prob: new equivalent IOCA}}{p_{N_k}^*}| \nonumber \\
        & \ + \lim_{k \rightarrow \infty} | \costeval{\eqref{prob: new equivalent IOCA}}{p_{N_k}^*}- \costeval{\eqref{prob: new equivalent IOCA}}{p_{\inf}^*}|. \label{eq: first limit in the convergence of approximated optimal solution}
    \end{align}
    \revision{The first limit on the right-hand side of \eqref{eq: first limit in the convergence of approximated optimal solution} is zero for the following reason. For all $p \in \mathcal{P}(p_{\max},a_{\max})$, $\costeval{\eqref{prob: equivalent finite dim approx integrated optimization problem}_N}{p}$ converges to $\costeval{\eqref{prob: new equivalent IOCA}}{p}$ pointwise as the dimension of approximation $N\rightarrow \infty$ (see Lemma~\ref{lemma: continuity of the total cost wrt guidance}-1). Furthermore, since the sequence of approximated PDE cost $\{\innerproduct{Z_N(0)}{\Pi_N(0)Z_N(0)} \}_{N=1}^{\infty}$ is a monotonically increasing sequence, the sequence $\{ \costeval{\eqref{prob: equivalent finite dim approx integrated optimization problem}_N}{p}\}_{N=1}^{\infty}$ is a monotonically increasing sequence for each $p$ on the compact set $\mathcal{P}(p_{\max},a_{\max})$. By Dini's Theorem \cite[Theorem~7.13]{rudin1964principles}, $|\costeval{\eqref{prob: equivalent finite dim approx integrated optimization problem}_N}{p} - \costeval{\eqref{prob: new equivalent IOCA}}{p}| \rightarrow 0$ uniformly on $\mathcal{P}(p_{\max},a_{\max})$ as $N \rightarrow \infty$. 
    By Moore-Osgood Theorem \cite[Theorem~7.11]{rudin1964principles}, this uniform convergence and the convergence $p_{N_k}^* \rightarrow p_{\inf}^*$ as $k \rightarrow \infty$ (see \eqref{eq: convergence of p_Nk to p_inf}) imply that $\lim_{k \rightarrow \infty} \costeval{\eqref{prob: new equivalent IOCA}}{p_{N_k}^*} = \lim_{j \rightarrow \infty} \lim_{k \rightarrow \infty} \optcosteval{\eqref{prob: equivalent finite dim approx integrated optimization problem}_j}{p_{N_k}^*}$, in which the iterated limit equals the double limit \cite[p.~140]{taylor1985general}%https://www.google.com/books/edition/General_Theory_of_Functions_and_Integrat/pczdngEACAAJ?hl=en&gbpv=1&bsq=140
    , i.e., 
    \begin{align*}
        \lim_{j \rightarrow \infty} \lim_{k \rightarrow \infty} \optcosteval{\eqref{prob: equivalent finite dim approx integrated optimization problem}_j}{p_{N_k}^*} = & \ \lim_{\substack{j \rightarrow \infty \\k \rightarrow \infty}} \optcosteval{\eqref{prob: equivalent finite dim approx integrated optimization problem}_{j}}{p_{N_k}^*} \\
        = & \ \lim_{k \rightarrow \infty }\optcosteval{\eqref{prob: equivalent finite dim approx integrated optimization problem}}{p_{N_k}^*}.
    \end{align*}
    }

    The second limit on the right-hand side of \eqref{eq: first limit in the convergence of approximated optimal solution} is zero due to Lemma~\ref{lemma: continuity of the total cost wrt guidance}-2. Hence, it follows from \eqref{eq: first limit in the convergence of approximated optimal solution} that $\lim_{k \rightarrow \infty} \optcosteval{\eqref{prob: equivalent finite dim approx integrated optimization problem}}{p_{N_k}^*} = \costeval{\eqref{prob: new equivalent IOCA}}{p_{\inf}^*}$, which implies
    \begin{align}
        \liminf_{N \rightarrow \infty} \optcosteval{\eqref{prob: equivalent finite dim approx integrated optimization problem}}{p_N^*} = & \  \lim_{k \rightarrow \infty} \optcosteval{\eqref{prob: equivalent finite dim approx integrated optimization problem}}{p_{N_k}^*} \nonumber \\
        = & \ \costeval{\eqref{prob: new equivalent IOCA}}{p_{\inf}^*} \nonumber \\
        \geq & \ \optcosteval{\eqref{prob: new equivalent IOCA}}{p^*}. \label{eq: intermediate step 4 for solution convergence} 
    \end{align}
    Therefore, we conclude $\lim_{N \rightarrow \infty} \optcosteval{\eqref{prob: equivalent finite dim approx integrated optimization problem}}{p_N^*} = \optcosteval{\eqref{prob: new equivalent IOCA}}{p^*}$ from \eqref{eq: intermediate step 3 for solution convergence} and \eqref{eq: intermediate step 4 for solution convergence}.
    
    Next, we show \eqref{eq: convergence of the approximated optimal guidance}, i.e., $ |\costeval{\eqref{prob: new equivalent IOCA}}{p_N^*} - \optcosteval{\eqref{prob: new equivalent IOCA}}{p^*}| \rightarrow 0$ as $N \rightarrow \infty$. 
    \revision{We start with $\optcosteval{\eqref{prob: new equivalent IOCA}}{p^*} \leq \costeval{\eqref{prob: new equivalent IOCA}}{p_N^*}$ for all $N$, which implies that
    \begin{equation}\label{eq: intermediate step 15}
        \optcosteval{\eqref{prob: new equivalent IOCA}}{p^*} \leq \liminf_{N \rightarrow \infty} \costeval{\eqref{prob: new equivalent IOCA}}{p_N^*}.
    \end{equation}
    To prove \eqref{eq: convergence of the approximated optimal guidance}, what remains to be shown is $\optcosteval{\eqref{prob: new equivalent IOCA}}{p^*} \geq \limsup_{N \rightarrow \infty} \costeval{\eqref{prob: new equivalent IOCA}}{p_N^*}$. Choose a convergent subsequence $\{\costeval{\eqref{prob: new equivalent IOCA}}{p_{N_j}^*} \}_{j = 1}^{\infty}$ such that $\lim_{j \rightarrow \infty} \costeval{\eqref{prob: new equivalent IOCA}}{p_{N_j}^*} = \limsup_{N \rightarrow \infty} \costeval{\eqref{prob: new equivalent IOCA}}{p_N^*}$. Since $\{p_{N_j}^* \}_{j=1}^{\infty} \subset \mathcal{P}(p_{\max},a_{\max})$ is uniformly equicontinuous and uniformly bounded, by Arzel\`a-Ascoli Theorem \cite{royden2010real}, the sequence has a (uniformly) convergent subsequence which we denote with the same indices $ N_j$ to simplify notation. Denote the limit of $\{p_{N_j}^* \}_{j=1}^{\infty}$ by $p_{\sup}^*$ such that
    \begin{equation}\label{eq: pointwise convergence to p_sup}
        \lim_{j \rightarrow \infty} \norm{p_{N_j}^* - p_{\sup}^*}_{C([0,t_f];\mathbb{R}^m)} = 0.
    \end{equation}
    Due to the continuity of $\costeval{\eqref{prob: new equivalent IOCA}}{\cdot}$ (see Lemma~\ref{lemma: continuity of the total cost wrt guidance}-1), we have 
    \begin{equation*}
    \costeval{\eqref{prob: new equivalent IOCA}}{p_{\sup}^*} = \lim_{j \rightarrow \infty} \costeval{\eqref{prob: new equivalent IOCA}}{p_{N_j}^*} =  \limsup_{N \rightarrow \infty }\costeval{\eqref{prob: new equivalent IOCA}}{p_N^*}.
    \end{equation*}
    It follows that
    \begin{align}
        & \costeval{\eqref{prob: new equivalent IOCA}}{p_{\sup}^*}  \nonumber \\
        \leq & |\costeval{\eqref{prob: new equivalent IOCA}}{p_{\sup}^*} - \optcosteval{\eqref{prob: new equivalent IOCA}}{p^*}| + \optcosteval{\eqref{prob: new equivalent IOCA}}{p^*} \nonumber \\
        =& |\costeval{\eqref{prob: new equivalent IOCA}}{p_{\sup}^*} - \lim_{N \rightarrow \infty} \optcosteval{\eqref{prob: equivalent finite dim approx integrated optimization problem}}{p_N^*}| + \optcosteval{\eqref{prob: new equivalent IOCA}}{p^*} \nonumber \\
        =& |\costeval{\eqref{prob: new equivalent IOCA}}{p_{\sup}^*} - \lim_{j \rightarrow \infty} \optcosteval{\eqref{prob: equivalent finite dim approx integrated optimization problem}}{p_{N_j}^*}| + \optcosteval{\eqref{prob: new equivalent IOCA}}{p^*}. \label{eq: intermediate step 13}
    \end{align}
    % First notice that for all $p \in \mathcal{P}(p_{\max},a_{\max})$, $\costeval{\eqref{prob: equivalent finite dim approx integrated optimization problem}_N}{p}$ converges to $\costeval{\eqref{prob: new equivalent IOCA}}{p}$ pointwise as the dimension of approximation $N$ goes to infinity (see \eqref{eq: intermediate step 1 for solution convergence}). Furthermore, 
    Since the sequence of approximated PDE cost $\{\langle {Z_N(0)},$ ${\Pi_N(0) Z_N(0)} \rangle \}_{N=1}^{\infty}$ is monotonically increasing, the sequence $\{ \costeval{\eqref{prob: equivalent finite dim approx integrated optimization problem}_N}{p}\}_{N=1}^{\infty}$ is a monotonically increasing sequence for each $p$ on the compact set $\mathcal{P}(p_{\max},a_{\max})$. Since $\lim_{N \rightarrow \infty}\costeval{\eqref{prob: equivalent finite dim approx integrated optimization problem}_N}{p} = \costeval{\eqref{prob: new equivalent IOCA}}{p}$ for all $p \in \mathcal{P}(p_{\max},a_{\max})$ (see Lemma~\ref{lemma: continuity of the total cost wrt guidance}-1), by Dini's Theorem \cite[Theorem~7.13]{rudin1964principles}, the limit holds uniformly on $\mathcal{P}(p_{\max},a_{\max})$ as $N \rightarrow \infty$. 
    By Moore-Osgood Theorem \cite[Theorem~7.11]{rudin1964principles}, this uniform convergence and the convergence $p_{N_j}^* \rightarrow p_{\sup}^*$ as $j \rightarrow \infty$ (see \eqref{eq: pointwise convergence to p_sup}) imply that
    \begin{align}\label{eq: convergence due to the continuity of original cost function 2}
        \costeval{\eqref{prob: new equivalent IOCA}}{p_{\sup}^*} = \lim_{k \rightarrow \infty} \lim_{j \rightarrow \infty}  \optcosteval{\eqref{prob: equivalent finite dim approx integrated optimization problem}_k}{p_{N_j}^*}.
    \end{align}
    Furthermore, the iterated limit equals the double limit \cite[p.~140]{taylor1985general}%https://www.google.com/books/edition/General_Theory_of_Functions_and_Integrat/pczdngEACAAJ?hl=en&gbpv=1&bsq=140
    , i.e.,
    \begin{align}
        \lim_{k \rightarrow \infty} \lim_{j \rightarrow \infty}  \optcosteval{\eqref{prob: equivalent finite dim approx integrated optimization problem}_k}{p_{N_j}^*} = & \lim_{\substack{j \rightarrow \infty \\ k \rightarrow \infty}} \optcosteval{\eqref{prob: equivalent finite dim approx integrated optimization problem}_k}{p_{N_j}^*} \nonumber \\
        = & \lim_{j \rightarrow \infty}\optcosteval{\eqref{prob: equivalent finite dim approx integrated optimization problem}}{p_{N_j}^*}. \label{eq: intermediate step 16}
    \end{align}
    Hence, combining \eqref{eq: intermediate step 13}--\eqref{eq: intermediate step 16}, we have $ \optcosteval{\eqref{prob: new equivalent IOCA}}{p^*} \geq \costeval{\eqref{prob: new equivalent IOCA}}{p_{\sup}^*}  = \limsup_{N \rightarrow \infty} \costeval{\eqref{prob: new equivalent IOCA}}{p_N^*}$, from which and \eqref{eq: intermediate step 15} we conclude the desired convergence $\lim_{N \rightarrow \infty}\costeval{\eqref{prob: new equivalent IOCA}}{p_N^*} = \optcosteval{\eqref{prob: new equivalent IOCA}}{p^*}$.}
\qed 
\end{pf*}

\bibliographystyle{abbrv}
\bibliography{reference}

\end{document}

\begin{frontmatter}
\date{15 June 2021}
\title{Supplementary material for: \\
Optimal control of a 2D diffusion-advection process with a team of mobile actuators under jointly optimal guidance\thanksref{footnoteinfo}} % Title, preferably not more 

\thanks[footnoteinfo]{This paper was not presented at any IFAC 
meeting. Corresponding author Sheng Cheng.}

\author[Sheng]{Sheng Cheng}\ead{cheng@terpmail.umd.edu},    % Add the 
\author[Derek]{Derek A. Paley}\ead{dpaley@umd.edu}               % e-mail address 
% \author[Baiae]{Publius Maro Vergilius}\ead{vergilius@culture.ir}  % (ead) as shown

\address[Sheng]{University of Maryland, College Park}  % Please supply                                              
\address[Derek]{University of Maryland, College Park}             % full addresses
% \address[Baiae]{The White House, Baiae}        % here.

\begin{keyword}                           % Five to ten keywords, 
Infinite-dimensional systems; Multi-agent systems; Modeling for control optimization; Guidance navigation and control. 
% chosen from the IFAC 
\end{keyword}                             % keyword list or with the 
                                          % help of the Automatica 
                                          % keyword wizard

\begin{abstract}                          % Abstract of not more than 200 words.
This supplement provides proofs for Lemmas~2.3, 2.5, and D.1 in [3].
\end{abstract}

\end{frontmatter}

\begin{customlem}{2.3}\label{lemma: continuity of equivalent cost wrt actuator location}
Suppose $\mathcal{Z}_0 \in \mathcal{H}$. Let assumptions (A1)--(A3) hold with $q = 1$ and $\Pi \in C([0,t_f];\mathcal{J}_1(\mathcal{H}))$ be defined as in \cite[(12)]{Cheng2020optimalInReivew}. If assumption (A4) holds, then the mapping $K: C([0,t_f];\mathbb{R}^n) \rightarrow \mathbb{R}^+$ such that $K(\xi) \defeq \innerproduct{\mathcal{Z}_0}{\Pi(0) \mathcal{Z}_0}$ is continuous.
\end{customlem}

\begin{pf*}{Proof}
Without loss of generality, consider the case of one mobile actuator, i.e., $m_a = 1$. The case of multiple actuators follows naturally. We first show a consequence of the input operator $\mathcal{B}$ being continuous with respect to location. Consider two actuator states, $\xi_1$ and $\xi_2 \in C([0,t_f];\mathbb{R}^n)$. For any $\phi \in \mathcal{H} = L^2(\Omega)$ and all $t \in [0,t_f]$,
\begin{align}
     & |\mathcal{B}^{\star}(M\xi_1(t),t) \phi - \mathcal{B}^{\star} (M\xi_2(t),t)  \phi| \nonumber \\
     \leq &    \norm{\mathcal{B}(M\xi_1(t),t) - \mathcal{B}(M \xi_2(t),t)}_{L^2(\Omega)} \norm{\phi}_{L^2(\Omega)} \nonumber \\
    \leq & l\left(|M(\xi_1(t) - \xi_2(t))|_2 \right) \norm{\phi}_{L^2(\Omega)}, \label{eq: inequality 1 in proving continuity of PDE cost wrt actuator location} 
\end{align}
where we use the fact that $\mathcal{B}(\cdot,\cdot)$ is the integral kernel of $\mathcal{B}^{\star}(\cdot,\cdot)$. Hence,
\begin{align}\label{eq: intermediate step 1}
    & \norm{\mathcal{B}^{\star}(M\xi_1(t),t) - \mathcal{B}^{\star} (M\xi_2(t),t)}_{\mathcal{L}(\mathcal{H};\mathbb{R})} \nonumber \\
    \leq & l\left(|M(\xi_1(t) - \xi_2(t))|_2 \right).
\end{align}

Since $\mathbb{R}$ is finite-dimensional, there exists $c_1 > 0$ such that\cite[Proof of Lemma 4.3]{burns2015infinite}
\begin{align}\label{eq: intermediate step 2}
    & \norm{\mathcal{B}^{\star}(M\xi_1(t),t) - \mathcal{B}^{\star} (M\xi_2(t),t)}_{\mathcal{J}_1(\mathcal{H};\mathbb{R})}\nonumber \\
    \leq & c_1 \norm{\mathcal{B}^{\star}(M\xi_1(t),t) - \mathcal{B}^{\star} (M\xi_2(t),t)}_{\mathcal{L}(\mathcal{H};\mathbb{R})}   %\\
    % \leq & c_1 l\left(|M(\xi_1(t) - \xi_2(t))|_2 \right).
\end{align}
For brevity, we shall use $\mathcal{B}_1(t)$ for $\mathcal{B}(M\xi_1(t),t)$ and $\mathcal{B}_2(t)$ for $\mathcal{B}(M\xi_2(t),t)$. Now,
\begin{align}
    & \norm{\mathcal{B}_1(t) R^{-1} \mathcal{B}_1^{\star}(t) - \mathcal{B}_2(t) R^{-1} \mathcal{B}_2^{\star}(t)}_{\mathcal{J}_1(\mathcal{H})} \nonumber \\
    \leq & \norm{\mathcal{B}_1(t) R^{-1} }_{\mathcal{J}_1(\mathbb{R};\mathcal{H})} \norm{\mathcal{B}_1^{\star} (t) - \mathcal{B}_2^{\star}(t)}_{\mathcal{J}_1(\mathcal{H};\mathbb{R})} \nonumber \\
    & + \norm{R^{-1} \mathcal{B}_2^{\star}(t)}_{\mathcal{J}_1(\mathcal{H};\mathbb{R})} \norm{\mathcal{B}_1(t) - \mathcal{B}_2(t)}_{\mathcal{J}_1(\mathbb{R};\mathcal{H})} \nonumber \\
    =& (\norm{\mathcal{B}_1(t) R^{-1} }_{\mathcal{J}_1(\mathbb{R};\mathcal{H})}+ \norm{R^{-1} \mathcal{B}_2^{\star}(t)}_{\mathcal{J}_1(\mathcal{H};\mathbb{R})} ) \nonumber \\
    & \norm{\mathcal{B}_1^{\star} (t) - \mathcal{B}_2^{\star}(t)}_{\mathcal{J}_1(\mathcal{H};\mathbb{R})}
  \nonumber \\
    \leq & c_2  \norm{\mathcal{B}_1^{\star} (t) - \mathcal{B}_2^{\star}(t)}_{\mathcal{J}_1(\mathcal{H};\mathbb{R})} \nonumber \\
    \leq & c_2 c_1 l(|M(\xi_1(t)-\xi_2(t))|_2) 
    \label{eq: intermediate step 4}
\end{align}
for some $c_2 >0$ where the last inequality follows from \eqref{eq: intermediate step 1} and \eqref{eq: intermediate step 2}. %(The equality in the third line holds roughly because for $A \in \mathcal{L}(X,Y) $, we have $\text{dom}(A^{\star}) = Y^{\star}$ and $\norm{A^{\star}}_{\mathcal{L}(Y^{\star},X^{\star})} = \norm{A}_{\mathcal{L}(X,Y)}$).

We now continue to prove that $K:C([0,t_f];\mathbb{R}^n) \rightarrow \mathbb{R}$ is a continuous mapping. For brevity, we use $\Pi_{1}(0)$ and $\Pi_{2}(0)$ for the Riccati operator associated with trajectory $\xi_1$ and $\xi_2$, respectively. We also suppress the usage of the time argument of the integrand in the following derivation. \revision{We start with $K(\xi_1) - K(\xi_2)$:
\begin{align*}
    & K(\xi_1)-K(\xi_2) \nonumber \\ %\innerproduct{\mathcal{Z}_0}{\int_0^{t_f} \mathcal{S}^{\star}(\tau) \left( 
    % \begin{aligned}
    %     & \Pi_{1}(\tau) \mathcal{B}_1(\tau) R^{-1} \mathcal{B}_1^{\star}(\tau) \Pi_{1}(\tau) \nonumber \\
    %     -&\Pi_{2}(\tau) \mathcal{B}_2(\tau) R^{-1} \mathcal{B}_2^{\star}(\tau) \Pi_{2}(\tau) 
    % \end{aligned}
    % \right) \mathcal{S}(\tau) \dd \tau \mathcal{Z}_0 } \nonumber \\
    % = & 
    = &  \innerproduct{\mathcal{Z}_0}{  \int_0^{t_f} \mathcal{S}^{\star} \Pi_{1} \mathcal{B}_1 R^{-1} \big(\mathcal{B}_1^{\star} \Pi_{1} - \mathcal{B}_2^{\star} \Pi_{2}  \big) \mathcal{S} \dd \tau \mathcal{Z}_0 } \nonumber \\
    & + \innerproduct{\mathcal{Z}_0}{ \int_0^{t_f} \mathcal{S}^{\star} \big( \Pi_{1}  \mathcal{B}_1 - \Pi_{2} \mathcal{B}_2 \big) R^{-1}  \mathcal{B}_2^{\star} \Pi_{2}   \mathcal{S} \dd \tau \mathcal{Z}_0 }.% \label{eq: intermediate step 3} 
\end{align*}
% \begin{figure*}[t]
% \begin{align}
%     K(\xi_1)-K(\xi_2) = & %\innerproduct{\mathcal{Z}_0}{\int_0^{t_f} \mathcal{S}^{\star}(\tau) \left( 
%     % \begin{aligned}
%     %     & \Pi_{1}(\tau) \mathcal{B}_1(\tau) R^{-1} \mathcal{B}_1^{\star}(\tau) \Pi_{1}(\tau) \nonumber \\
%     %     -&\Pi_{2}(\tau) \mathcal{B}_2(\tau) R^{-1} \mathcal{B}_2^{\star}(\tau) \Pi_{2}(\tau) 
%     % \end{aligned}
%     % \right) \mathcal{S}(\tau) \dd \tau \mathcal{Z}_0 } \nonumber \\
%     % = & 
%     \innerproduct{\mathcal{Z}_0}{\int_0^{t_f} \mathcal{S}^{\star}(\tau) \Pi_{1}(\tau) \mathcal{B}_1(\tau) R^{-1} \big(\mathcal{B}_1^{\star}(\tau) \Pi_{1}(\tau) - \mathcal{B}_2^{\star}(\tau) \Pi_{2}(\tau)  \big) \mathcal{S}(\tau) \dd \tau \mathcal{Z}_0 } \nonumber \\
%     & + \innerproduct{\mathcal{Z}_0}{\int_0^{t_f} \mathcal{S}^{\star}(\tau) \big( \Pi_{1}(\tau)  \mathcal{B}_1(\tau) - \Pi_{2}(\tau) \mathcal{B}_2(\tau) \big) R^{-1}  \mathcal{B}_2^{\star}(\tau) \Pi_{2}(\tau)   \mathcal{S}(\tau) \dd \tau \mathcal{Z}_0 }. \label{eq: intermediate step 3} 
% \end{align}
% %\hrulefill
% \end{figure*}
Our analysis continues with the first term on the right-hand side%of \eqref{eq: intermediate step 3}
, because the second term can be analyzed similarly using the following derivation: %in \eqref{eq: intermediate step 7}.
\begin{align}
    & \innerproduct{\mathcal{Z}_0}{ \int_0^{t_f} \mathcal{S}^{\star} \Pi_{1} \mathcal{B}_1 R^{-1} \big(\mathcal{B}_1^{\star} \Pi_{1} - \mathcal{B}_2^{\star} \Pi_{2}  \big) \mathcal{S} \dd \tau \mathcal{Z}_0 } \nonumber \\
    = & \innerproduct{\mathcal{Z}_0}{ \int_0^{t_f} \mathcal{S}^{\star} \Pi_{1} \mathcal{B}_1 R^{-1} \mathcal{B}_1^{\star} \big( \Pi_{1} - \Pi_{2}  \big) \mathcal{S} \dd \tau \mathcal{Z}_0 } \nonumber \\
    & + \innerproduct{\mathcal{Z}_0}{ \int_0^{t_f} \mathcal{S}^{\star} \Pi_{1} \mathcal{B}_1 R^{-1} \big(\mathcal{B}_1^{\star} - \mathcal{B}_2^{\star} \big) \Pi_{2}  \mathcal{S} \dd \tau \mathcal{Z}_0 }. \label{eq: intermediate step 7}
\end{align}
% \begin{figure*}[t]
% \begin{align}
%     & \innerproduct{\mathcal{Z}_0}{\int_0^{t_f} \mathcal{S}^{\star}(\tau) \Pi_{1}(\tau) \mathcal{B}_1(\tau) R^{-1} \big(\mathcal{B}_1^{\star}(\tau) \Pi_{1}(\tau) - \mathcal{B}_2^{\star}(\tau) \Pi_{2}(\tau)  \big) \mathcal{S}(\tau) \dd \tau \mathcal{Z}_0 } \nonumber \\
%     = & \innerproduct{\mathcal{Z}_0}{\int_0^{t_f} \mathcal{S}^{\star}(\tau) \Pi_{1}(\tau) \mathcal{B}_1(\tau) R^{-1} \mathcal{B}_1^{\star}(\tau) \big( \Pi_{1}(\tau) - \Pi_{2}(\tau)  \big) \mathcal{S}(\tau) \dd \tau \mathcal{Z}_0 } \nonumber \\
%     & + \innerproduct{\mathcal{Z}_0}{\int_0^{t_f} \mathcal{S}^{\star}(\tau) \Pi_{1}(\tau) \mathcal{B}_1(\tau) R^{-1} \big(\mathcal{B}_1^{\star}(\tau) - \mathcal{B}_2^{\star}(\tau) \big) \Pi_{2}(\tau)  \mathcal{S}(\tau) \dd \tau \mathcal{Z}_0 }. \label{eq: intermediate step 7}
% \end{align}
% \end{figure*}
Take absolute values on both sides, we get
\begin{align*}
    & |\innerproduct{\mathcal{Z}_0}{ \int_0^{t_f} \mathcal{S}^{\star} \Pi_{1} \mathcal{B}_1 R^{-1} \big(\mathcal{B}_1^{\star} \Pi_{1} - \mathcal{B}_2^{\star} \Pi_{2}  \big) \mathcal{S} \dd \tau \mathcal{Z}_0 } | \nonumber \\
    \leq & \norm{\mathcal{Z}_0}^2_{\mathcal{H}}  \int_0^{t_f} ( \norm{\mathcal{S}^{\star}}_{\mathcal{L}(\mathcal{H})} \norm{\Pi_{1} \mathcal{B}_1 R^{-1} \mathcal{B}_1^{\star}}_{\mathcal{L}(\mathcal{H})} \nonumber \\
    & \phantom{aaaaaaaaaaaaaaaaaaaa} \norm{\Pi_{1} - \Pi_{2}}_{\mathcal{L}({\mathcal{H})}}\norm{\mathcal{S}}_{\mathcal{L}(\mathcal{H})} ) \dd \tau \nonumber \\
    & + \norm{\mathcal{Z}_0}^2_{\mathcal{H}} \int_0^{t_f}( \norm{\mathcal{S}^{\star}}_{\mathcal{L}(\mathcal{H})} \norm{\Pi_{1} \mathcal{B}_1 R^{-1}}_{\mathcal{L}(\mathbb{R};\mathcal{H})} \norm{\Pi_{2}}_{\mathcal{L}(\mathcal{H})} \nonumber \\
    & \phantom{aaaaaaaaaaaaaaaaaaa}  \norm{\mathcal{B}_1^{\star} - \mathcal{B}_2^{\star} }_{\mathcal{L}(\mathcal{H};\mathbb{R})} \norm{\mathcal{S}}_{\mathcal{L}(\mathcal{H})} )\dd \tau. %\label{eq: intermediate step 8}
\end{align*}
% \begin{figure*}[t]
% \begin{align}
%     & |\innerproduct{\mathcal{Z}_0}{\int_0^{t_f} \mathcal{S}^{\star}(\tau) \Pi_{1}(\tau) \mathcal{B}_1(\tau) R^{-1} \big(\mathcal{B}_1^{\star}(\tau) \Pi_{1}(\tau) - \mathcal{B}_2^{\star}(\tau) \Pi_{2}(\tau)  \big) \mathcal{S}(\tau) \dd \tau \mathcal{Z}_0 } | \nonumber \\
%     \leq & \norm{\mathcal{Z}_0}^2_{\mathcal{H}} \int_0^{t_f} \norm{\mathcal{S}^{\star}(\tau)}_{\mathcal{L}(\mathcal{H})} \norm{\Pi_{1}(\tau) \mathcal{B}_1(\tau) R^{-1} \mathcal{B}_1^{\star}(\tau)}_{\mathcal{L}(\mathcal{H})} \norm{\Pi_{1}(\tau) - \Pi_{2}(\tau)}_{\mathcal{L}({\mathcal{H})}}\norm{\mathcal{S}(\tau)}_{\mathcal{L}(\mathcal{H})} \dd \tau \nonumber \\
%     & + \norm{\mathcal{Z}_0}^2_{\mathcal{H}} \int_0^{t_f} \norm{\mathcal{S}^{\star}(\tau)}_{\mathcal{L}(\mathcal{H})} \norm{\Pi_{1}(\tau) \mathcal{B}_1(\tau) R^{-1}}_{\mathcal{L}(\mathbb{R};\mathcal{H})} \norm{\Pi_{2}(\tau)}_{\mathcal{L}(\mathcal{H})} \norm{\mathcal{B}_1^{\star}(\tau) - \mathcal{B}_2^{\star}(\tau) }_{\mathcal{L}(\mathcal{H};\mathbb{R})} \norm{\mathcal{S}(\tau)}_{\mathcal{L}(\mathcal{H})} \dd \tau. \label{eq: intermediate step 8}
% \end{align}
% \end{figure*}
Since there exists~$c_{t_f}$ such that $\norm{\mathcal{S}(t)}_{\mathcal{L}(\mathcal{H})} \leq c_{t_f}$ and $\norm{\mathcal{S}^{\star}(t)}_{\mathcal{L}(\mathcal{H})} \leq c_{t_f}$ for all $t \in [0,t_f]$, it follows that
\begin{align*}
    & |\innerproduct{\mathcal{Z}_0}{ \int_0^{t_f} \mathcal{S}^{\star} \Pi_{1} \mathcal{B}_1 R^{-1} \big(\mathcal{B}_1^{\star} \Pi_{1} - \mathcal{B}_2^{\star} \Pi_{2}  \big) \mathcal{S} \dd \tau \mathcal{Z}_0 } | \nonumber  \\
    \leq &  \norm{\mathcal{Z}_0}^2_{\mathcal{H}} c_{t_f}^2
     \int_0^{t_f} ( \norm{\Pi_{1} \mathcal{B}_1 R^{-1} \mathcal{B}_1^{\star}}_{\mathcal{L}(\mathcal{H})} \norm{\Pi_{1} - \Pi_{2}}_{\mathcal{L}({\mathcal{H})}}  \nonumber \\
    &  + \norm{\Pi_{1} \mathcal{B}_1 R^{-1}}_{\mathcal{L}(\mathbb{R};\mathcal{H})} \norm{\Pi_{2}}_{\mathcal{L}(\mathcal{H})} \norm{\mathcal{B}_1^{\star} - \mathcal{B}_2^{\star} }_{\mathcal{L}(\mathcal{H};\mathbb{R})} ) \dd \tau. %    \label{eq: intermediate step 9}
\end{align*}
% \begin{figure*}[t]
% \begin{align}
%     & |\innerproduct{\mathcal{Z}_0}{\int_0^{t_f} \mathcal{S}^{\star}(\tau) \Pi_{1}(\tau) \mathcal{B}_1(\tau) R^{-1} \big(\mathcal{B}_1^{\star}(\tau) \Pi_{1}(\tau) - \mathcal{B}_2^{\star}(\tau) \Pi_{2}(\tau)  \big) \mathcal{S}(\tau) \dd \tau \mathcal{Z}_0 } | \nonumber  \\
%     \leq &  \norm{\mathcal{Z}_0}^2_{\mathcal{H}} c_{t_f}^2
%     \left( \int_0^{t_f}  \norm{\Pi_{1}(\tau) \mathcal{B}_1(\tau) R^{-1} \mathcal{B}_1^{\star}(\tau)}_{\mathcal{L}(\mathcal{H})} \norm{\Pi_{1}(\tau) - \Pi_{2}(\tau)}_{\mathcal{L}({\mathcal{H})}} \dd \tau \right. \nonumber \\
%     & \left. + \int_0^{t_f}  \norm{\Pi_{1}(\tau) \mathcal{B}_1(\tau) R^{-1}}_{\mathcal{L}(\mathbb{R};\mathcal{H})} \norm{\Pi_{2}(\tau)}_{\mathcal{L}(\mathcal{H})} \norm{\mathcal{B}_1^{\star}(\tau) - \mathcal{B}_2^{\star}(\tau) }_{\mathcal{L}(\mathcal{H};\mathbb{R})}  \dd \tau \right). 
%     \label{eq: intermediate step 9}
% \end{align}
% \end{figure*}
Since $\mathcal{J}_q(\mathcal{H}) \hookrightarrow \mathcal{L}(\mathcal{H})$ \cite{Burns2015Solutions}, there exists $c_3 > 0$ such that $\norm{\Pi_{1}(\tau) - \Pi_{2}(\tau)}_{\mathcal{L}(\mathcal{H})} \leq  c_3 \norm{\Pi_{1}(\tau) - \Pi_{2}(\tau)}_{\mathcal{J}_q(\mathcal{H})}$ for ${1 \leq q < \infty}$. Hence, the following inequality holds for $q = 1$:
\begin{align}
    & |\innerproduct{\mathcal{Z}_0}{ \int_0^{t_f} \mathcal{S}^{\star} \Pi_{1} \mathcal{B}_1 R^{-1} \big(\mathcal{B}_1^{\star} \Pi_{1} - \mathcal{B}_2^{\star} \Pi_{2}  \big) \mathcal{S} \dd \tau \mathcal{Z}_0 } | \nonumber \\
    \leq & \norm{\mathcal{Z}_0}^2_{\mathcal{H}} c_{t_f}^2 
          \int_0^{t_f} (c_3 \norm{\Pi_{1} \mathcal{B}_1 R^{-1} \mathcal{B}_1^{\star}}_{\mathcal{L}(\mathcal{H})} \norm{\Pi_{1} - \Pi_{2}}_{\mathcal{J}_1(\mathcal{H})}  \nonumber \\
        & +   \norm{\Pi_{1} \mathcal{B}_1 R^{-1}}_{\mathcal{L}(\mathbb{R};\mathcal{H})} \norm{\Pi_{2}}_{\mathcal{L}(\mathcal{H})} \norm{\mathcal{B}_1^{\star} - \mathcal{B}_2^{\star} }_{\mathcal{L}(\mathcal{H};\mathbb{R})} ) \dd \tau \nonumber \\
    % \leq & \norm{\mathcal{Z}_0}^2_{\mathcal{H}} c_{t_f}^2 \left( 
    %     c_3 \int_0^{t_f}  \norm{\Pi_{1}(\tau) \mathcal{B}_1(\tau) R^{-1} \mathcal{B}_1^{\star}(\tau)}_{\mathcal{L}(\mathcal{H})} \norm{\Pi_{1}(\tau) - \Pi_{2}(\tau)}_{\mathcal{J}_1(\mathcal{H})} \dd \tau \right. \nonumber \\
    %     & + \left. \int_0^{t_f}  \norm{\Pi_{1}(\tau) \mathcal{B}_1(\tau) R^{-1}}_{\mathcal{L}(\mathbb{R};\mathcal{H})} \norm{\Pi_{2}(\tau)}_{\mathcal{L}(\mathcal{H})}   \dd \tau \norm{\mathcal{B}_1^{\star} - \mathcal{B}_2^{\star} }_{L^{\infty}([0,t_f];\mathcal{L}(\mathcal{H};\mathbb{R}))}
    % \right) \nonumber \\
    \leq & \norm{\mathcal{Z}_0}^2_{\mathcal{H}} c_{t_f}^2 \Big(
        c_3  \int_0^{t_f}  \norm{\Pi_{1} \mathcal{B}_1 R^{-1} \mathcal{B}_1^{\star}}_{\mathcal{L}(\mathcal{H})}\dd \tau \nonumber \\
        & \phantom{aaaaaaaaaaaaaaaaaaa} \esssup_{t \in [0,t_f]} \norm{\Pi_{1}(t) - \Pi_{2}(t)}_{\mathcal{J}_1(\mathcal{H})}   \nonumber \\
        & +  \int_0^{t_f}  \norm{\Pi_{1} \mathcal{B}_1 R^{-1}}_{\mathcal{L}(\mathbb{R};\mathcal{H})} \norm{\Pi_{2}}_{\mathcal{L}(\mathcal{H})}   \dd \tau \nonumber \\
        &  \phantom{aaaaaaaaaa} \esssup_{t \in [0,t_f]}\norm{\mathcal{B}_1^{\star}(t) - \mathcal{B}_2^{\star}(t) }_{\mathcal{L}(\mathcal{H};\mathbb{R})}\Big). \label{eq: intermediate step 10}
\end{align}
}
We now wish to bound $\esssup_{t \in [0,t_f]} \norm{\Pi_{1}(t) - \Pi_{2}(t)}_{\mathcal{J}_1(\mathcal{H})}$ and $\esssup_{t \in [0,t_f]}\norm{\mathcal{B}_1^{\star}(t) - \mathcal{B}_2^{\star}(t) }_{\mathcal{L}(\mathcal{H};\mathbb{R})}$.

Notice that $\Pi(\cdot)$ is continuous on $[0,t_f]$ into $\mathcal{J}_1(\mathcal{H})$ \cite[Lemma~2.1]{Cheng2020optimalInReivew}
\begin{align}
    & \esssup_{t \in [0,t_f]} \norm{\Pi_{1}(t) - \Pi_{2}(t)}_{\mathcal{J}_1(\mathcal{H})} \nonumber \\
    = & \sup_{t\in [0,t_f]} \norm{\Pi_{1}(t) - \Pi_{2}(t)}_{\mathcal{J}_1(\mathcal{H})}. \label{eq: intermediate step 11}
\end{align}
By \cite[(12)]{Cheng2020optimalInReivew}%\eqref{eq: new mild solution of operator Riccati equation}
, the mapping $\Pi: [0,t_f] \rightarrow \mathcal{J}_1(\mathcal{H})$ varies continuously in $\sup_{t \in [0,t_f]} \norm{\cdot}_{\mathcal{J}_1(\mathcal{H})}$-norm with respect to $\bar{\mathcal{B}} \bar{\mathcal{B}}^{\star}(\cdot)$~\cite{Burns2015Solutions}. Hence, there exists $c_4 >0$ such that 
\begin{align}
     & \sup_{t\in [0,t_f]} \norm{\Pi_{1}(t) - \Pi_{2}(t)}_{\mathcal{J}_1(\mathcal{H})}  \nonumber \\
    \leq & \sup_{t\in [0,t_f]} c_4 \norm{\bar{\mathcal{B}}_1 \bar{\mathcal{B}}_1^{\star}(t) - \bar{\mathcal{B}}_2 \bar{\mathcal{B}}_2^{\star}(t)}_{\mathcal{J}_1(\mathcal{H})}. \label{eq: intermediate step 12}
\end{align}
Recall $\mathcal{B}(\cdot) R^{-1} \mathcal{B}^{\star}(\cdot) = \bar{\mathcal{B}} \bar{\mathcal{B}}^{\star}(\cdot)$ and combine \eqref{eq: intermediate step 4}, \eqref{eq: intermediate step 11}, and \eqref{eq: intermediate step 12}:
\begin{align}
    & \esssup_{t \in [0,t_f]} \norm{\Pi_{1}(t) - \Pi_{2}(t)}_{\mathcal{J}_1(\mathcal{H})} \nonumber \\
    \leq & \sup_{t\in [0,t_f]} c_4 c_1 c_2 l(|M(\xi_1(t) - \xi_2(t))|_2).\label{eq: intermediate step 5}
\end{align}
It remains to bound $\esssup_{t \in [0,t_f]}\norm{\mathcal{B}_1^{\star}(t) - \mathcal{B}_2^{\star}(t) }_{\mathcal{L}(\mathcal{H};\mathbb{R})}$. By \eqref{eq: intermediate step 1},
\begin{align}
     & \esssup_{t \in [0,t_f]}\norm{\mathcal{B}_1^{\star}(t) - \mathcal{B}_2^{\star}(t) }_{\mathcal{L}(\mathcal{H};\mathbb{R})} \nonumber \\
     \leq & \esssup_{t \in [0,t_f]} l(|M(\xi_1(t) - \xi_2(t))|_2).\label{eq: intermediate step 6}
\end{align}
Finally, plugging \eqref{eq: intermediate step 5} and \eqref{eq: intermediate step 6} into \eqref{eq: intermediate step 10}, it follows that $|K(\xi_1)-K(\xi_2)| \rightarrow 0$ as $\sup_{t \in [0,t_f]}|\xi_1(t)-\xi_2(t)|\rightarrow 0$, which concludes the continuity of the mapping $K(\cdot)$. \qed 

\end{pf*}

\begin{customlem}{2.5}
Suppose $Z_{0,N} \in \mathcal{H}_N$. Let assumptions (A5)--(A7) hold and $\Pi_{N}(t)$ be defined as in {\cite[(15)]{Cheng2020optimalInReivew}}%\eqref{eq: mild solution of approximate Riccati}
. If assumption (A8) holds, then the mapping $K_N: C([0,t_f];\mathbb{R}^n) \rightarrow \mathbb{R}^+$ such that $K_N(\xi) \defeq \innerproduct{Z_{0,N}}{\Pi_{N}(0) Z_{0,N}}$ is continuous.
\end{customlem}

\begin{pf*}{Proof}
    Since the norm defined on $\mathcal{H}_N$ is inherited from that of $\mathcal{H}$, the proof follows from the same derivation of Lemma~\ref{lemma: continuity of equivalent cost wrt actuator location}.
\qed 
\end{pf*}

\revision{
\begin{customlem}{D.1}\label{lemma: continuity of the total cost wrt guidance}
Consider problem \cite[(P1)]{Cheng2020optimalInReivew} %\eqref{prob: new equivalent IOCA}
and its approximation \cite[(AP1)]{Cheng2020optimalInReivew}%\eqref{prob: equivalent finite dim approx integrated optimization problem}
. If assumptions (A4)--(A7) and (A9)--(A12) hold, then the following two implications hold:
\newline 1. For $p \in C([0,t_f];P)$, $\lim_{N \rightarrow \infty} |\costeval{{(AP1)}}{p} - \costeval{\text{(P1)}}{p}| = 0$, where $N$ is the dimension of approximation applied in \cite[(AP1)]{Cheng2020optimalInReivew}%\eqref{prob: finite dim approx integrated optimization problem}
.
\newline 2. The mapping $J_{(P1)}: C([0,t_f];P) \rightarrow \mathbb{R}^+$ is continuous, where $\costeval{(P1)}{p} = \innerproduct{\mathcal{Z}_0}{\Pi(0) \mathcal{Z}_0} + J_{\text{m}}(\xi,p)$. Here, the actuator state $\xi$ follows the dynamics \cite[(6)]{Cheng2020optimalInReivew}%\eqref{eq: general dynamics of the mobile actuator}
steered by the guidance $p$, and $\Pi(0)$ follows \cite[(12)]{Cheng2020optimalInReivew} % \eqref{eq: new operator Riccati equation} 
with the actuator state $\xi$.
\end{customlem}
}

\begin{pf*}{Proof}
    1. We first prove for $p \in C([0,t_f];P)$, 
    \begin{equation}\label{eq: intermediate step 1 for solution convergence}
    \lim_{N \rightarrow \infty} |\costeval{(AP1)%\eqref{prob: equivalent finite dim approx integrated optimization problem}
    }{p} - \costeval{(P1)
    %\eqref{prob: new equivalent IOCA}
    }{p}| = 0.
    \end{equation}
    The limit \eqref{eq: intermediate step 1 for solution convergence} is established for the following reason: By \cite[Theorem~2.4]{Cheng2020optimalInReivew}%\ref{thm: convergence of riccati operator}
    , $\norm{\Pi(0) - \Pi_{N}(0)}_{\mathcal{J}_q(\mathcal{H})} \rightarrow 0$ as $N \rightarrow \infty$. Since $\mathcal{J}_q(\mathcal{H}) \hookrightarrow \mathcal{L}(\mathcal{H})$ \cite{Burns2015Solutions}, it follows that 
    \begin{equation}\label{eq: intermediate step 2 for solution convergence}
        \lim_{N \rightarrow \infty} \norm{\Pi(0) - \Pi_{N}(0)}_{\mathcal{L}(\mathcal{H})} = 0.
    \end{equation}
    The convergence \eqref{eq: intermediate step 2 for solution convergence} and the fact
    \begin{equation}
        \lim_{N \rightarrow \infty} \norm{Z_{0,N} - \mathcal{Z}_0}_{\mathcal{H}} = 0
    \end{equation}
    imply that 
    \begin{equation}
        \lim_{N \rightarrow \infty} |\innerproduct{Z_{0,N}}{\Pi_{N}(0) Z_{0,N}} - \innerproduct{\mathcal{Z}_0}{\Pi(0) \mathcal{Z}_0}| = 0.
    \end{equation}
    The limit in \eqref{eq: intermediate step 1 for solution convergence} follows naturally since $\costeval{(AP1)%\eqref{prob: equivalent finite dim approx integrated optimization problem}
    }{p} - \costeval{(P1)%\eqref{prob: new equivalent IOCA}
    }{p} = \innerproduct{Z_{0,N}}{\Pi_{N}(0) Z_{0,N}} - \innerproduct{\mathcal{Z}_0}{\Pi(0) \mathcal{Z}_0}$.
    
    2. Consider the continuous mapping $T:L^2([0,t_f];\mathbb{R}^m) \rightarrow C([0,t_f];\mathbb{R}^n)$ as defined in the proof of \cite[Theorem~3.1]{Cheng2020optimalInReivew} %\ref{thm: existence of a solution of IOCA}
    such that $(Tp)(t) \defeq \xi(t) = e^{\alpha t} \xi_0 + \int_0^t e^{\alpha (t-\tau)} \beta p(\tau) \dd \tau $ for $t \in [0,t_f]$. Since the admissible guidance is in $C([0,t_f];\mathbb{R}^m) \subset L^2([0,t_f];\mathbb{R}^m)$, the continuity of $T$ is preserved, i.e., there exists $d_2 > 0$ such that for $p_1,p_2 \in C([0,t_f];\mathbb{R}^m)$
    \begin{align}
        \norm{Tp_1 - Tp_2}_{C([0,t_f];\mathbb{R}^n)} \leq  d_2 \norm{p_1 - p_2}_{C([0,t_f];\mathbb{R}^m)}. \label{eq: boundedness of trajectory norm}
    \end{align}
    Furthermore, define the mapping $\bar{J}_{\text{m}}: C([0,t_f];\mathbb{R}^m) \rightarrow \mathbb{R}^+$ by
\begin{equation}
    \bar{J}_{\text{m}}(p) \defeq J_{\text{m}}(Tp,p). \label{eq: define mobility cost as a mapping}
\end{equation}
    
    Define mappings $G: C([0,t_f];\mathbb{R}^m) \rightarrow \mathbb{R}^+$, $H: C([0,t_f];\mathbb{R}^n) \rightarrow \mathbb{R}^+$, and $H_f: C([0,t_f];\mathbb{R}^n) \rightarrow \mathbb{R}^+$ such that
    \begin{align}
        G(p) = & \   \int_0^{t_f} g(p(t),t) \dd t, \\
        H(p) = & \   \int_0^{t_f} h(Tp(t),t) \dd t, \\
        H_f(p) = & \  h_f(Tp(t_f)).
    \end{align}
    % where the mapping $T: C([0,t_f];\mathbb{R}^m) \rightarrow C([0,t_f];\mathbb{R}^n)$ is defined in \eqref{eq: define trajectory state as a mapping}.
    We first show the mapping $\bar{J}_{\text{m}}$ is continuous by showing that the mappings $G$, $H$, and $H_f$ are continuous since $\bar{J}_{\text{m}}(p) = G(p) + H(p) + H_f(p)$.

    Let $p_1,p_2 \in \mathcal{P}(p_{\max},a_{\max})$. Both the set of admissible guidance's values $P_0 \defeq \cup_{t \in [0,t_f]} \{p(t): p \in \mathcal{P}(p_{\max},a_{\max}) \}$ and the interval $[0,t_f]$ are closed and bounded (hence compact). Since $g: P_0 \times [0,t_f] \rightarrow \mathbb{R}^+$ is continuous, by the Heine-Cantor Theorem \cite[Proposition~5.8.2]{sutherland2009introduction}, $g$ is uniformly continuous, i.e., for all $\epsilon > 0$ there exists $\delta > 0$ such that for all $t \in [0,t_f]$
    \begin{equation}
        |p_1(t)-p_2(t)| < \delta \Rightarrow |g(p_1(t),t) - g(p_2(t),t)| < \epsilon.
    \end{equation}
    Hence, it follows that
    \begin{gather}
        \norm{p_1-p_2}_{C([0,t_f];\mathbb{R}^m)} =  \sup_{t \in [0,t_f]} |p_1(t)-p_2(t)| < \delta, \nonumber \\
        \Rightarrow |g(p_1(t),t) - g(p_2(t),t)| <  \epsilon, \quad \forall t \in [0,t_f].
    \end{gather}
    Therefore, for all $\epsilon>0$ there exists $\delta >0$ such that $\norm{p_1-p_2}_{C([0,t_f];\mathbb{R}^m)} < \delta $ implies
    \begin{equation}
          \int_0^{t_f} |g(p_1(t),t) - g(p_2(t),t)|  \dd t < \epsilon t_f,
    \end{equation}
    which concludes the continuity of the mapping $G$.

    Since the continuous image of a compact set is compact \cite[Proposition~5.5.1]{sutherland2009introduction}, the image set $T(\mathcal{P}(p_{\max},a_{\max}))$ is compact, i.e., the set $\Xi \defeq \{\xi \in C([0,t_f];\mathbb{R}^n): \xi = Tp, p \in \mathcal{P}(p_{\max},a_{\max}) \}$ is compact. The compactness of $\Xi$ implies that the set of actuator state's values $\xi(t)$, $\Xi_0 \defeq \cup_{t \in [0,t_f]} \{\xi(t)|\xi \in \Xi \}$, is closed. Furthermore, since $\norm{Tp}_{C([0,t_f];\mathbb{R}^n)}$ is bounded (see \eqref{eq: boundedness of trajectory norm}) and $\Xi_0$ is finite dimensional, the set $\Xi_0$ is compact. The compactness of $\Xi_0$ and continuity of the function $h: \Xi_0 \times [0,t_f] \rightarrow \mathbb{R}^+$ implies that $ h$ is uniformly continuous by the Heine-Cantor Theorem \cite[Proposition~5.8.2]{sutherland2009introduction}. Hence, for all $\epsilon > 0$ there exists $\delta > 0$ such that if $\norm{p_1-p_2}_{C([0,t_f];\mathbb{R}^m)} < \delta/d$, which implies $\norm{Tp_1-Tp_2}_{C([0,t_f];\mathbb{R}^n)} < \delta$, then
    \begin{equation}
           \int_0^{t_f} |h(Tp_1(t),t) - h(Tp_2(t),t)| \dd t \leq \epsilon t_f,
    \end{equation}
    which concludes the continuity of the mapping $H$.
    
    The mapping $H_f$ is continuous because for all $\epsilon > 0$ there exists $\delta >0$ such that if $\norm{p_1 - p_2}_{C([0,t_f];\mathbb{R}^m)} < \delta/d $, which implies $\sup_{t \in [0,t_f]} |Tp_1(t) - Tp_2(t)|< \delta$, then 
    \begin{equation}
        |Tp_1(t_f) - Tp_2(t_f)| < \delta,
    \end{equation}
    and, furthermore, $|H_f(p_1) - H_f(p_2)| = |h_f(Tp_1(t_f)) - h_f(Tp_2(t_f))| < \epsilon$ due to the continuity of $h_f$. Hence, we conclude the continuity of $\bar{J}_{\text{m}}$.
    
    The cost function of \cite[(P1)]{Cheng2020optimalInReivew} %\eqref{prob: new equivalent IOCA} 
    is the sum of two parts: the PDE cost $\innerproduct{\mathcal{Z}_0}{\Pi(0) \mathcal{Z}_0}$, cast as a continuous mapping $K: C([0,t_f];\mathbb{R}^n) \rightarrow \mathbb{R}^+$ (see Lemma~\ref{lemma: continuity of equivalent cost wrt actuator location}) and the mobility cost $J_{\text{m}}(\xi,p)$, cast as a continuous mapping $\bar{J}_{\text{m}}: C([0,t_f];\mathbb{R}^m) \rightarrow \mathbb{R}^+$. 
    Due to the continuity of the mapping $T$ (see \eqref{eq: boundedness of trajectory norm}), there exists $d_3>0$ such that
    \begin{equation}\label{eq: eq1}
        |K(Tp_1) - K(Tp_2)|_1 \leq d_3 \norm{p_1 - p_2}_{C([0,t_f];\mathbb{R}^m)}.
    \end{equation}
    The continuity of the mapping $J_{(P1)}$ follows from \eqref{eq: eq1} and the continuity of $\bar{J}_{\text{m}}$. \qed
\end{pf*}

\bibliographystyle{abbrv}
\bibliography{refSup}